\def\block{\hbox{${\vcenter{\vbox{\hrule height 0.4pt\hbox{\vrule width 0.4pt 
height 6pt \kern 5pt\vrule width 0.4pt}\hrule height 0.4pt }}}$}}
\def\qed{\hfill\block\medskip}
\font\eightit=cmti8

\headline={\ifnum\pageno>1\eightit\hfill spherical 2-categories and 4-manifold invariants\hfill\fi}

\def\joinrel{\mathrel{\mkern-3mu}}
\def\Relbar{\mathrel=}
\def\longeq{\Relbar\joinrel\Relbar}

\def\Not{Notation  \the\chno.\the\thmno}
\def\Def{Definition \the\chno.\the\thmno}
\def\Lem{Lemma \the\chno.\the\thmno}
\def\Thm{Theorem \the\chno.\the\thmno}
\def\Cor{Corollary \the\chno..\the\thmno}
\def\Ex{Example \the\chno.\the\thmno}

\def\qd{\hbox{dim}_q}
\def\rk{\hbox{rk}}
\def\tr{\hbox{tr}}

\def\Tr{\hbox{Tr}}
\def\id{\hbox{id}}

\def\ev{\hbox{ev}}

\def\Hom{\hbox{Hom}}
\def\THom{\hbox{2Hom}}

\def\End{\hbox{End}}
\def\TEnd{\hbox{2End}}
\def\Hilb{\hbox{Hilb}}
\def\hom{\hbox{hom}}
\def\Vect{\hbox{Vect}}
\def\TVect{\hbox{2Vect}}
\def\THilb{\hbox{2Hilb}}
\voffset=2\baselineskip

\input epsf
\magnification 1200

\newcount\chno
\newcount\thmno
\newcount\equano
\newcount\refno

\chno=0
\thmno=0
\equano=0
\refno=0 

\centerline{\bf Spherical 2-categories and 4-manifold invariants}
\bigskip
\centerline{\sl Marco Mackaay}
\medskip
\centerline{Area Departamental de Matematica}
\centerline{Universidade do Algarve, UCEH}
\centerline{Campus de Gambelas}
\centerline{8000 Faro}
\centerline{Portugal}
\centerline{e-mail: mmackaay@ualg.pt}
\bigskip
\beginsection Introduction\par
In this paper I give a formal construction of a 4-manifold invariant 
out of what I call a {\it non-degenerate finitely semi-simple 
semi-strict spherical 2-category of non-zero dimension}. 

For some time now people have 
had the feeling that 4-manifold invariants and certain kinds of monoidal 
2-categories have a relation with each other similar to that of 3-manifold 
invariants and certain kinds of monoidal categories. 

The first evidence for this feeling can be found in the work of 
Crane and collaborators. In [16] Crane  
and Frenkel give a formal construction of 4-manifold invariants out 
of Hopf categories and indicate where one should look for such 
algebraic objects, namely in the work on crystal bases by Saito and 
the work on canonical bases and perverse sheaves by Lusztig. 
One could argue that it should be possible 
to use the 2-category of representations of a Hopf category instead 
of the Hopf category itself; the reason for this thought being that for the construction of 3-manifolds one 
can use Hopf algebras, as Kuperberg [29] and Chung, Fukuma and Shapere 
did [15], or the category of representations 
of Hopf algebras as people like Turaev and Viro [37], 
Yetter [40], and Barrett and 
Westbury [9] did. Recently Neuchl [31] showed that the representations of 
a Hopf category do form a monoidal $2$-category indeed.

In [18] Crane and Yetter proposed a construction of $4$-manifold 
invariants out of the semi-simple sub-quotient of the category of 
finite dimensional representations of the quantum group $U_q(sl(2))$ 
for $q$ a principal $4r^{th}$ root of unity. In [17] Crane, Kauffman and 
Yetter generalized this construction for any finitely semi-simple tortile 
category and gave detailed proofs. These Crane-Yetter invariants 
can be seen as a special case in which the authors use a 
2-category of a certain type with only one object.  

Another  
piece of evidence for the aforementioned ``feeling'' is the work of 
Baez and collaborators.
In a series of 
papers [1,2,3,4] they have tried to convince people of the importance of 
$n$-categories. Among other important applications $n$-categories seem to 
be the right context in which to study Topological Quantum Field Theories, 
which in this context can be defined as $n$-functors 
from the $n$-category of $n$-cobordisms to the $n$-category of 
$n$-vector spaces. If we take a $3$-category with one object 
we get a monoidal $2$-category, which is roughly speaking the case we 
are concerned with in this paper. 

Based on 
the work of Carter, Rieger and Saito [13,14], several people [5,6,22,28,30] 
have shown that 
braided 2-categories with duals form the right algebraic context in which to 
study invariants of 2-tangles. This should be closely related 
to the construction of 4-manifold invariants, at least formally. 

So altogether one could say that there is more than enough reason 
to believe that it is possible to construct $4$-manifold invariants 
out of certain kinds of $2$-categories. But nowhere in the literature can 
one find a paper with an explicit 
construction. This paper is meant to fill this gap.

In this paper I use triangulations of 4-manifolds   
for the construction of a state sum. I will show that this state sum 
is independent of the chosen triangulation by using Pachner's 
theorem [32], that relates triangulations of piece-wise linear homeomorphic 
manifolds. The whole construction should be  
considered as a lift to the fourth dimension of Barrett's and Westbury's [9] 
construction of 3-manifold  
invariants out of non-degenerate finitely semi-simple spherical categories. 
For that reason I have given the kind of 2-categories that 
I use the name above. 

The simplest example of a $2$-category as used in this paper is the 
semi-strict version of 2Hilb, the $2$-category of 2-Hilbert spaces. In [3] 
Baez defined the weak version of 2Hilb. In this paper I define the 
completely coordinatized version, $\THilb_{\rm cc}$, which is semi-strict.  

In section 6 I show how 
any finite group gives rise to 
a $2$-category of the right kind and I give an explicit formula for the 
invariant. This invariant looks like a four dimensional version of the 
Dijkgraaf-Witten invariant [21]. 

It is likely that 
the 2-category of representations of the right kind of Hopf category 
will be such a 2-category too. 
But we have only one Hopf category that has been worked out in detail, namely 
the categorification $C(D(G))$ of the quantum double of a finite 
group [19]. Probably $C(D(G))$-mod, the $2$-category of finite dimensional 
representations of $C(D(G))$, is an example of a non-degenerate 
finitely semi-simple semi-strict spherical $2$-category. But I have not 
worked out the details yet. 
Carter, Kauffman and Saito are working out the Crane-Frenkel invariant for this particular example [12]. 
In [16] Crane and Frenkel indicate that it is possible to construct more examples of 
Hopf categories $C(U_0(g))$ using the crystal bases of quantum groups at $q=0$. 
However, they do not define the categorification of the antipode for a 
general Hopf category. In the aforementioned case of $C(D(G))$ this antipode 
looks to be straight forward. In the case of $C(U_0(g))$ it is not known how to define the 
right categorification of the antipode. 

Also Carter, Kauffman and Saito found that it is necessary to impose extra conditions on the 
cocycles defining the structure isomorphisms of the Hopf category $C(D(G))$ in order to 
obtain invariance of the state sum under permutation of the vertices of a chosen triangulation 
of the 4-manifold. This was not foreseen in [16] and, as far as I know, it is not known 
how to obtain invariance under permutation of vertices for an arbitrary Hopf 
category.

In a future paper Crane and Yetter show 
how to build a monoidal 2-category out of the modules of a quantum group 
at $q=0$ 
using 
their crystal bases. Here they avoid the Hopf categories and build 
the $2$-categories directly, which makes a construction of $4$-manifold 
invariants out of a certain kind of $2$-categories, as presented in 
this paper, even more desirable. It is definitely a good 
place to look for interesting examples of monoidal $2$-categories,  
but it is not likely that these 2-categories  
are already the ones we are looking for. It is like having the 
non-finitely semi-simple ``trivially'' spherical category $U(g)$-mod, 
where $g$ is a finite dimensional semi-simple Lie algebra, and having 
an abstract Turaev-Viro-like construction of $3$-manifold invariants 
which requires finitely semi-simple spherical categories. 
In that case the missing link comes from the deformation theory of $U(g)$ 
which shows that there are deformations $U_q(g)$, 
where $q$ is a certain 
root of unity, such that there is a certain non-degenerate quotient of the category 
of {\sl tilting modules} of $U_q(g)$ that is a finitely semi-simple 
spherical category. For details about this see [9] and some references 
therein. Likewise, in the case of $2$-categories, one should  
first define the deformation theory of monoidal $2$-categories analogously 
to what Crane and Yetter have done for monoidal categories [20]. Then 
one has to find actual deformations of the aforementioned 2-categories 
and finally one has to ``melt'' them, i.e. get back to generic $q$.  
These last two problems are of a very deep nature and certainly far beyond  
the scope of this paper. In the meanwhile it is worthwile, I think, to study 
the kind of 2-category we are looking for from an abstract point of 
view. 

\advance\chno by 1
\equano=0
\thmno=0
\beginsection\the\chno. The basic idea\par
Throughout this paper a manifold means an oriented piece-wise linear 
compact 4-manifold without boundary. A triangulated manifold $(M,T)$ is a manifold $M$ 
together with a given simplicial 4-complex $T$, the triangulation, such that its 
underlying PL-manifold is PL-homeomorphic 
to $M$. Throughout this paper we will always assume that there is 
a total ordering on the vertices of the triangulation 
of a manifold. 
\noindent A combinatorial isomorphism between two triangulated manifolds $(M,T)$ 
and $(M',T')$ will always mean an 
isomorphism between the simplicial complexes $T$ and $T'$. A simplicial isomorphism 
between two triangulated manifolds $(M,T)$ and 
$(M',T')$ will always mean a combinatorial isomorphism which preserves the ordering 
on the vertices.  
I also want to fix some   
notation. The letter ${\bf F}$ will always denote a fixed field of 
characteristic 0 and any vector  
space in this paper will be a finite dimensional vector space over 
${\bf F}$. My notation for the simplices follows Barrett's and Westbury's 
convention in [9]. The standard $n$-simplex $(012\ldots n)$ with 
vertices $\{0,1,2,\ldots, n\}$ has the standard orientation $(+)$. 
The opposite orientation is denoted by $(-)$. The standard $4$-simplex 
$+(01234)$ has boundary 
$$(1234)-(0234)+(0134)-(0124)+(0123).$$ The $4$-simplex $-(01234)$ has 
the same boundary but with the opposite signs. The sign with which 
a tetrahedron appears in the boundary of a $4$-simplex I will call 
the induced orientation of the tetrahedron. The total ordering on 
the vertices of the simplicial complex defining the triangulation 
of a manifold induces an ordering on all $4$-simplices and 
the ordering on the vertices of a $4$-simplex induces an orientation 
on its underlying polytope. If the orientation of the underlying 
polytope of a $4$-simplex induced 
by the orientation of the manifold is equal to the orientation 
of this polytope induced by the total ordering on the vertices of 
the $4$-simplex, then our convention will be that the $4$-simplex 
has the positive orientation, as a simplex, and if the two 
induced orientations of the polytope are opposite, 
we take the $4$-simplex to be negatively oriented, as a simplex.

\advance\thmno by 1
\proclaim Remark~\the\chno.\the\thmno. Notice that in a triangulated 
manifold   
each tetrahedron lies in the boundary of exactly two 
$4$-simplices, occuring once with induced orientation $(+)$ and once 
with $(-)$.\par 

Let $(M,T)$ be a triangulated manifold. For the definition of my state sum I need two 
sets of labels, $\cal E$ and $\cal F$ respectively. The edge $(ij)$ 
with vertices $i,j$ is labelled with $e_{ij}\in {\cal E}$, 
the face $(ijk)$ with vertices $i,j,k$ is 
labelled with $f_{ijk}\in {\cal F}$. 
Let $T((0123),e,f)$ be the labelled standard oriented 
tetrahedron $+(0123)$ . We do not take the different orientations of 
the edges and the faces into account for the labelling yet. 
\advance\thmno by 1
\proclaim\Not. The state space of 
this labelled tetrahedron is 
a certain vector space $$2H((0123),e,f).$$ 
The state space of $-(0123)$ is defined to be the dual vector 
space $$2H((0123),e,f)^*.$$\par

\advance\thmno by 1
\proclaim\Not. Let $+(01234)$ be the standard oriented $4$-simplex with
the triangle $(ijk)$ labelled by $f_{ijk}\in {\cal F}$ and the edge 
$(ij)$ by 
$e_{ij}\in {\cal E}$. Denote the state space of the labelled tetrahedron $(ijkl)$ by $2H(ijkl)$ 
just as a shorthand. The partition function of this labelled $4$-simplex 
is a certain linear map
$$Z(+01234)\colon 2H(0234)\otimes 2H(0124)\to  2H(1234)\otimes 2H(0134)\otimes 
2H(0123).$$
The partition function of $-(01234)$ is also a certain linear map 
$$Z(-01234)\colon 2H(1234)\otimes 2H(0134)\otimes 
2H(0123)\to 2H(0234)\otimes 2H(0124).$$\par

\noindent Notice that $Z(\pm(01234))$ is defined for a fixed labelling of $T$, 
although this dependence does not show up in the notation. This is a deliberate 
choice, or a deliberate flaw in the notation, that I want to allow myself 
in order to write down formulas that are not too polluted by a high number of sub- and superscripts. I think that the context will leave no doubt of 
what depends on what in my formulas. 

Assume that $M=(M,T)$ is labelled with a fixed labelling $\ell$. 
Notice that the 
ordering on the vertices of $T$ induces a natural total ordering on the 
tetrahedra, by 
means of the boundary operator, within each 
$4$-simplex. For example in the ordered $4$-simplex $(abcde)$ the ordering is 
given by $1.(bcde),$ $2.(acde),$ $3.(abde),$ $4.(abce)$, $5.(abcd)$, which is 
independent of the orientation. In the same way we get a fixed ordering 
on the triangles within each tetrahedron etc. Fix also a total ordering 
on the $4$-simplices, for example the one induced by the total ordering 
on the vertices of the whole triangulation. Take out of each $4$-simplex 
the tetrahedra that 
appear with a negative sign in its boundary in the induced order 
described above. Together with the chosen ordering on the $4$-simplices  
this fixes an ordering of all the tetrahedra of $M$. Notice that each  
tetrahedron appears exactly once in this way by remark 1.1.  
\advance\thmno by 1
\proclaim\Def. Let $V(M,T,\ell)$ be the tensor product of the 
state spaces of all 
the tetrahedra of $T$ in which the ordering of the factors is as 
described above. The tensor 
product of the respective partition functions applied to $V(M,T,\ell)$ has 
its image in a permuted tensor product $V(M,T,\ell)'$. Again remark 1.1 is 
vital. 
Now compose this linear map with 
the linear map $P(M,T,\ell)\colon V(M,T,\ell)'\to V(M,T,\ell)$ induced 
by the standard transposition $P\colon x\otimes y\to y\otimes x$. The result is a linear map 
$L(M,T,\ell)\colon 
V(M,T,\ell)\to V(M,T,\ell)$. The element $Z(M,T,\ell)\in {\bf F}$ is 
defined to be the 
trace of $L(M,T,\ell)$.\par
\noindent Notice that $Z(M,T,\ell)$ does not depend on the ordering 
on the $4$-simplices, because it is defined by a conjugation invariant trace. 
In the next section we will prove that 
$Z(M,T,\ell)$ is a combinatorial invariant of $M$.

The state sum $I(M,T)$ is a certain weighted sum over all 
labellings of the numbers $Z(M,T,\ell)$.

So far we have only sketched the basic idea of our approach without 
telling anyone where to get these state spaces and these partition 
functions. In section 3 we will show how they appear 
naturally out of a certain kind of $2$-categories. Therefore 
we have to study this kind of $2$-categories in the next section 
first. 

\advance\chno by 1
\equano=0
\thmno=0
\beginsection\the\chno. Spherical 2-categories\par
In this paragraph I define what I call a spherical 2-category. The 
underlying 2-category will always be assumed to be strict. This means 
that the composition is strictly associative and the composite of a 
1-morphism with an identity 1-morphism is equal to the 1-morphism itself.
I will denote the composite of two $1$-morphisms $f,g$ by $fg$, the vertical composite of two $2$-morphisms $\alpha,\beta$ 
by $\alpha\cdot
\beta$ and their horizontal composite by $\alpha\circ\beta$. If $f,f'\colon A\to 
B$ and $g,g'\colon B\to C$ are $1$-morphisms and $\alpha\colon f\to f'$ 
and $\beta\colon g\to g'$ are $2$-morphisms, then $f\circ\beta\colon fg\to fg'$ denotes  
$1_f\circ\beta$ and $\alpha\circ g\colon fg\to f'g$ denotes $\alpha\circ 
 1_g$. In this paper $\Hom(A,B)$ will always denote the category whose 
objects are all $1$-morphisms 
with source $A$ and target $B$ and whose morphisms are all 2-morphisms between such $1$-morphisms, where the composition of the morphisms is defined by the 
vertical composition of the $2$-morphisms. $\THom(f,g)$ will always denote the 
set of $2$-morphisms with source $f$ and target $g$. When the source and 
target are equal we will also use the notations $\End(A)=\Hom(A,A)$ 
and $\TEnd(f)=\THom(f,f)$.   

So let $C$ be a strict 2-category. This is not too restrictive an assumption  
because a weak 2-category can always be strictified, see [27]. We also assume that $C$ has a semi-strict 
monoidal structure. Loosely speaking this means that for every pair 
of objects $A,B$ in $C$ there is a unique object $A\otimes B$. For 
every object $A$ and every $1$-morphism $f\colon X\to Y$ there is a unique 
$1$-morphism $A\otimes f\colon A\otimes X\to A\otimes Y$ and a 
unique 1-morphism $f\otimes A\colon X\otimes A\to Y\otimes A$. 
For every 
object $A$ and every $2$-morphism $\alpha\colon f\Rightarrow g$ there is 
a unique $2$-morphism $A\otimes\alpha\colon A\otimes f\Rightarrow A\otimes g$ 
and a unique $2$-morphism $\alpha\otimes A\colon f\otimes A\Rightarrow 
g\otimes 
A$. Also 
there is an identity object $I$ such that $I\otimes X=X\otimes I=X$ 
for all objects, $1$- and $2$-morphisms $X$. All 
the usual structural 2-isomorphisms are identities except one: given 
a pair of $1$-morphisms $f\colon A\to C,\ g\colon B\to D$ in $C$ there 
is a $2$-isomorphism called the {\sl tensorator} 
$$\otimes_{f,g}\colon (f\otimes B)(C\otimes g)\Rightarrow (A\otimes g)(f\otimes D).$$ 
This $2$-isomorphism is required to satisfy some conditions. These 
conditions guarantee that $\otimes_{f,g}$ behaves well under the tensor 
product and composition. The obvious condition which tells us how to 
obtain $\otimes_{fg,h}$ by pasting $\otimes_{f,h}$ and $\otimes_{g,h}$,  
and analogously how to obtain $\otimes_{f,gh}$ by pasting $\otimes_{f,g}$ and 
$\otimes_{f,h}$, resembles the condition defining a braiding in a monoidal 
category. For the exact definition 
of a semi-strict monoidal structure see [27]. 
\advance\thmno by 1
\proclaim\Def. A {\rm semi-strict monoidal $2$-category} is a strict 
$2$-category with a semi-strict monoidal structure.\par

\noindent Let $C$ always be a semi-strict monoidal $2$-category. Again 
this is a legitimate assumption, since every weak monoidal $2$-category 
is equivalent in a well defined sense to a semi-strict one [26]. 
The following 
lemma is well known. For a proof see [27].
\advance\thmno by 1
\proclaim\Lem. Let $I$ be the identity object in $C$. Then the 
category ${\rm End}(I)$ is a braided monoidal category with the tensor 
product being defined by 
the horizontal composition of $1$- and $2$-morphisms and the braiding 
being defined by the tensorator $\otimes_{\bullet,\bullet}$.\par 

In order to get to the spherical condition I first have to define duality  
in $C$. For this I will copy the definition given by Baez and Langford in 
[5,6,30]. 
\advance\thmno by 1
\proclaim\Def. $C$ is called a {\rm semi-strict monoidal $2$-category 
with duals} if it 
is equipped with the following structures:
\item{1.}For every $2$-morphism  $\alpha\colon f\Rightarrow g$ there is a 
$2$-morphism $\alpha^*\colon g\Rightarrow f$ called the {\it dual} of $\alpha$.
\item{2.} For every $1$-morphism $f\colon A\to B$ there is a $1$-morphism $f^*\colon
B\to A$ called the {\it dual} of $f$, and $2$-morphisms $i_f\colon 
1_A\Rightarrow ff^*$  and $e_f\colon f^*f\Rightarrow 1_B$, called  the {\it unit} and 
{\it counit} of $f$, respectively.
\item{3.} For any object $A$, there is an object $A^*$ called the 
{\it dual} of $A$, $1$-morphisms $i_A\colon I\to A\otimes A^*$ and 
$e_A\colon A^*\otimes A\to I$ called the {\it unit} and {\it counit} 
of $A$, respectively, and a $2$-morphism $T_A\colon (i_A\otimes  A)(A\otimes e_A)\Rightarrow 1_A$ called the {\it triangulator} of $A$.

\noindent We say that a $2$-morphism $\alpha$ is {\it unitary} if it 
is invertible and $\alpha^{-1}=\alpha^*$. Given a $2$-morphism $\alpha
\colon f\Rightarrow g$, we define the {\it adjoint} $\alpha^{\dag}
\colon g^*\Rightarrow f^*$ by 
$$\alpha^{\dag}=(g^*\circ i_f)\cdot(g^*\circ\alpha\circ f^*)\cdot(e_g\circ f^*).$$
\noindent In addition, the structures above are also required to 
satisfy the following conditions:
\item{1.} $X^{**}=X$ for any object, $1$-morphism or $2$-morphism.
\item{2.} $1^*_{X}=1_X$ for any object or $1$-morphism $X$.
\item{3.} For all objects $A,B$, $1$-morphisms $f,g$, and $2$-morphisms 
$\alpha,\beta$ for which both sides of the following  equations are 
well-defined, we have
$$(\alpha\cdot\beta)^*=\beta^*\cdot\alpha^*,$$
$$(\alpha\circ\beta)^*=\alpha^*\circ\beta^*,$$
$$(fg)^*=g^*f^*,$$
$$(A\otimes\alpha)^*=A\otimes\alpha^*,$$
$$(A\otimes f)^*=A\otimes f^*,$$
$$(\alpha\otimes A)^*=\alpha^*\otimes A,$$
$$(f\otimes A)^*=f^*\otimes A,$$
$$(A\otimes B)^*=B^*\otimes A^*.$$

\item{4.} For all $1$-morphisms $f,g$ the $2$-isomorphism $\otimes_{f,g}$ 
is unitary.
\item{5.} For any object or $1$-morphism $X$ we have $i_{X^*}=e^*_{X}$ 
and $e_{X^*}=i^*_{X}$.
\item{6.} For any object $A$, the $2$-morphism $T_A$ is unitary.
\item{7.} For any objects $A$ and $B$ we have  
$$i_{A\otimes B}=i_A\cdot(A\otimes i_B\otimes A^*),$$
$$e_{A\otimes B}=(B^*\otimes e_A\otimes B)\cdot e_B,$$
$$T_{A\otimes B}=[(i_A\otimes A\otimes B)(A\otimes \bigotimes\nolimits^{-1}_{i_{B},e_{A}}\otimes B)(A\otimes B\otimes e_B)]\cdot [(T_A\otimes B)\circ(A\otimes T_B)].$$
\item{8.} $T_I=1_{1_I}$.
\item{9.} For any object $A$ and $1$-morphism $f$, we have
$$i_{A\otimes f}=A\otimes i_f,$$
$$i_{f\otimes A}=i_f\otimes A,$$
$$e_{A\otimes f}=A\otimes e_f,$$
$$e_{f\otimes A}=e_f\otimes A.$$
\item{10.} For any $1$-morphisms $f,g$, $i_{fg}=i_f\cdot(f\circ i_g\circ 
f^*)$ and 
$e_{fg}=(g^*\circ e_f\circ g)\cdot e_g$.
\item{11.} For any $1$-morphism $f$, $(i_f\circ f)\cdot (f\circ e_f)=1_f$ and 
$(f^*\circ i_f)\cdot (e_f\circ f^*)=1_{f^*}$.
\item{12.} For any $2$-morphism $\alpha$, $\alpha^{\dag *}=\alpha^{*\dag}$.
\item{13.} For any object $A$ we have
$$[i_A\circ (A\otimes T^{\dag}_{A^*})]\cdot [\bigotimes\nolimits^{-1}_{i_A,i_A}(A\otimes 
e_A\otimes A^*)]\cdot [i_A\circ (T_A\otimes A^*)]=1_{i_{A}}.$$
\par

\noindent Note that the first two identities in condition 7 imply 
$i_I=e_I=1_I$ so that condition 8 makes sense; the source of $T_I$ 
is equal to $1_I$. It would be worthwhile to study weaker notions of duality 
and prove a coherence theorem that allows one to strictify $2$-categories 
with such duality up to a semi-strict monoidal $2$-category with duals 
as defined above. In [7] Barrett and Westbury prove such a coherence 
theorem for monoidal categories with duals. In this paper we always work 
with semi-strict monoidal 2-categories with duals as defined above. 

For our purpose we   
need a little bit of extra structure. Let $f\colon A\to B$ be a $1$-morphism 
in $C$. We define the $1$-morphism ${^{\#}f}\colon B^*\to A^*$ by 
$$B^*\buildrel i_{A^*}\otimes B^*\over\longrightarrow A^*\otimes A
\otimes B^*\buildrel A^*\otimes f\otimes B^*\over\longrightarrow 
A^*\otimes B\otimes B^*\buildrel A^*\otimes e_{B^*}\over\longrightarrow
A^*.$$
Analogously we define $f^{\#}\colon B^*\to A^*$ by 
$$B^*\buildrel B^*\otimes i_A\over\longrightarrow B^*\otimes A\otimes 
A^*\buildrel B^*\otimes f\otimes A^*\over\longrightarrow 
B^*\otimes B\otimes A^*\buildrel e_B\otimes A^*\over\longrightarrow
A^*.$$
Notice that given a $2$-morphism $\alpha\colon f\Rightarrow g$ in $C$ 
we also have the $2$-morphism $\alpha^{\#}\colon f^{\#}\Rightarrow g^{\#}$ 
defined by 
$${^{\#}\alpha}=1_{(i_{A^*}\otimes B^*)}\circ(A^*\otimes\alpha\otimes B^*)
\circ 1_{(A^*\otimes e_{B^*})},$$ 
\noindent and the 2-morphism ${^{\#}\alpha}\colon {^{\#}f}\Rightarrow {^{\#}g}$ defined by 
$$\alpha^{\#}=1_{(B^*\otimes i_A)}\circ(B^*\otimes\alpha\otimes A^*)\circ 
1_{(e_B\otimes A^*)}.$$
\noindent It is not difficult to show that for any 1-morphisms 
$f\colon A\to B$ and $g\colon B\to C$ 
there exist unitary $2$-morphisms 
$${^{\#}(fg)}\cong {^{\#}g}{^{\#}f}$$ 
\noindent and 
$$(fg)^{\#}\cong g^{\#}f^{\#}.$$
\noindent The first 2-isomorphism is given by 
$$(1_{(i_{A^*}\otimes C^*)(A^*\otimes f\otimes C^*)}\circ (A^*\otimes 
T_B^{\dag}\otimes C^*)\circ 1_{(A^*\otimes g\otimes C^*)(A^*\otimes e_{C^*})})
$$
$$\cdot (1_{(i_{A^*}\otimes C^*)(A^*\otimes f\otimes C^*)(A^*\otimes B\otimes 
i_{B^*}\otimes C^*)}\circ(A^*\otimes\bigotimes\nolimits_{e_{B^*},g}\otimes C^*)\circ 1_{(A^*\otimes e_{C^*})})$$
$$\cdot (\bigotimes\nolimits_{(i_{A^*}\otimes C^*)(A^*\otimes f\otimes C^*)
,(i_{B^*}\otimes C^*)(B^*\otimes g\otimes C^*)}\circ 1_{(A^*\otimes e_{B^*}
\otimes C\otimes C^*)(A^*\otimes e_{C^*})})$$
$$\cdot 
(1_{(i_{B^*}\otimes C^*)(B^*\otimes g\otimes C^*)}\circ \bigotimes\nolimits
_{(i_{A^*}\otimes B^*)(A^*\otimes f\otimes B^*)(A^*\otimes e_{B^*}),e_{C^*}}).$$ 
\noindent The second 2-isomorphism is given by 
$$(1_{(C^*\otimes i_A)(C^*\otimes f\otimes A^*)}\circ (C^*\otimes T_B^{-1}
\otimes A^*)\circ 1_{(C^*\otimes g\otimes A^*)(e_C\otimes A^*)})$$
$$\cdot (1_{(C^*\otimes i_A)(C^*\otimes f\otimes A^*)(C^*\otimes i_B\otimes 
B\otimes A^*)}\circ (C^*\otimes \bigotimes\nolimits_{g,e_B}^{-1}\otimes A^*)
\circ 1_{(e_C\otimes A^*)})$$
$$\cdot (\bigotimes\nolimits_{(C^*\otimes i_B)(C^*\otimes g\otimes B^*),
(C^*\otimes i_A)(c^*\otimes f\otimes A^*)}^{-1}\circ 1_{(C^*\otimes C\otimes 
e_B\otimes A^*)(e_C\otimes A^*)})$$
$$\cdot (1_{(C^*\otimes i_B)(C^*\otimes g\otimes B^*)}\circ 
\bigotimes\nolimits_{e_C,(B^*\otimes i_A)(B^*\otimes f\otimes A^*)(e_B\otimes 
A^*)}^{-1}).$$

\noindent If $A=B=C=I$ we have ${^{\#}f}=f^{\#}=f$ and ${^{\#}g}=g^{\#}=g$ 
and the $2$-isomorphisms above become simply equal 
to $\otimes_{f,g}$ and $\otimes_{g,f}^{-1}$ respectively, since $i_I=e_I=1_I$, $T_I=1_{1_I}$ and $\otimes_{X,1_I}=\otimes_{1_I,X}=
1_X$ for any 1-morphism $X$. 
 
The following definition reminds us  
of the definition of a pivotal category (for the definition of a pivotal 
category see [24]).
\advance\thmno by 1
\proclaim\Def. A {\rm semi-strict pivotal $2$-category} is a semi-strict 
monoidal $2$-category with duals $C$ such that for any $1$-morphism  
$f\colon A\to B$ there exists a unitary $2$-morphism
$$\phi_f\colon {^{\#}f}\Rightarrow f^{\#}$$
\noindent satisfying the following conditions:
\item{1)}For any 2-morphism $\alpha\colon f\Rightarrow g$, we have 
$${^{\#}\alpha}\cdot\phi_g=\phi_f\cdot {\alpha^{\#}},$$
\item{2)}$\phi_{f}^{\dag}=\phi_{f^*},$
\item{3)}The following diagram commutes
$$\vbox{\offinterlineskip
\halign{\hfil#\hfil&\quad\hfil#\hfil\quad&\hfil#\hfil\cr
${^{\#}(fg)}$&$\buildrel \phi_{fg}\over{\hbox to 30pt
{\rightarrowfill}}$&$(fg)^{\#}$\cr
\noalign{\vskip5pt}
$\Bigg\downarrow$&&$\Bigg\downarrow$\cr
\noalign{\vskip5pt}
${^{\#}g}{^{\#}f}$&$\buildrel \phi_{g}\circ\phi_{f}\over{\hbox to 30pt{
\rightarrowfill}}$&$g^{\#}f^{\#}$
\cr}}$$ 
\item{4)}For any $2$-morphism $\alpha$, we have 
$${^{\#}(\alpha^{\dag})}={^{\dag}(\alpha^{\#})}\ \ \hbox{and}\ \ 
({^{\#}\alpha})^{\dag}=({^{\dag}\alpha})^{\#}.$$\par

\noindent Note that, by definition, we have 
$${^{\#}(f^*)}=(f^{\#})^*\ \ \hbox{and}\ \ ({^{\#}f})^*=(f^*)^{\#}$$
for any $1$-morphism $f$, so these conditions make sense. Note also that 
for $A=B=C=I$ the third condition becomes 
$$\phi_{fg}=\otimes_{f,g}(\phi_g\circ\phi_f)\otimes_{g,f}.$$
\noindent Together with the second condition this reminds us of the conditions 
defining a {\sl twist} in a ribbon category (see [35]). As a matter of fact 
lemma 
2.14 shows that this is 
exactly what we should have in mind. 
Note finally that, by definition, any $2$-morphism $\alpha$ satisfies 
$$(\alpha^{\dag})^*={^{\dag}(\alpha^*)}\ \ \hbox{and}\ \ 
({^{\dag}\alpha})^*=(\alpha^*)^{\dag},$$
\noindent where $${^{\dag}\alpha}=(i_{f^*}\circ g^*)\cdot(f^*\circ\alpha\circ g^*)
\cdot(f^*\circ e_{g^*}).$$
\noindent This, together with conditions 1 and 12 in the definition 
of duality (def.~2.3), implies the equalities 
$${^{\dag}\alpha}={^{\dag}(\alpha^{**})}=((\alpha^*)^{\dag})^*=
(\alpha^{\dag})^{**}=\alpha^{\dag}.$$
\noindent Condition 11 in the same definition implies 
$${(^{\dag}\alpha)^{\dag}}={^{\dag}(\alpha^{\dag})}=\alpha,$$
\noindent so we also get 
$${^{\dag\dag}\alpha}=\alpha^{\dag\dag}={(^{\dag}\alpha)^{\dag}}=\alpha.$$ 

The next thing to define is the notion of a {\it trace-functor} in a semi-strict  
pivotal $2$-category. Let $C$ be such a $2$-category for the rest of 
this section. 
\advance\thmno by 1
\proclaim\Def. For any object $A$ in $C$ there is a {\rm left trace 
functor} $\Tr_L\colon \End(A)\to\End(I)$. For any object $f\in\End(A)$ 
define $\Tr_L(f)$ by 
$$I\buildrel i_{A^*}\over\longrightarrow A^*\otimes A\buildrel A\otimes f\over
\longrightarrow A^*\otimes A\buildrel e_A\over\longrightarrow I.$$
For a morphism $\alpha\colon f\Rightarrow g$ in $\End(A)$ define $\Tr_L(\alpha)$ by
$$1_{i_{A^*}}\circ(A^*\otimes\alpha)\circ 1_{e_A}.$$ 
Analogously there is a {\rm right trace functor} $\Tr_R$ defined by 
$$I\buildrel i_{A}\over\longrightarrow A\otimes A^*\buildrel f\otimes A\over
\longrightarrow A\otimes A^*\buildrel e_{A^*}\over\longrightarrow I,$$
\noindent and
$$1_{i_{A}}\circ(\alpha\otimes A^*)\circ 1_{e_{A^*}}.$$\par

\noindent The following lemma
 is the analogue of the well known 
lemma that says that traces are conjugation invariant in pivotal 
categories.

\advance\thmno by 1
\proclaim\Lem. For any two objects $A,B$ in $C$ and any two 
$1$-morphisms 
$f\colon A\to B$ and $g\colon B\to A$ we have a unitary $2$-morphism 
$\theta_{f,g}^R\colon \Tr_R(fg)\Rightarrow \Tr_R(gf)$. For any 
$1$-morphisms $f'\colon A\to B$ and $g'\colon B\to A$ and any 
$2$-morphisms $\alpha\colon f\Rightarrow f'$ and $\beta\colon g\Rightarrow g'$ we have 
$$\theta_{f,g}^R\cdot\Tr_R(\beta\circ\alpha)=\Tr_R(\alpha\circ\beta)\cdot
\theta_{f',g'}^R.$$
\noindent This last property of $\theta_R$ can be called its {\rm naturality 
property}.\par
\noindent {\bf Proof}. The proof is identical to the proof in the case of pivotal 
categories except that the essential identities are now $2$-isomorphisms. 
Explicitly we find
$$\theta_{f,g}^R=(1_{i_A(f\otimes A^*)}\circ(T_B^{\dag}\otimes A^*)\circ 
1_{(g\otimes A^*)e_{A^*}})\cdot (1_{i_A(f\otimes A^*)(B\otimes i_{B^*}\otimes 
A^*)}\circ \bigotimes\nolimits_{e_{B^*},(g\otimes A^*)e_{A^*}})$$
$$\cdot (1_{i_A(f\otimes 
A^*)}\circ \phi_g\circ 1_{e_{B^*}})\cdot (\bigotimes\nolimits_{i_A(f\otimes A^*),i_B(g\otimes B^*)}\circ 1_{(B\otimes 
e_A \otimes B^*)e_{B^*}})$$
$$\cdot (1_{i_B(g\otimes B^*)(i_A\otimes A\otimes B^*)}
\circ (\bigotimes\nolimits_{f,e_A}\otimes B^*)\circ 1_{e_{B^*}})\cdot 
(1_{i_B(g\otimes B^*)}\circ(T_A\otimes B^*)\circ 1_{(f\otimes B^*)e_{B^*}}).$$
Since all the 2-morphisms in this formula are unitary, it follows that 
$\theta(f,g)$ is unitary also. The ``naturality'' of $\theta$ follows 
immediately from the formula above and the ``naturality'' of the triangulator 
and the tensorator.\qed

Of course there is also a 
family of natural unitary 2-morphism $\theta_L$ for the left trace functor. 
Now if we take $g=f^*$ then we 
can write 
$$\Tr_R(ff^*)=i_A(f\otimes A^*)(f^*\otimes A^*)i^*_A=i_A(f\otimes A^*)
(i_A(f\otimes A^*))^*,$$ 
so we have the $2$-morphism
$$\cap_f^R=i_{i_A(f\otimes A^*)}\colon 1_I\Rightarrow \Tr_R(ff^*).$$
\noindent Analogously we have the $2$-morphisms 
$$\cup_f^R=i^*_{i_A(f\otimes A^*)}\colon \Tr_R(ff^*)\Rightarrow 1_I.$$ 
We will call these $2$-morphisms {\sl cap} and {\sl cup} respectively. 
Analogously there exist caps and cups with respect to the left trace functor:
$$\cap_f^L=i_{i_{A^*}(A^*\otimes f)}\colon 1_I\Rightarrow \Tr_L(ff^*)$$ 
\noindent and
$$\cup_f^L=i^*_{i_{A^*}(A^*\otimes f)}\colon \Tr_L(f^*f)\Rightarrow 1_I.$$
These names are of course inspired by what these $2$-morphisms 
stand for in the 2-category of 2-tangles (see [5,6,30]). For our purposes we 
want these cups and 
caps to be compatible with the pivotal condition. 

\advance\thmno by 1
\proclaim\Def. A semi-strict pivotal $2$-category is called 
{\rm consistent} if the following conditions are satisfied
\item{1)} $\cap_f^R\cdot\theta_{f,f^*}^R=\cap_{f^*}^R$
\item{2)} $\theta_{f,f^*}^R\cdot\cup_{f^*}^R=\cup_f^R$
\item{3)} $\cap_f^L\cdot\theta_{f,f^*}^L=\cap_{f^*}^L$
\item{4)} $\theta_{f,f^*}^L\cdot\cup_{f^*}^L=\cup_f^L.$\par

\noindent In the sequel pivotal $2$-categories will always be assumed to be consistent.  
\advance\thmno by 1
\proclaim\Def. A {\rm semi-strict spherical $2$-category} is a semi-strict  
consistent pivotal $2$-category $C$ such that  
for any object $A$ in $C$ and 
any $1$-morphism $f\in\End(A)$ we have a unitary $2$-morphism  
$\sigma_f\colon\Tr_L(f)\Rightarrow \Tr_R(f)$. For any $1$-morphism 
$g\in\End(A)$ and any $2$-morphism $\alpha\colon f\Rightarrow g$ these 
$2$-isomorphisms are required to satisfy
$$\sigma_f\cdot\Tr_R(\alpha)=\Tr_L(\alpha)\cdot\sigma_g.$$
Furthermore we require the following identities to be satisfied:  
$$\cap_{1_{A}}^L\cdot\sigma_{1_A}=\cap_{1_{A}}^R\ \ \hbox{and}\ \ \sigma_{1_A}
\cdot \cup_{1_{A}}^R=\cup_{1_{A}}^L.$$\par

\noindent The last two conditions just mean that  
 ``cupping'' and ``capping'' are compatible with the spherical 
condition. 

In a semi-strict spherical $2$-category $C$ we can define a  
symmetric pairing 
$$\langle \cdot,\cdot\rangle \colon \THom(f,g)\otimes\THom(g,f)\to \TEnd(1_I).$$

\advance\thmno by 1
\proclaim\Def. Let $A,B$ be any objects in $C$ and $f,g\colon A\to B$ be any 
1-morphisms. 
For any $\alpha\in\THom(f,g)$ and any $\beta\in\THom(g,f)$ we 
define a {\rm pairing} by 
$$\langle \alpha,\beta\rangle =\cap_f^R\cdot\Tr_R((\alpha\cdot\beta)\circ 
1_{f^*})\cdot \cup_f^R.$$\par 

\noindent Notice that the analogous definition of a pairing 
, using the left trace functor,  
gives exactly the same pairing because $C$ is assumed to be spherical. 
Another way to understand this pairing is by noting that for any 
$2$-morphism $\alpha\colon f\Rightarrow f$, with $f\colon A\to B$ an 
arbitrary $1$-morphism, there is ``a kind of 
trace'' defined by 
$$1_A\buildrel i_f \over\Rightarrow ff^*\buildrel\alpha\circ f^*\over
\Rightarrow ff^*\buildrel i^*_f\over\Rightarrow 1_A.$$ In order to 
get a real trace we have to map this to $\TEnd(1_I)$. The right 
way to do this is by first applying the trace-functor to the ``traced'' 
$2$-morphism above and then close this up by adding cups and caps. This leads 
precisely to our definition of the pairing. The 
following lemmas show that the pivotal condition implies that 
one gets the same trace if one uses 
$$1_B\buildrel i_{f^*} \over\Rightarrow f^*f\buildrel f^*\circ\alpha
\over\Rightarrow f^*f\buildrel i^*_{f^*}\over\Rightarrow 1_B.$$

\advance\thmno by 1
\proclaim\Lem. $\langle \alpha,\beta\rangle =\langle \beta,\alpha\rangle $ for any 
$\alpha\in\THom(f,g)$ and any $\beta\in\THom(g,f)$.\par
\noindent{\bf Proof}. The proof is identical to the proof in the 
case of pivotal categories. Just as in that case it is easy to show  
that 
$$\langle \alpha,\beta\rangle =\cap_f^R\cdot\Tr_R(\alpha\circ{^{\dag}
\beta})\cdot \cup_f^R$$
$$=\cap_f^R\cdot\Tr_R(\alpha\circ
\beta^{\dag})\cdot \cup_f^R=\langle \beta,\alpha\rangle .$$
\qed

\advance\thmno by 1
\proclaim\Lem. $\langle \alpha,\beta\rangle ^*=\langle \alpha^*,\beta^*\rangle $ for any 
$\alpha\in\THom(f,g)$ and any $\beta\in\THom(g,f)$.\par 
\noindent{\bf Proof}. From the definition of the pairing we see 
$$\langle \alpha,\beta\rangle ^*=\langle \beta^*,\alpha^*\rangle .$$
\noindent Now use the previous lemma.
\qed

\advance\thmno by 1
\proclaim\Lem. $\langle \alpha,\beta\rangle =\langle \alpha^{\dag},\beta^{\dag}\rangle $ for any 
$\alpha\in\THom(f,g)$ and any $\beta\in\THom(g,f)$.\par
\noindent{\bf Proof}.
From the pivotal condition we get 
$$\langle \alpha,\beta\rangle =\cap_f^R\cdot\Tr_R(1_f\circ 
(\beta^{\dag}\cdot\alpha^{\dag}))\cdot \cup_f^R.$$
\noindent Now use $\Tr_R(ff^*)\cong\Tr_R(f^*f)$ and compatibility 
with cupping and capping. We get 
$$\cap_{f^*}^R\cdot\Tr_R((\beta^{\dag}\cdot\alpha^{\dag})\circ 1_f)
\cdot \cup_{f^*}^R=\langle \beta^{\dag},\alpha^{\dag}\rangle .$$
\noindent Finally use lemma 2.10 again.
\qed

\advance\thmno by 1
\proclaim\Lem. $\langle \alpha^{\#},\beta^{\#}\rangle =\langle {^{\#}\alpha},{^{\#}\beta}\rangle =
\langle \alpha,\beta\rangle $ for any 
$\alpha\in\THom(f,g)$ and any $\beta\in\THom(g,f)$.\par
\noindent{\bf Proof}. If $\alpha=\beta=1_{1_{A}}$, then this follows 
immediately from condition 13 in definition 2.3 and the spherical condition. 
The general case now follows from this particular case by using the pivotal 
condition.\qed
 
The following lemma shows that our setup is really a generalization 
of the setup of Crane, Kauffman and Yetter in [17] (see [35] for the 
definition of a ribbon category). 
\advance\thmno by 1
\proclaim\Lem. Let $C$ be a semi-strict consistent pivotal $2$-category, 
then 
${\rm End}(I)$ is a ribbon category.\par
\noindent{\bf Proof}. We already know that $\End(I)$ is a braided 
monoidal category (see lemma 2.2). Duality of course follows from the 
duality on $C$. The dual object of $f\in\End(I)$ is $f^*\in\End(I)$. 
The evaluation on $f$ is defined by $i_f$ and the coevaluation by 
$e_f$ and the identities they should satisfy are exactly those of 
condition~11 in  
Definition~2.3. The dual morphism of $\alpha\in\THom(f,g)$ is $\alpha^{\dag}\in
\THom(g^*,f^*)$. 

$\End(I)$ is pivotal since we have imposed the condition  
$\alpha^{\dag}={^{\dag}\alpha}$, which is equivalent to condition 12 in definition 2.3, as 
we explained. 

The ribbon structure in $\End(I)$ is defined by the 
family of unitary $2$-morphisms 
$$\phi_f\colon f={^{\#}f}\Rightarrow f^{\#}=f$$ 
\noindent indexed by all $f\in\End(I)$, which define the pivotal structure 
of $C$ (see definition 2.4). 
As we already mentioned below definition 2.4, conditions 2 and 3 in that 
definition are exactly the conditions that make $\phi_f$ into a twist. 
Note that the consistency of the pivotal structure as defined in definition 
2.7 restricted to the case $f\in\End(I)$ translates into the 
condition that $i_{f^*}$ and $e_{f^*}$, which define the structure of 
the right dual $f^*$, be compatible with the structure of the left dual 
$f^*$, defined by $i_f$ and $e_f$, and the twist $\phi_f$. In order to see 
this 
one should realise that in this case $\theta_{f,f^*}^R$ is equal to 
$(f\otimes \phi_{f^*})\cdot \otimes_{f,f^*}$ (see the formula in the proof of 
lemma 2.6).\qed
      
Before we go on let us have a look at an example of a spherical 
$2$-category. In [3] Baez defines the $2$-category of 2-Hilbert 
spaces, 2Hilb, and equips it with a tensor product and duals. Unfortunately 
2Hilb is not a semi-strict monoidal 2-category with duals: the composition 
of the 1-morphisms is not strict, the monoidal 
structure is not semi-strict and the duals do not satisfy all the required 
identities up to the nose either. Let us recall the definitions from 
[3] and some of Baez's results 
and then indicate how to obtain a semi-strict version of $\THilb$, the 
2-category of completely coordinatized 2-Hilbert spaces, $\THilb_{\rm cc}$. 
 We also define duals in $\THilb_{\rm cc}$ and show that the pivotal and the 
spherical conditions are 
satisfied. It is not difficult to show that $\THilb_{\rm cc}$ is equivalent to $\THilb$, 
but we omit the proof here. 
\medskip
\advance\thmno by 1
\noindent{\bf Example \the\chno.\the\thmno.} The objects in 2Hilb, 
which are called 2-Hilbert spaces, 
are all finite dimensional abelian $H^*$-categories. 
An abelian $H^*$-category $H$ is 
an abelian category such that the 
Hom-spaces are Hilbert spaces and composition is bilinear, and additionally 
$H$ is equipped with antilinear maps $*\colon \hom(x,y)\to \hom(y,x)$ 
for all objects $x,y$ in $H$ such that
\item{1.} $f^{**}=f,$
\item{2.} $(fg)^*=g^*f^*,$
\item{3.} $\langle fg,h\rangle =\langle g,f^*h\rangle ,$ 
\item{4.} $\langle fg,h\rangle =\langle f,hg^*\rangle $

\noindent for all $f\colon x\to y$, $g\colon y\to z$ and $h\colon x\to z$. 
It is shown in [3] that any abelian $H^*$-category $H$ is semi-simple as an 
abelian category and if it is finitely semi-simple than its dimension 
is defined as the number of objects in a basis of $H$. It is also shown 
that any basis of a 2-Hilbert space 
has the same number of objects, so its dimension is well defined. Furthermore 
Baez shows that two 2-Hilbert spaces are equivalent if and only if they have 
the same dimension.  

A $1$-morphism $F\colon H\to H'$ in 2Hilb is an exact functor such that 
$F\colon \hom(x,y)\to \hom(F(x),F(y))$ is linear and $F(f^*)=F(f)^*$ for 
all $f\in \hom(x,y)$. 

A $2$-morphism $\alpha\colon F\Rightarrow F'$ in 2Hilb is a natural 
transformation.

This all looks a little abstract, but, since we can always choose bases in 
2-Hilbert spaces, $1$-morphisms correspond to matrices with integer 
coefficients. This correspondence is reliable because any 2-Hilbert space 
is unitarily equivalent to a skeletal 2-Hilbert space, which is one where 
isomorphic objects are equal. So given a basis $\{a_i\}$ in a 
2-Hilbert space $A$ and a basis $\{b_i\}$ in a 2-Hilbert space $B$ any 
$1$-morphism $F\colon A\to B$ can be presented by the matrix $(F_i^j)$ with 
coefficients $F_i^j\in {\bf N}$, where 
$$F(a_i)=\bigoplus_j F_i^jb_j.$$ 
\noindent $2$-Morphisms now correspond to matrices of complex matrices. If 
  $(F_i^j)$ presents the $1$-morphism $F\colon A\to B$ and 
$(G_i^j)$ 
the $1$-morphism $G\colon A\to B$, then we can write a $2$-morphism 
$\alpha\colon F\to G$ as the matrix $(\alpha_i^j)$, where $\alpha_i^j$ 
is a $G_i^j\times F_i^j$ matrix with complex coefficients. Note that 
in our notation $(X_i^j)$ denotes the matrix with coefficients $X_i^j$.  

The tensor product is a little bit complicated in 2Hilb and we will not define 
it here in a basis invariant way. Roughly speaking the tensor product of two 
2-Hilbert spaces $A$ and $B$ can be obtained in the obvious way: 
choose a basis $\{a_i\}$ in $A$ and a basis $\{b_i\}$ in $B$ and ``define'' 
$A\otimes B$ as the $2$-Hilbert space with basis $\{a_i\otimes b_j\}$ and 
define $\hom(a_i\otimes b_j,a_i\otimes b_j)=\hom(a_i,a_i)\otimes 
\hom(b_j,b_j)$. It is obvious that Hilb, 
the category of finite dimensional 
Hilbert spaces, is the identity object. 
The braiding comes from the ordinary transposition of factors in the tensor product and we will not say more about it because it is not important 
for our purpose. 

Let us now have a look at the duality in 2Hilb. 
The dual of a 2-Hilbert 
space $H$ is the 2-Hilbert space $H^*=\Hom(H,\Hilb)$. There is always a 
unique dual basis in $H^*$ for each basis in $H$, up to isomorphism of course. 
The coevaluation $i_H$ and evaluation $e_H$ are now defined as usual. 
Let us assume that $H$ is skeletal. Given a basis $\{h_i\}$ in $H$ and 
its dual basis $\{h^i\}$ in $H^*$ we define $i_H\colon \Hilb\to H\otimes H^*$ 
by $${\bf C}\mapsto \bigoplus_{i} h_i\otimes h^i$$
\noindent and $e_H\colon H^*\otimes H\to \Hilb$ by 
$$h^i\otimes h_j\mapsto \delta^i_j {\bf C},$$ 
\noindent where $\delta^i_j$ is the Kronecker-delta. If $H$ is skeletal, 
then $T_H$, the triangulator, is trivial. Since, as already mentioned, any 
2-Hilbert space is unitarily equivalent to a skeletal one, this defines the 
triangulator in general.

The dual of a $1$-morphism $F\colon A\to B$ presented by the matrix 
$(F_i^j)$ is defined by the transpose of this matrix, i.e. 
$$F^*(b_j)=\bigoplus_{i}F_i^ja_i.$$
\noindent Baez shows that this really defines a left and right adjoint to 
$F$, so the 2-coevaluation $i_F$ and the 2-evaluation $e_F$ are easy to 
define. The Hilbert space $\hom(F(a_i),F(a_i))$ is isomorphic to  
$\hom(a_i,FF^*(a_i))$ for every $i$, so $i_F\colon 1_A\Rightarrow FF^*$ 
is simply defined as the natural transformation corresponding to 
the identity on $F(a_i)$ for each 
basis element $a_i$ under this isomorphism. 
In the same way we obtain the 2-evaluation $e_F\colon F^*F\Rightarrow 1_B$ by 
the isomorphism $\hom(F^*(b_i),F^*(b_i))\cong \hom(F^*F(b_i),b_i)$. Note 
that here we continue to use the convention under which $FF^*$ means first 
$F$ and then $F^*$ and not the other way around as is the more usual 
convention for functors and natural transformation, though not for 2-categories. 
In order to understand what $i_H$ is in more concrete terms one 
should note that the matrix corresponding to $FF^*$ has diagonal 
coefficients that are sums of squares of coefficients of $F$ and 
that $1_A$ is just a diagonal matrix with all diagonal coefficients 
equal to $1$. So $i_H$ is defined by the sums of the coevaluation 
maps on the terms in the diagonal coefficients and by the zero map 
for all non-diagonal coefficients of $FF^*$. In the same way we 
see that $e_H$ is defined by the ``ordinary'' evaluation maps on 
the terms of the diagonal coefficients of $F^*F$, which are squares 
also of course, and by the zero map on all the non-diagonal 
coefficients. 

The dual of a $2$-morphism $\alpha\colon F\Rightarrow G$ presented by 
the matrix $(\alpha_i^j)$ is the $2$-morphism $\alpha^*\colon G
\Rightarrow F$ presented by the matrix $((\alpha^{*})_i^j)$, where 
$(\alpha^*)_i^j$ is the adjoint of $\alpha_i^j$, obtained by taking the 
transpose and then the complex conjugate of each coefficient. A little 
thinking shows that $\alpha^{\dag}$ is presented by the adjoint of the 
whole matrix $(\alpha_i^j)$. It is now obvious that 
$\alpha^{\dag *}=\alpha^{*\dag}$ corresponds to the matrix 
$((\alpha^{\dag *})_i^j)$ where $(\alpha^{\dag *})_i^j=\alpha_j^i$ 

We now define $\THilb_{\rm cc}$. The objects are the non-negative integers. 
A 1-morphism between $n$ and $m$ is an $m\times n$ matrix with non-negative 
integer coefficients. A 2-morphism between two 1-morphisms $(F_i^j), 
(G_i^j)\colon n\to m$ is an  
$m\times n$ matrix $(\alpha_i^j)$, where $\alpha_i^j$ is a $G_i^j\times 
F_i^j$ matrix with complex coefficients. The various compositions are 
defined as for the matrices representing 1- and 2-morphisms in 2Hilb. 

The tensor product of two objects $n$ and $m$ is defined by $n\otimes m=nm$, 
the identity object being $1$. The tensor product of a 1- or 2-morphism $(X_i^j)$ 
and an object $n$ is defined by the matrix $((X\otimes n)_{ik}^{jl})$ with 
coefficients $(X\otimes n)_{ik}^{jl}=X_i^j\delta_k^l$. Likewise we define 
$(n\otimes X)_{ik}^{jl}=\delta_i^jX_k^l$. Note that for any 1-morphisms 
$(F_i^j)\colon n\to m$ and $(G_k^l)\colon r\to s$ 
we have 
$$[(F\otimes r)(m\otimes G)]_{ik}^{jl}=F_i^jG_k^l$$ 
\noindent and
$$[(n\otimes G)(F\otimes s)]_{ik}^{jl}=G_k^lF_i^j.$$ 
\noindent The tensorator is non-trivial: 
$$(\bigotimes\nolimits_{F,G})_{ik}^{jl}=P_{F_i^jG_k^l},$$ 
\noindent where $P_{F_i^jG_k^l}$ is the $F_i^jG_k^l\times F_i^jG_k^l$ 
permutation matrix with coefficients $(P_{F_i^jG_k^l})_{xy}^{wz}=\delta_x^z\delta_y^w$. 
This matrix clearly corresponds to the permutation of the factors in 
the tensor product $F_i^j\otimes G_k^l=F_i^jG_k^l$. It is clear that $\THilb_{\rm cc}$ is a 
semi-strict 2-category. Note that so far the definition of $\THilb_{\rm cc}$ 
coincides with the definition of $\TVect_{\rm cc}$, the totally coordinatized 
version of $\TVect$ defined by Kapranov and Voevodsky 
(see section 5.21 in [27]). 
Additionally we are going to define duals on $\THilb_{\rm cc}$, which is 
the extra structure that distinguishes $\THilb_{\rm cc}$ from 
$\TVect_{\rm cc}.$      

The dual of the object $n$ is $n$ itself. The coevaluation $i_n\colon 
1\to n^2$ is defined by the $n^2\times 1$ matrix with coefficients $(i_n)_i^1
=1$ if $i=(k-1)n+k$ for all $1\leq k\leq n$ and $(i_n)_i^1=0$ otherwise. 
In order to find this matrix we used the equivalence $\Hom(1,n\otimes n)
\simeq \Hom(n,n)$ and the fact that $i_n$ corresponds to the $n\times n$ 
identity matrix under this equivalence. Loosely speaking one just renumbers   
the coefficients of this identity matrix according to the following 
principle: the coefficient on the $k$-th row and the $l$-th column of 
the identity matrix becomes the coefficient on the $((k-1)n+l)$-th row of $i_n$. In an analogous way we see that the evaluation $e_n$ is given by the 
$1\times n^2$ matrix that is the transpose of $i_n$. The rest of the 
duality structure on $\THilb_{\rm cc}$ can be defined as for the matrices 
representing 1- and 2-morphisms in $\THilb$. With these definitions 
$\THilb_{\rm cc}$ becomes a semi-strict monoidal 2-category with duals as 
defined in definition 2.3.   
It is easy to check that $\THilb_{\rm cc}$ is a semi-strict spherical 
$2$-category. As a matter of fact it turns out that for any 
1-morphism $F\colon n\to m$ the 1-morphism $F^{\#}$ 
also corresponds to the transpose of $(F_i^j)$. 
\noindent The same is true for ${^{\#}F}$, so we see that $F^{\#}={^{\#}F}$, 
which suffices to conclude that 2Hilb is pivotal. For any 2-morphism 
$\alpha\colon F\Rightarrow G$ we find that $\alpha^{\#}$ corresponds to 
the matrix 
$((\alpha^{\#})_i^j)$ 
where $(\alpha^{\#})_i^j=\alpha_j^i$, as for $\alpha^{\dag *}$. 
Note that this does not mean that $\alpha^{\#}$ is equal to 
$\alpha^{\dag *}$, because 
$\alpha^{\#}$ is a $2$-morphism from $F^{\#}$ to $G^{\#}$ and 
$\alpha^{\dag *}$ is a $2$-morphism from $F^*$ to $G^*$. 

The left trace-functor 
applied to $FF^*$ gives the same result as when it is applied to 
$F^*F$, since this 
is the ordinary trace of the matrices corresponding to $FF^*$ and $F^*F$ 
respectively. For the same reason compatibility with cupping and capping is 
guaranteed. 
As a matter of fact it is easy to see directly what the cup and cap 
are for $\Tr_L(FF^*)=\Tr_L(F^*F)$. The 2-morphism $1_I$ is just 
equal to $1$, corresponding to the identity on ${\bf C}$, and 
$\Tr_L(FF^*)=\Tr_L(F^*F)$ is just the sum of the squares of the 
coefficients of $F$, so the cup $1_I\Rightarrow\Tr_L(FF^*)$ is 
nothing but the sum of the respective coevaluations and the 
cap $\Tr_L(FF^*)\Rightarrow 1_I$ is nothing but the sum of the 
respective evaluations.   

The left and the right trace functor are equal when applied to $1$-morphisms 
in $\THilb_{\rm cc}$, again because these are just ordinary 
matrix traces, so the spherical condition is also satisfied. 
\medskip

Having defined a semi-strict spherical $2$-category we now have to go 
a little further and add some linear structure. In order to define 
this linear structure we can work with a semi-strict $2$-category $C$ 
first. 
\advance\thmno by 1
\proclaim\Def. $C$ is called {\rm Vect-linear} if ${\rm Hom}(A,B)$ for any two objects $A,B$ in $C$ 
is a $2$-vector space of finite rank and the composition and the 
tensor product are ${\rm Vect}$-bilinear, with ${\rm Vect}$ the 
ring category of finite dimensional vector spaces over the fixed  
field ${\bf F}$.\par
\noindent For the definitions of a $2$-vector space and its rank 
and the definition of a ring category 
see [27].   
The following lemma is straight forward.

\advance\thmno by 1
\proclaim\Lem. The category ${\rm End}(I)$ is a braided monoidal ring 
category and ${\rm 2End}(1_I)$ is a commutative ring.\par

\noindent Each category $\Hom(A,B)$ becomes an $\End(I)$-module category 
with the action 
$$\End(I)\times\Hom(A,B)\buildrel\triangleright\over\to\Hom(A,B)$$
\noindent on objects defined by 
$$f\times g\mapsto f\triangleright g=(f\otimes A)(I\otimes g).$$ 
\noindent The isomorphism $(f_1f_2)\triangleright\ g\cong f_1
\triangleright (f_2\triangleright g)$ is defined by $\otimes_{f_1f_2,g}$ 
in the definition of a semi-strict monoidal $2$-category, see condition  
$(\rightarrow\rightarrow\otimes\cdot)$ in [27].
Notice that 
$$f\times g\mapsto (I\otimes g)(f\otimes B)$$
\noindent defines another action. This action is isomorphic to the 
one we have chosen by means of the isomorphisms $\otimes_{f,g}$.
The action on morphisms is defined in an obvious way now. 
For the rest of this paper I will assume that $\End(I)$ is equivalent  
to $\Vect$ and that $\TEnd(1_I)$ is 
a field isomorphic to ${\bf F}$. It is now easy to prove that the 
composition and the tensor product 
in $C$ are $\End(I)$-bilinear. It is obvious that the action of 
${\bf F}=\TEnd(1_I)$ on $\THom(f,g)$ for any $f,g\colon A\to B$ and for 
any $A$ and $B$ is the one induced by the action above. 
Notice that it makes sense to 
write $\Hom(A,B)\oplus\Hom(C,D)$ as the direct sum of two 
$\End(I)$-module categories and $\Hom(A,B)\otimes\Hom(C,D)$ as the tensor product of two 
$\End(I)$-module categories, see [27].

The following condition that we should impose on our $2$-categories concerns 
the non-degeneracy of the pairing defined in 2.9. 
Assume that $C$ is a Vect-linear semi-strict spherical $2$-category. 
Notice that the pairing is bi-linear. 
\advance\thmno by 1
\proclaim\Def. We call $C$ {\rm non-degenerate} if the pairing defined in 
2.9 is non-degenerate.\par

\noindent As in the case of ${\bf F}$-linear spherical categories [7,9] 
one can always take a non-degenerate quotient of 
an additive semi-strict spherical $2$-category. 
\advance\thmno by 1
\proclaim\Lem. Let $C$ be as above. Let $J$ be the {\rm Vect}-linear subcategory with 
the same objects and $1$-morphisms, but with $\THom_J(f,g)$ being the 
sub-vector space of $\THom_C(f,g)$ defined by 
$$\THom_J(f,g)=\lbrace \alpha\in\THom_C(f,g)| \langle \alpha,\beta\rangle =0\quad\forall
\beta\in\THom_C(g,f)\rbrace.$$ Then $C/J$ is a {\rm Vect}-linear semi-strict  
spherical $2$-category.\par
\noindent{\bf Proof}. It is clear that the vertical composition of 
$2$-morphisms is well defined in $C/J$. 

Let us now prove that the horizontal composition is well defined also. 
Let $f,g\colon X\to Y$ and $f',g'\colon Y\to Z$ be $1$-morphisms in $C$. 
For any $\alpha\in\THom_J(f,g)$, any $\beta\in\THom_C(f',g')$ and any 
$\gamma\in\THom_C(gg',ff')$ we have 
$$\langle \alpha\circ\beta,\gamma\rangle =\langle \alpha,(g\circ i_{f'})\cdot (g\circ\beta
\circ (f')^*)\cdot (\gamma\circ(f')^*)\cdot(f\circ i^*_{f'})\rangle =0.$$ 

Next we show that the tensor product is well defined in $C/J$. 
If $\alpha\in\THom_J(f,g)$, then $\alpha\otimes W\in\THom_J(f\otimes W, 
g\otimes W)$ and $W\otimes\alpha\in\THom_J(W\otimes f,W\otimes g)$ for 
any object $W$. 
The proof of these facts is not difficult, but is a bit cumbersome 
because we have to keep track of both the vertical and the horizontal  
composition. Writing out everything carefully gives  
$$\langle \alpha\otimes W,\beta\rangle $$
$$=\langle \alpha,(g\circ i_{(B\otimes i_W)})\cdot 
(\otimes_{g,i_W}\circ(B\otimes 1_{i^*_W}))$$
$$\cdot((A\otimes 1_{i_W})
\circ(\beta\otimes W^*)\circ(B\otimes 1_{i^*_W}))\cdot (\otimes^{-1}_{f,i_W}
\circ(B\otimes 1_{i^*_W}))\cdot(f\circ i^*_{(B\otimes i_W)})\rangle =0.$$
\noindent Notice that the long formula defining the second 
$2$-morphism 
in the second pairing is really a $2$-morphism from $g$ to 
$f$, so that its pairing with $\alpha$ makes sense. 
This shows our first assertion. The proof of the second is 
analogous and we leave the details to the reader. 

Lemma 2.11 shows that $\alpha^*\in\THom_J(g,f)$ if and only if 
$\alpha\in\THom_J(f,g)$, lemma 2.12 shows that $\alpha^{\dag}\in
\THom_J(g^*,f^*)$ if and only if $\alpha\in\THom_J(f,g)$ and lemma 
2.13 shows that $\alpha^{\#}\in\THom_J(f^{\#},g^{\#})$ if and only if 
$\alpha\in\THom_J(f,g)$. 

Now take the quotients in $C/J$ of all the structural $2$-morphisms 
involved in the definition of the tensor product, the duality and the 
pivotal and the spherical condition. Then $C/J$ becomes a {\rm Vect}-linear 
non-degenerate semi-strict spherical $2$-category.\qed    

Now let us define semi-simplicity for Vect-linear $2$-categories. A 
non-zero object $A$ in $C$ is an object for which we have $\End(A)\ne 0$. 
\advance\thmno by 1
\proclaim\Def. Let $C$ be a {\rm Vect}-linear semi-strict monoidal $2$-category. 
We say that $C$ is {\rm finitely semi-simple} if the following condition is  
satisfied:
\item{} There is a finite set of non-equivalent non-zero 
objects ${\cal E}$ such that for any pair of objects $A,B$ in $C$ we have  
$$\bigoplus_{X\in{\cal E}}\Hom(A,X)\otimes\Hom(X,B)\simeq\Hom(A,B).$$ 
\noindent The equivalence is given by the obvious composition of 
$1$- and $2$-morphisms.\par

\noindent Note that each Hom-space is finitely semi-simple as an ${\bf F}$-linear  
category because it is a $2$-vector space of finite rank. Let us fix a 
basis ${\cal F}_{A,B}\in \Hom(A,B)$ for any objects $A$ and $B$. We will 
usually 
refer to this semi-simplicity of the Hom-spaces as ``vertical 
semi-simplicity''. The semi-simplicity in the definition above we will always 
refer to as ``horizontal semi-simplicity''. Note that $\THilb_{\rm cc}$ 
is finitely semi-simple with 
the unique simple non-zero object 1. It is also non-degenerate, because 
the pairing of a 2-morphism $\alpha$ with its dual is equal to the sum of the 
squares of the absolute values of the complex coefficients of the matrix 
representing $\alpha$, which of course is non-zero if $\alpha$ is non-zero.   

\advance\thmno by 1
\proclaim\Def. Let $C$ be a {\rm Vect}-linear semi-strict monoidal $2$-category. 
An object $A$ in $C$ is called 
{\rm simple} if $\rk(\End(A))=1$.\par
\noindent Here $\rk(\End(A))$ is the rank of $\End(A)$ as a 2-vector space, 
see [27]. 
\advance\thmno by 1
\proclaim\Lem. Assume that $C$ is a {\rm Vect}-linear semi-strict monoidal 
finitely semi-simple $2$-category. An object $A$ in $C$ is simple 
if and only if $A$ is equivalent to one in ${\cal E}$.\par
\noindent{\bf Proof}. The proof is identical to the one that proves the 
analogous statement about finitely semi-simple categories.   
It follows from the following facts:
$$\rk(X\oplus Y)=\rk(X)+\rk(Y),\ \ \rk(X\otimes Y)=\rk(X)\rk(Y),$$
for any 2-vector spaces $X$ and $Y$.
\noindent These identities can be found in [27].\qed

Let us define the quantum dimension of objects and $1$-morphisms 
in a finitely semi-simple non-degenerate semi-strict spherical $2$-category. 
\advance\thmno by 1
\proclaim\Def. Let $A$ be an object in $C$. Then we define its 
{\rm quantum dimension} to be 
$$\qd(A)=\langle 1_{1_A},1_{1_A}\rangle .$$\par
\advance\thmno by 1
\proclaim\Def. Let $f$ be a $1$-morphism in $C$. Then we define 
its {\rm quantum dimension} to be 
$$\qd(f)=\langle 1_f,1_f\rangle .$$\par

\advance\thmno by 1
\proclaim\Lem. For any simple objects $A,B,C,D,E$ and any $1$-morphisms 
$f\colon A\to B\otimes C$, $g\colon B\to D\otimes E$ and 
$h\colon C\to D\otimes E$ we have  
$$\qd(f(g\otimes C))=\qd(f)\qd(B)^{-1}\qd(g),$$
\noindent and
$$\qd(f(B\otimes h))=\qd(f)\qd(C)^{-1}\qd(h).$$\par
\noindent{\bf Proof}. Let us first look closely at $\qd(f)$. The 
$1$-morphism $ff^*\colon A\to A$ is just a ``multiple'' of the identity, 
because $A$ is simple. ``Multiple'' here means that there is a vector space 
$V(f)$ such that $ff^*=V(f)1_A$. So we get 
$$\qd(f)=\dim(V(f))\qd(A).$$ Of course the same holds for $g$ and $h$. 

Using $f(g\otimes C)(g^*\otimes C)f^*=V(g)ff^*$, we get
$$\qd(f(g\otimes C))=\qd(f)\dim(V(g))=\qd(f)\qd(B)^{-1}\qd(g).$$ 
\noindent In the same way we get 
$$\qd(f(B\otimes h))=\qd(f)\dim(V(h))=\qd(f)\qd(C)^{-1}\qd(h).$$
\qed 

\noindent Note also that $\qd(A^*)=\qd(A)$ for any object $A$ by the spherical 
condition, that $\qd(A\otimes B)=\qd(A)\otimes\qd(B)$ for any objects $A$ and 
$B$ by the pivotal and the spherical conditions, that $\qd(f^*)=\qd(f)$ for any $1$-morphism $f$ by lemma 2.12 and 
the fact that $1_f^{\dag}=1_{f^*}$, and that $\qd(f^{\#})=\qd({^{\#}f})=
\qd(f)$ for any $1$-morphism $f$ by lemma 2.13. 

Finally we have to define the dimension of $C$.

\advance\thmno by 1
\proclaim\Def. Let $C$ be a finitely semi-simple non-degenerate 
semi-strict spherical $2$-category. Its {\rm dimension} $K$ is defined as the 
number of equivalence classes of non-zero simple objects in $C$.\par

\noindent This definition may seem rather surprising at first. In the 
proof of invariance under the $1\rightleftharpoons 5$ Pachner move in lemma 
5.4 and the proof of the auxiliary lemma 5.6 we show that this defines the 
right weight for the vertices of the triangulation of a manifold. Probably it 
has to do with the fact that we define a simple object to be an object $A$ 
such that $\End(A)\simeq\Vect$, which has dimension 1 as a category, and not just 
any finitely semi-simple non-degenerate spherical category. In our concluding 
remarks we will say a little more about this.

\advance\chno by 1
\equano=0
\thmno=0
\beginsection\the\chno. The definitions of $Z_C(\pm(ijklm))$ and $I_C(M,T)$\par
In this section $M=(M,T)$ still denotes a triangulated manifold. We also 
assume that there is a total ordering on the vertices of the simplicial 
complex defining the triangulation $T$ and a total ordering on the  
$4$-simplices of the same complex. Let $C$ be a finitely semi-simple  
semi-strict non-degenerate spherical $2$-category with a fixed finite basis of 
simple objects $\cal E$ and for any objects $A$ and $B$ a fixed finite basis 
of simple $1$-morphisms ${\cal F}_{A,B}\in\Hom(A,B)$. The linear maps 
$Z_C(\pm(ijklm))$, the number $Z_C(M,T,\ell)$ and the state sum 
$I_C(M,T)$ obviously depend on the 
given $2$-category $C$, but we will suppress the subscript $C$ at all places 
where this does not lead to any confusion. 

We now label the edges $(ij)$ of the 
triangulation with simple objects $e_{ij}$ in $\cal E$ and the triangles 
$(ijk)$ with simple $1$-morphisms $f_{ijk}\in{\cal F}_{ijk}
={\cal F}_{e_{ik},e_{jk}\otimes e_{ij}}$, the chosen basis in $\Hom(e_{ik}, e_{jk}
\otimes e_{ij})=H(ijk)$. We will use the following notation for 
the different compositions of the $1$-morphisms:
$$f_{(ijk)l}=f_{ikl}(e_{kl}\otimes f_{ijk}),$$
$$f_{i(jkl)}=f_{ijl}(f_{jkl}\otimes e_{ij}),$$
$$f_{((ijk)l)m}=f_{(ikl)m}(e_{lm}\otimes e_{kl}\otimes f_{ijk})=
f_{ilm}(e_{lm}\otimes f_{ikl})(e_{lm}\otimes e_{kl}\otimes f_{ijk}),$$
$$f_{i(j(klm))}=f_{i(jkm)}(f_{klm}\otimes e_{jk}\otimes e_{ij})=
f_{ijm}(f_{jkm}\otimes e_{ij})(f_{klm}\otimes e_{jk}\otimes e_{ij}).$$

\noindent Let $+(ijklm)$ be the positively oriented standard $4$-simplex 
labelled as described above. We are going to define the linear map 
$$Z(+(ijklm))\colon \THom(f_{(ikl)m},f_{i(klm)})\otimes
\THom(f_{(ijk)m},f_{i(jkm)})\to$$
$$\THom(f_{(jkl)m},f_{j(klm)})\otimes 
\THom(f_{(ijl)m},f_{i(jlm)})\otimes\THom(f_{(ijk)l},f_{i(jkl)})$$
\noindent using the pairing $\langle \cdot,\cdot\rangle $ 
described in the previous section. Consider the linear map
$$\THom(f_{(ikl)m},f_{i(klm)})\otimes
\THom(f_{(ijk)m},f_{i(jkm)})\to \THom(f_{((ijk)l)m},f_{i(j(klm))})$$ 
defined by 
$$\beta\otimes\delta\mapsto $$
$$(\beta\circ (e_{lm}\otimes e_{kl}\otimes f_{ijk}))\cdot 
(f_{ikm}\circ\otimes_{f_{klm},f_{ijk}})\cdot (\delta\circ (f_{klm}\otimes 
e_{jk}\otimes e_{ij})).$$ Let us write the right-hand side as 
$\beta\delta$ as a shorthand. Consider also the linear map 
$$\THom(f_{j(klm)},f_{(jkl)m})\otimes 
\THom(f_{i(jlm)},f_{(ijl)m)})\otimes\THom(f_{i(jkl)},f_{(ijk)l)})\to $$
$$\THom(f_{i(j(klm))},f_{((ijk)l)m})$$ defined by 
$$\alpha\otimes\gamma\otimes\epsilon\mapsto$$
$$(f_{ijm}\circ(\alpha\otimes e_{ij}))\cdot(\gamma\circ(e_{lm}\otimes 
f_{ijk}\otimes e_{ij}))\cdot(f_{ilm}\circ(e_{lm}\otimes\epsilon)).$$ 
Let us write the right-hand side as $\alpha\gamma\epsilon$ 
as a shorthand. Now define the linear map 
$$\THom(f_{j(klm)},f_{(jkl)m})\otimes\THom(f_{(ikl)m},f_{i(klm)})\otimes
\THom(f_{i(jlm)},f_{(ijl)m})$$
$$\otimes\THom(f_{(ijk)m},f_{i(jkm)})
\otimes\THom(f_{i(jkl)},f_{(ijk)l})\to {\bf F}$$
by
$$\alpha\otimes\beta\otimes\gamma\otimes\delta\otimes\epsilon\mapsto$$
$$\langle \beta\delta,\alpha\gamma\epsilon\rangle .$$ 
Using $\THom(g,f)=\THom(f,g)^*$, following from the non-degeneracy of 
$\langle \cdot,\cdot\rangle $, this gives us $Z(+(ijklm))$.

\pageinsert
\line{\hfil\epsfbox{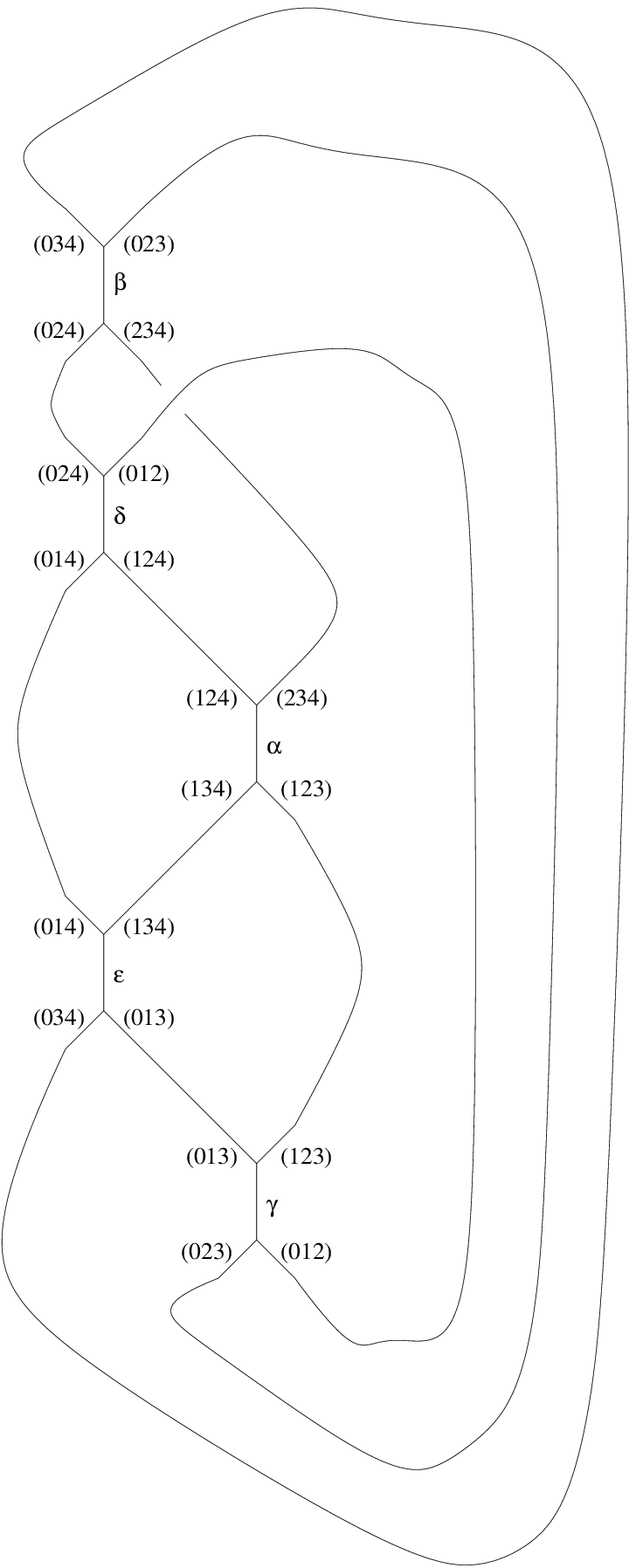}\hfil}
\vglue 5pt
\centerline{Picture 1}
\endinsert

Before we go on, let us have a look at a diagrammatic 
picture of $Z(+(01234))$. In picture 1 we have depicted a 
diagram that resembles a $15j$ symbol as defined in [17].
Here we have taken the dual 1-skeleton of the boundary of a 
4-simplex, a graph with 5 4-valent vertices and 10 edges, and ``split'' 
each 4-valent vertex into 2 trivalent vertices connected by an edge which 
we call a ``dumbbell''. Each dumbbell in the resulting graph corresponds to 
a tetrahedral face of the 4-simplex; the remaining 10 edges in this 
graph correspond to the triangular faces of the 4-simplex. Here we have 
labelled the dumbbells by the 2-morphisms $\alpha,$ $\beta,$ $\gamma,$ 
$\delta,$ $\epsilon,$ and labelled the remaining edges by triples 
$(ijk)$ corresponding to the 1-morphisms $f_{ijk}.$ Finally, the crossing 
is just $\otimes_{f_{234},f_{012}}$ (the tensorator that appears in the 
definition of $\beta\delta$).  
In this diagram we do not see the objects involved and so we 
do not see the effect of the trace functor and the final caps and 
cups either. But we do see the ``kind of trace'' mentioned after 
the definition of the pairing. 
Remember, if $f,g\colon A\to B$ are $1$-morphisms and 
$\zeta\colon f\Rightarrow g$ and $\xi\colon g\Rightarrow f$ 
are $2$-morphisms, then one obtains $\langle \zeta,\xi\rangle $ by first 
taking a ``kind of trace'', namely 
$$1_{A}\buildrel i_f\over\Rightarrow ff^*\buildrel \zeta\circ 
f^* \over\Rightarrow gf^*\buildrel \xi\circ f^*\over
\Rightarrow ff^*\buildrel i^*_{f}\over\Rightarrow 1_{A}.$$ This 
is exactly the diagram we see when we take $\zeta=\beta\delta$ 
and $\xi=\alpha\gamma\epsilon$, the respective $2$-morphisms 
involved in the definition of $Z(+(01234))$. The properties 
satisfied by the duality of the $1$-morphisms, condition 11 in 
definition 2.3, and the properties of the composition of 
$\otimes_{\cdot,\cdot}$, 
analogous to the properties of the braiding in a braided category, 
allow us to apply isotopies to this diagram, 
corresponding to Reidemeister moves $0$, $2$ and $3$, without changing 
the linear map $Z(+(01234))$ defined by the diagram. Since 
our $2$-category $C$ is also assumed to be spherical, we can 
also exchange closing strands on the left hand side for closing 
strands to the right-hand side of the diagram. The pivotal 
condition alone justifies these manipulations of the diagrams, 
but the spherical condition is needed for the 
different ways of nesting the cups and the caps that can not 
be seen in the diagrams. These are the ones we put on our 
pairing after applying the trace functor. This nice behaviour 
under different 
ways of nesting turns out to be necessary for our purposes, 
as can be seen in the proof of lemma 4.4. There this kind 
of diagram allows us to prove invariance of our invariant under 
any permutation of the vertices of the chosen triangulation. 
The formulas, which I had written out completely by hand at 
first, are far too large to fit on ordinary sheets of paper. 
So, although these diagrams do not define a complete diagrammatic 
calculus for my invariants, since they do not show the objects, 
they are certainly very helpful 
to see what is going on. Better diagrams would be very welcome, 
but so far I have not been able to invent them.  

Let us now define  
$$Z(-(ijklm))\colon \THom(f_{(jkl)m},f_{j(klm)})\otimes 
\THom(f_{(ijl)m},f_{i(jlm)})\otimes\THom(f_{(ijk)l},f_{i(jkl)})\to $$
$$\THom(f_{(ikl)m},f_{i(klm)})\otimes\THom(f_{(ijk)m},f_{i(jkm)}).$$
Consider the linear map 
$$\THom(f_{(jkl)m},f_{j(klm)})\otimes 
\THom(f_{(ijl)m},f_{i(jlm)})\otimes\THom(f_{(ijk)l},f_{i(jkl)})\to $$
$$\THom(f_{((ijk)l)m},f_{i(j(klm))})$$ defined by 
$$\alpha\otimes\gamma\otimes\epsilon\mapsto$$
$$(f_{ilm}\circ(e_{lm}\otimes\epsilon))\cdot(\gamma\circ(e_{lm}\otimes 
f_{jkl}\otimes e_{ij}))\cdot(f_{ijm}\circ(\alpha\otimes e_{ij})).$$
\noindent Call this $\epsilon\gamma\alpha$. Consider also the 
linear map
$$\THom(f_{i(klm)},f_{(ikl)m})\otimes
\THom(f_{i(jkm)},f_{(ijk)m})\to \THom(f_{i(j(klm))},f_{((ijk)l)m})$$ 
by 
$$\beta\otimes\delta\mapsto$$
$$(\delta\circ(f_{klm}\otimes e_{jk}\otimes e_{ij}))\cdot(f_{024}
\circ\otimes^{-1}_{f_{klm},f_{ijk}})\cdot(\beta\circ(e_{lm}\otimes e_{kl}
\otimes f_{ijk})).$$
\noindent Call this $\delta\beta$. Now define the linear map
$$\THom(f_{(jkl)m},f_{j(klm)})\otimes\THom(f_{i(klm)},f_{(ikl)m})\otimes
\THom(f_{(ijl)m},f_{i(jlm)})$$
$$\otimes\THom(f_{i(jkm)},f_{(ijk)m})
\otimes\THom(f_{(ijk)l},f_{i(jkl)})\to {\bf F}$$
defined by
$$\alpha\otimes\beta\otimes\gamma\otimes\delta\otimes\epsilon\mapsto$$
$$\langle \delta\beta,\epsilon\gamma\alpha\rangle .$$ 
Using $\THom(g,f)=\THom(f,g)^*$ again this gives us $Z(-(ijklm))$.
 
In the sequel let us write $2H(ijkl)$ for $\THom(f_{(ijk)l},f_{i(jkl)})$. 
Likewise let us write $2H(ijklm)$ for $\THom(f_{((ijk)l)m,i(j(klm))})$. 
The next lemma is the analogue of the Crossing Lemma 5.4 in [9].
\advance\thmno by 1
\proclaim\Lem ({\rm crossing}). The following diagram is commutative: 
$$\vbox{\offinterlineskip
\halign{\hfil#\hfil&\quad\hfil#\hfil&\hfil#\hfil\cr
$\bigoplus\limits_{f_{024}}2H(0234)\otimes 2H(0124)$&$
\buildrel\Phi_{01234}
\over
{\hbox to 50pt{\rightarrowfill}}$
&$\bigoplus\limits_{\scriptstyle e_{13},f_{013},
\atop\scriptstyle f_{123},f_{134}}
2H(1234)\otimes 2H(0134)\otimes 2H(0123)$\cr
$\Bigg\downarrow$ &&$\Bigg\downarrow$\cr
\noalign{\vskip7pt}
$2H(01234)$&$\longeq$&$2H(01234)$\cr}}$$ 
\noindent Here the vertical linear maps are the isomorphisms defined 
by the composition of the respective $2$-morphisms and $\Phi_{01234}$ is 
defined by 
$$\Phi_{01234}=\bigoplus_{\scriptstyle e_{13},f_{024},f_{013},\atop\scriptstyle
f_{123},f_{134}}Z(+(01234))\qd(e_{13})^{-1}\qd(f_{013})\qd(f_{123})\qd(f_{134}).$$ 
Furthermore the inverse of $\Phi_{01234}$ is given by  
$$\Phi^{-1}_{01234}=\bigoplus_{\scriptstyle e_{13},f_{024},f_{013},\atop\scriptstyle
f_{123},f_{134}}Z(-(01234))\qd(f_{024}).$$
\par
\noindent{\bf Proof}. First of all let us explain why the vertical 
maps are isomorphisms. Using the same notation as above we must 
show that for any $2$-morphism 
$\alpha\in\THom(f_{((012)3)4},f_{0(1(234))})$ there exist $2$-morphisms 
$\beta_{024}\in\THom(f_{(023)4},f_{0(234)})$ and $\gamma_{024}\in\THom(
f_{(012)4},f_{0(124)})$ such that 
$$\sum_{f_{024}}(\beta_{024}\circ(e_{34}\otimes e_{23}\otimes f_{012}))\cdot (f_{024}\circ  
\otimes_{f_{234},f_{012}})\cdot(\gamma_{024}\circ(f_{234}\otimes  
e_{12}\otimes e_{01}))=\alpha.$$ 
\noindent By (vertical) semi-simplicity we get the following 
(vertical) decomposition of $\alpha$:
$$\vbox{
\halign{\hfil#\hfil&\quad\hfil#\hfil\quad&#\cr
$f_{034}(e_{34}\otimes f_{023})(e_{34}\otimes e_{23}\otimes f_{012})$&
&\hfil $f_{034}(e_{34}\otimes f_{023})(e_{34}\otimes e_{23}
\otimes f_{012})$\hfil\cr 
\noalign{\vskip 5pt}
&&\hfil$\big\Downarrow$\hfil\cr
\noalign{\vskip 5pt}
&&\hfil$\sum\limits_{f_{024}}f_{024}(f_{234}\otimes e_{02})(e_{34}
\otimes e_{23}\otimes f_{012})$\hfil\cr
$\bigg\Downarrow$&$=$&\hfil\hfil\hskip 10pt$\big\Downarrow\   
\otimes^{-1}_{f_{234},f_{012}}$\hfil\cr
&&\hfil$\sum\limits_{f_{024}}f_{024}(e_{24}\otimes f_{012})(f_{234}
\otimes e_{12}\otimes e_{01})$\hfil\cr
&&\hfil$\big\Downarrow$\hfil\cr
\noalign{\vskip 5pt}
$f_{014}(f_{124}\otimes e_{01})(f_{234}\otimes e_{12}\otimes e_{01})$
&&\hfil$f_{014}(f_{124}\otimes e_{01})(f_{234}\otimes e_{12}
\otimes e_{01})$\hfil\cr}}$$
\noindent Notice that the $1$-morphisms 
$f_{024}(f_{234}\otimes e_{02})(e_{34}\otimes e_{23}\otimes 
f_{012})$ do not form a basis of simple $1$-morphisms of $\THom(
e_{04},e_{34}\otimes e_{23}\otimes e_{12}\otimes e_{01})$ in general, but they 
are generators by (horizontal) semi-simplicity and the fact that the  
$f_{024}$, the $f_{234}$ and the $f_{012}$ form bases of their 
respective Hom-spaces. So by summing the projections on the simple 
components of each $f_{024}(f_{234}\otimes e_{02})(e_{34}\otimes e_{23}\otimes 
f_{012})$ we get the decomposition in the diagram. By (horizontal) 
semi-simplicity we can now decompose the $2$-morphisms on the right 
side in the diagram above and we get the $\beta_{024}$ 
and $\gamma_{024}$ we were looking for. This proves that the 
left vertical map in the lemma is an isomorphism. In an analogous way one proves 
that the other vertical map in the lemma is an isomorphism. 

The commutativity of the diagram in the lemma now follows by taking 
arbitrary elements in $H(0234)$ and in $H(0124)$ respectively and 
pairing it with an arbitrary element in $H(01234)^*$. In order 
to get the multiplicative factors in the definition of $\Phi$ we have 
to use the identity 
$$\qd(f_{ikl}(e_{kl}\otimes f_{ijk}))=\qd(f_{ikl})\qd(e_{ik})^{-1}\qd(f_{ikl}).$$
This identity follows from lemma 2.25.

The inverse of $\Phi$ can be computed by reading the diagram the other 
way around and using similar arguments.\qed

We now define our state sum. Let $C$ be a finitely semi-simple 
non-degenerate semi-strict spherical $2$-category with non-zero 
dimension.
\advance\thmno by 1
\proclaim\Def. With the data and notations as above we define for 
every closed compact piece-wise linear $4$-manifold 
$M$ with triangulation $T$ the {\rm state sum}
$$I_C(M,T)=K^{-v}\sum_{\ell}Z_C(M,T,\ell)\prod_{e}\qd(\ell(e))^{-1}\prod_{f}
\qd(\ell(f)).$$
\noindent Here $v$ is the number of vertices in $T$. We sum over all the labellings $\ell$ and take the 
products over all the edges $e$ and all the faces $f$ in the 
triangulation.\par 

\advance\chno by 1
\thmno=0
\equano=0
\beginsection\the\chno. $Z_C(M,T,\ell)$ is a combinatorial invariant\par
In this section we show that $Z(M,T,\ell)=Z_C(M,T,\ell)$ is equal for all 
``isomorphic'' labellings and $Z(M,T,\ell)=Z(M',T',\ell')$ for any pair of 
triangulated manifolds $M$ and $M'$ with triangulations $T$ and $T'$ that are isomorphic 
under a combinatorial isomorphism and 
any pair of 
``compatible'' labellings $\ell$ and $\ell'$. 

First of all we have to show that $Z(M,T,\ell)$ does not depend 
on the choice of representatives $f_{ijk}\in{\cal F}_{ijk}$ nor 
on the choice of representatives $e_{ij}\in{\cal E}_{ij}$.
\advance\thmno by 1
\proclaim\Lem. Let $\ell$ and $\ell'$ be labellings with the 
same basis of objects ${\cal E}_{ij}$ but different 
bases of $1$-morphisms in ${\cal F}_{ijk}$. Then we have 
$Z(M,T,\ell)=Z(M,T,\ell')$.\par
\noindent{\bf Proof}. Let us denote the labels in $\ell$ by 
$f_{ijk}$ and the labels in $\ell'$ by $f'_{ijk}$. They are 
representatives of the same isomorphism classes so for any triple 
$ijk$ there is a $2$-isomorphism $\phi_{ijk}\colon f_{ijk}\Rightarrow 
f'_{ijk}$. These $2$-isomorphisms induce an isomorphism of 
vector spaces 
$$\THom(f_{(ijk)l},f_{i(jkl)})\cong\THom(f'_{(ijk)l},
f'_{i(jkl)})$$ given by 
$$\alpha\mapsto \phi^{-1}_{(ijk)l}\alpha\phi_{i(jkl)}.$$ 
\noindent Here $\phi_{(ijk)l}$ stands for $\phi_{ikl}(e_{kl}\otimes 
f_{ijk})$ and $\phi_{i(jkl)}$ for $\phi_{ijl}(\phi_{jkl}\otimes e_{ij})$.
Likewise there is an isomorphism between $\THom(f_{((ijk)l)m},f_{
i(j(klm))})$ and $\THom(f'_{((ijk)l)m},f'_{i(j(klm))})$. 

Now without writing out the explicit formulas, which is not very 
difficult but extremely tedious, one can see imediately the 
result of the lemma. The crossing lemma implies that the following 
diagram is commutative.
$$\vbox{\offinterlineskip
\halign{\hfil#\hfil&\quad#&\hfil#\hfil\cr
$\bigoplus\limits_{f_{ikm}}2H(iklm)\otimes 2H(ijkm)$&$\buildrel
 Z(ijklm)\over{\hbox to 40pt{\rightarrowfill}}$
&$\bigoplus\limits_{\scriptstyle e_{jl},f_{ijl},\atop
f_{jkl},f_{jlm}}2H(jklm)\otimes 2H(ijlm)\otimes 2H(ijkl)$\cr
\noalign{\vskip7pt}
$\Bigg\downarrow$&&$\Bigg\downarrow$\cr
\noalign{\vskip15pt}
$\bigoplus\limits_{f'_{ikm}}2H'(iklm)\otimes 2H'(ijkm)$&$
\buildrel Z'(ijklm)\over{\hbox to 40pt{\rightarrowfill}}$&
$\bigoplus\limits_{\scriptstyle 
e'_{jl},f'_{ijl},\atop
f'_{jkl},f'_{jlm}}2H'(jklm)\otimes 2H'(ijlm)\otimes 2H'(ijkl)$\cr}}$$

\noindent So $Z'(ijklm)$ is a conjugate of $Z(ijklm)$. Taking 
the respective traces over the tensor product of all the partition functions 
now shows that $Z(M,T,\ell)=Z(M,T,\ell')$.\qed

Now suppose we take a different basis of simple objects in $C$. 
In other words, suppose we have two labellings $\ell$ and $\ell'$ such 
that $\ell((ij))=e_{ij}$ is equivalent to $\ell'((ij))=e'_{ij}$ 
for every edge. Let us assume that there is an isomorphism 
$\phi_{ij}\colon e_{ij}\to e'_{ij}$ for every $i,j$ actually. 
In general $\phi_{ij}$ fails to be an isomorphism just by an invertible 
scalar, whose existence is not important for the following arguments. 
These isomorphisms induce a linear isomorphism 
$$\phi_{ijk}\colon\Hom(e_{ik},e_{jk}\otimes e_{ij})\to
\Hom(e'_{ik},e'_{jk}\otimes e'_{ij})$$
defined by 
$$f_{ijk}\mapsto \phi^{-1}_{ik}f_{ijk}(\phi_{jk}\otimes\phi_{ij})$$
\noindent for any $ijk$. 
Denote these $1$-morphisms by $f'_{ijk}$. It is clear that 
the simplicity of $f_{ijk}$ implies the simplicity of $f'_{ijk}$ 
for any $ijk$. Since we have proved  
that $Z(M,T,\ell)$ is independent of the choice of basis of simple $1$-morphisms 
in lemma 4.1 we may assume that $\ell'((ijk))=f'(ijk)$.

\advance\thmno by 1
\proclaim\Lem. Let $\ell$ and $\ell'$ be two labellings as 
described above. Then $Z(M,T,\ell)=Z(M,T,\ell')$.\par
\noindent{\bf Proof}. The identities 
$$f'_{ikl}(e'_{kl}\otimes f'_{ijk})=\phi^{-1}_{il}f_{ikl}
(\phi_{kl}\otimes\phi_{ik})(\phi^{-1}_{kl}\otimes\phi^{-1}_{ik}
(e_{kl}\otimes f_{ijk}))(\phi_{kl}\otimes\phi_{jk}\otimes\phi_{ij})
$$
$$=\phi^{-1}_{il}f_{ikl}(e_{kl}\otimes f_{ijk})(\phi_{kl}\otimes\phi_{jk}\otimes\phi_{ij})$$
\noindent and
$$f'_{ijl}(f'_{jkl}\otimes e'_{ij})=\phi^{-1}_{il}f_{ijl}
(\phi_{jl}\otimes\phi_{ij})(\phi^{-1}_{jl}\otimes\phi^{-1}_{ij}
(f_{jkl}\otimes e_{ij}))(\phi_{kl}\otimes\phi_{jk}\otimes\phi_{ij})$$
$$=\phi^{-1}_{ij}f_{ijl}(f_{jkl}\otimes e_{ij})(\phi_{kl}\otimes\phi_{jk}\otimes\phi_{ij})$$
\noindent show the existence of a linear isomorphism
$$2H(ijkl)=\THom(f_{(ijk)l},f_{i(jkl)})\cong\THom(f'_{(ijk)l},
f'_{i(jkl)})=2H'(ijkl).$$

Again by applying the crossing lemma we get the commutative 
diagram of the previous lemma, although the vertical arrows now represent 
different linear maps. So again $Z(ijklm)$ and $Z'(ijklm)$ 
are conjugates, which implies that $Z(M,T,\ell)$ and $Z(M,T,\ell')$ 
are equal.\qed  

Finally let us prove that $Z(M,T,\ell)$ does not depend on the 
ordering of the vertices in the triangulation $T$ of $M$. 

\advance\thmno by 1
\proclaim\Def. Let $\phi\colon T\to T'$ be a combinatorial 
isomorphism of two triangulations of $M$. Let the edge $(ij)$ of 
$T$ be labelled with a simple object $e_{ij}$ for all pairs of 
different vertices 
$i$ and $j$ of $T$ and the edge $(ij)$ 
of $T'$ with a simple object $e'_{ij}$ for all pairs of different 
vertices $i$ and $j$ of $T'$. Let a triangle $(ijk)$ 
of $T$ be labelled with a simple $1$-morphism $f_{ijk}$ for all 
triples of different vertices $i,j,k$ of $T$ and a 
triangle $(ijk)$ of $T'$ with a simple $1$-morphism $f'_{ijk}$ 
for all triples of different vertices $i,j,k$ of $T'$. We say 
that $\phi$ is {\rm compatible with the labellings} if the 
following conditions are satisfied:
\item{1)} $e'_{\phi(i)\phi(j)}=e_{ij}$ if $\phi$ preserves the orientation of the edge and 
$e'_{\phi(i)\phi(j)}=e_{ij}^*$ if it reverses the orientation.
\item{2)} If $\phi$ decomposes into the transposition $
(ijk)\mapsto (ikj)$ and a simplicial isomorphism, then 
$f'_{\phi(i)\phi(j)\phi(k)}$ is unitarily isomorphic to 
$(i_{e_{kj}}\otimes e_{ij})(e_{kj}\otimes f^*_{ijk})$. 
If the transposition in the decomposition 
of $\phi$ is $(ijk)\mapsto (kji)$, then 
$f'_{\phi(i)\phi(j)\phi(k)}$ is unitarily isomorphic to 
$(f_{ijk}^{\#})^*={^{\#}(f_{ijk}^*)}$.\par

\noindent If we represent $f_{ijk}$ and $f^*_{ijk}$ by the diagrams in picture 2, then 
we can represent $(i_{e_{kj}}\otimes e_{ij})(e_{kj}\otimes f^*_{ijk})$ as in picture 3 
and $(f_{ijk}^{\#})^*$ as in picture 4. Note that any combinatorial 
isomorphism $\phi$, when 
restricted to a triangle $(ijk)$ in $T,$ always decomposes into a 
permutation of $(ijk)$ and a simplicial isomorphism. Since 
$S_3$ is generated by the two transpositions in condition 2, we 
have really defined compatibility for any combinatorial isomorphism.  
Note also that in condition 2 we have only required the isomorphy 
of the corresponding $1$-morphisms. The reason for this is that we 
want the composition of a combinatorial isomorphism compatible with 
the labellings $\phi$ with its inverse $\phi^{-1}$ to be compatible with 
the labellings. This would not be the case if we required the corresponding 
$1$-morphisms to be equal because $^{\#}(f_{ijk}^{\#})$ is only isomorphic to 
$f_{ijk}$. For our purpose these isomorphisms do not matter, because 
we have already shown that $Z(M,T,\ell)$ is independent of the choice 
of representative simple $1$-morphisms in each isomorphism class. 
So, given a combinatorial isomorphism $\phi$ 
and a labelling of $T$, there is a unique compatible labelling only 
up to isomorphism.

\midinsert
\vglue 13pt
\line{\hfil\vbox{\epsfysize=2.5cm\epsfbox{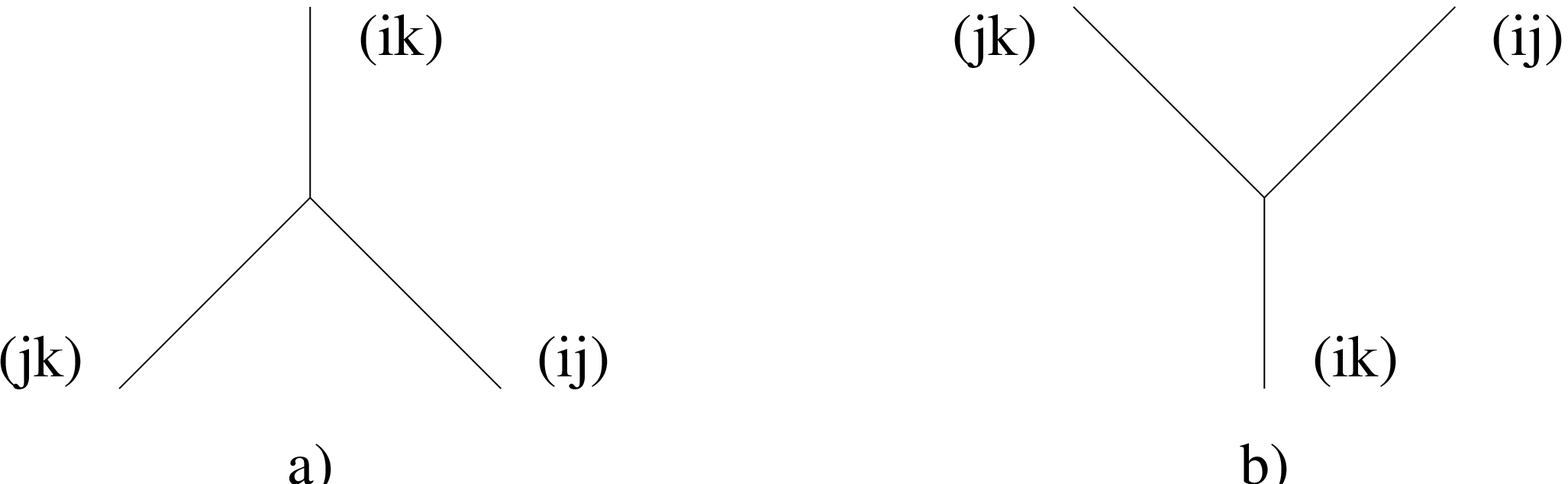}}\hfil}
\vglue 5pt
\centerline{Picture 2. a) $f_{ijk}$ b) $f^*_{ijk}$}
\vglue 13pt
\endinsert

\midinsert
\vglue 13pt
\line{\hfil\vbox{\epsfysize=2.7cm\epsfbox{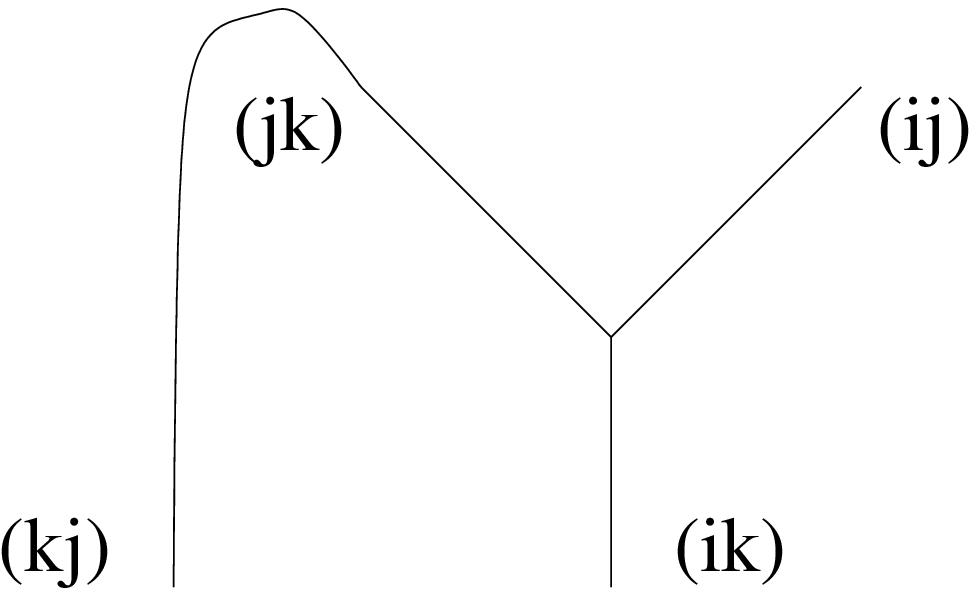}}\hfil}
\vglue 5pt
\centerline{Picture 3. $(i_{e_{kj}}\otimes e_{ij})(e_{kj}\otimes f^*_{ijk})$}
\vglue 13pt
\endinsert

\midinsert
\vglue 13pt
\line{\hfil\vbox{\epsfysize=3cm\epsfbox{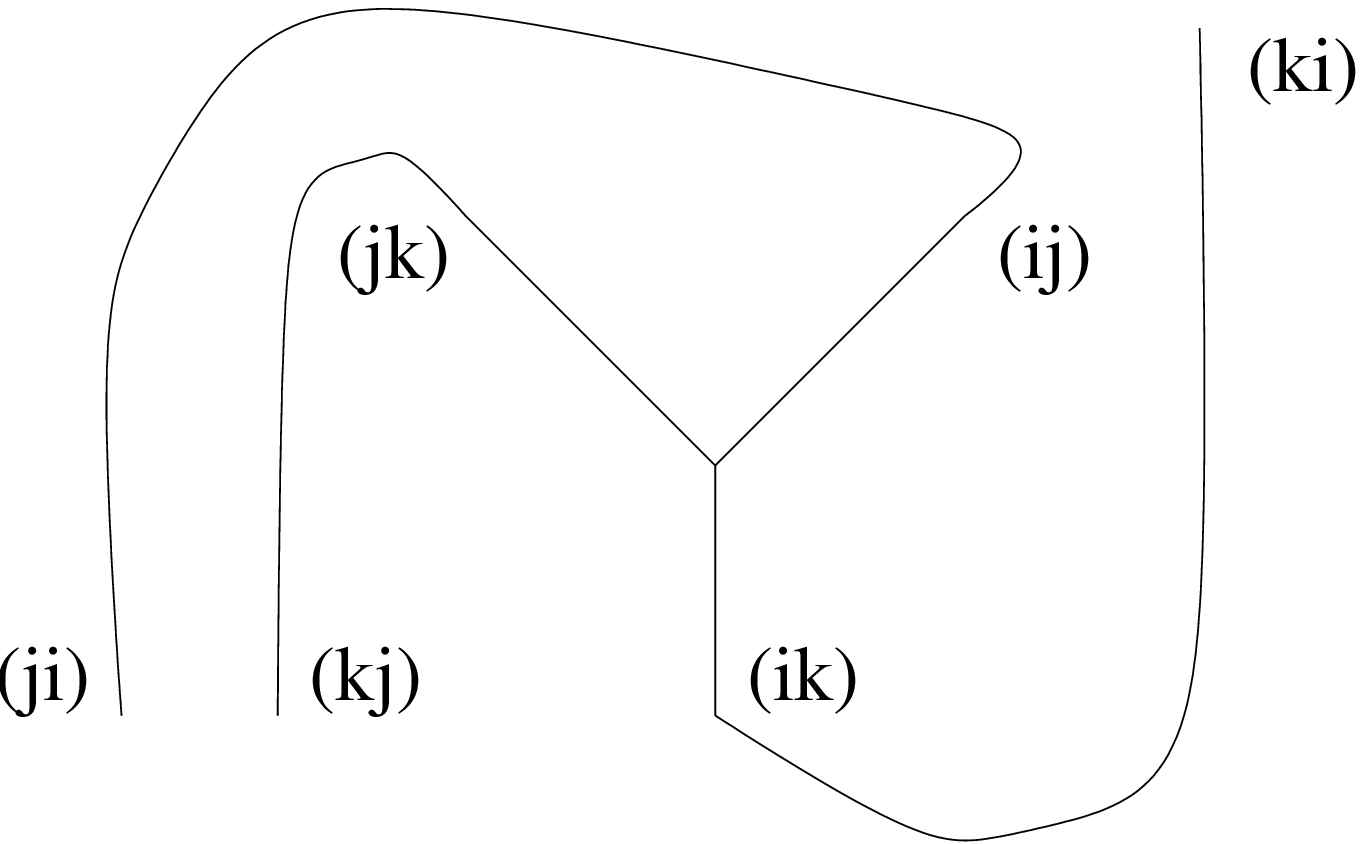}}\hfil}
\vglue 5pt
\centerline{Picture 4. $(f_{ijk}^{\#})^*$ }
\vglue 13pt
\endinsert

Since $Z(M,T,\ell)$ does not depend on the 
chosen ordering on the $4$-simplices themselves, as we have 
explained already at the end of section 1, we can 
restrict our attention to what happens when we 
permute the vertices of a $4$-simplex. 
The symmetric group 
on $5$ elements, $S_5,$ is generated by the transpositions 
interchanging $i$ and $5$ for $1\leq i\leq 4$, so we will 
restrict ourselves to showing that $Z(M,T,\ell)$ is invariant 
under the transposition $\sigma\colon (01234)\mapsto (01432)$. 
The cases of the 
other transpositions are similar. 
Remember the definition of the boundary of $-(01432)$ 
$$-\partial(01432)=-(1432)+(0432)-(0132)+(0142)-(0143).$$ 
\noindent We have to show that 
$\sigma$ induces linear 
isomorphisms $\sigma_{0432}\colon 2H(0234)\to 2H(0432)^*$ and  
$\sigma_{0142}\colon 2H(0124)\to 2H(0142)$ and linear isomorphisms 
$\sigma'_{1432}\colon 2H(1234)^*\to 2H(1432)$, 
$\sigma'_{0143}\colon 2H(0134)^*\to 2H(0143)$ and 
$\sigma'_{0132}\colon 2H(0123)^*\to 2H(0132)$  
such that the following diagram 

\advance\equano by 4
$$\vbox{\offinterlineskip
\halign{\hfil#\hfil&\hfil#\hfil&\hfil#\hfil\cr
$2H(+(01234))$&$\buildrel
 \sigma_{01432}\over{\hbox to 40pt{\rightarrowfill}}$
&$2H(-(01432))$\cr
\noalign{\vskip7pt}
$\Bigg\downarrow$&&$\Bigg\downarrow$\cr
\noalign{\vskip 5pt}
${\bf F}$&$=$&${\bf F}$\cr}}\eqno{(\the\chno.\the\equano)}$$

\noindent commutes. In this diagram $2H(+(01234))$ stands for 
$$2H(1234)^*\otimes 2H(0234)\otimes 2H(0134)^*\otimes 2H(0124)
\otimes 2H(0123)^*$$ 
\noindent and $2H(-(01234))$ for
$$2H(1432)\otimes 2H(0432)^*\otimes 2H(0132)\otimes 2H(0142)^*
\otimes 2H(0143).$$ The linear map $\sigma_{01432}$ is the 
tensor product of the $\sigma_{\sigma(i)\sigma(j)\sigma(k)\sigma(l)}$ 
and the $\sigma'_{\sigma(i)\sigma(j)\sigma(k)\sigma(l)}$ 
composed with some transpositions $P\colon x\otimes y\mapsto y\otimes x$ 
so that the respective factors in the tensor product appear in 
the right order. 
The vertical linear 
maps are the ones defining $Z(+(01234))$ and $Z(-(01432))$ 
respectively by means of the pairing (see section 3). 

First of all we will show what the $\sigma_{\sigma(i)\sigma(j)
\sigma(k)\sigma(l)}$ and 
the $\sigma'_{\sigma(i)\sigma(j)\sigma(k)\sigma(l)}$ are. 
Using the diagrammatic conventions as established by the diagrams in 
the pictures 2, 3 and 4, we are now going to give the definition of $\sigma'_{1432}$. 
Remember that $2H(1234)^*$ is the vector space of $2$-morphisms 
$$\alpha\colon f_{124}(f_{234}\otimes e_{12})\mapsto f_{134}(e_{34}
\otimes f_{123}).$$ 
Diagrammatically we denote such a $2$-morphism by the diagram 
in picture 5. 

\midinsert
\vglue 13pt
\line{\hfil\vbox{\epsfysize=2.5cm\epsfbox{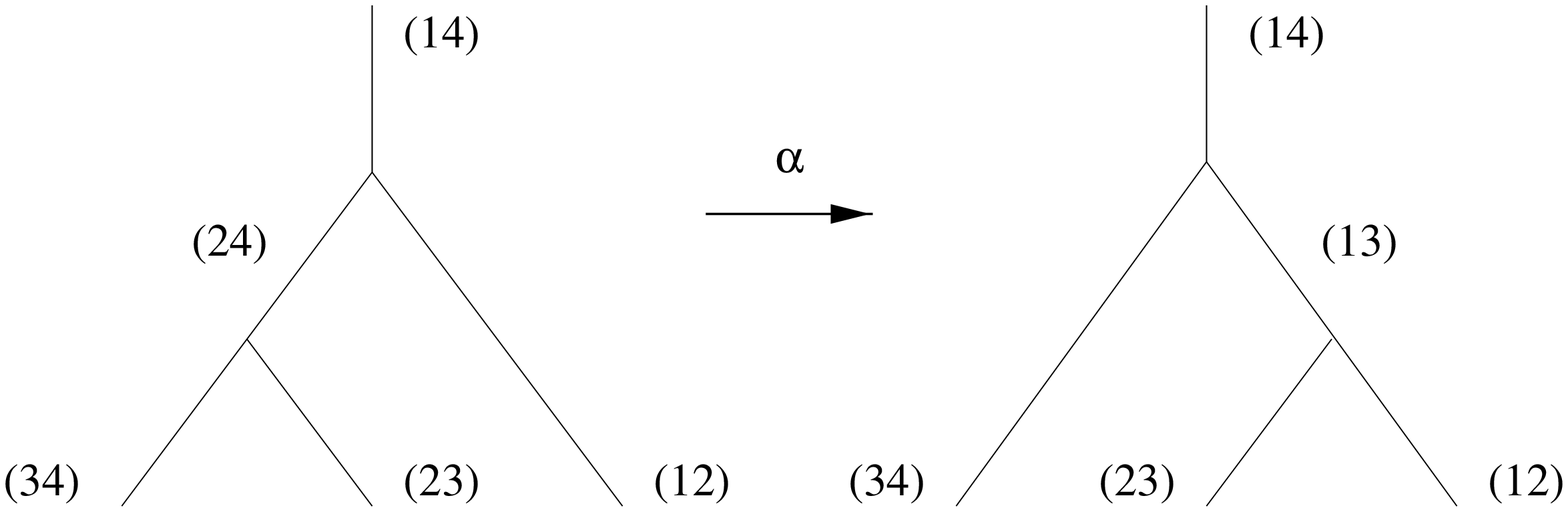}}\hfil}
\vglue 5pt
\centerline{Picture 5.}
\vglue 13pt
\endinsert

\noindent The vector space $2H(1432)$ was defined by 
$$\THom(f_{132}(e_{32}\otimes 
f_{143}),f_{142}(f_{432}\otimes e_{14})).$$
In picture 6 we show how to get a $1$-morphism 
$$e_{12}\mapsto e_{42}\otimes e_{14}\mapsto e_{32}\otimes e_{43}
\otimes e_{14}$$
\noindent from $f_{124}(f_{234}\otimes e_{12})$. 

\midinsert
\vglue 13pt
\line{\hfil\vbox{\epsfysize=3cm\epsfbox{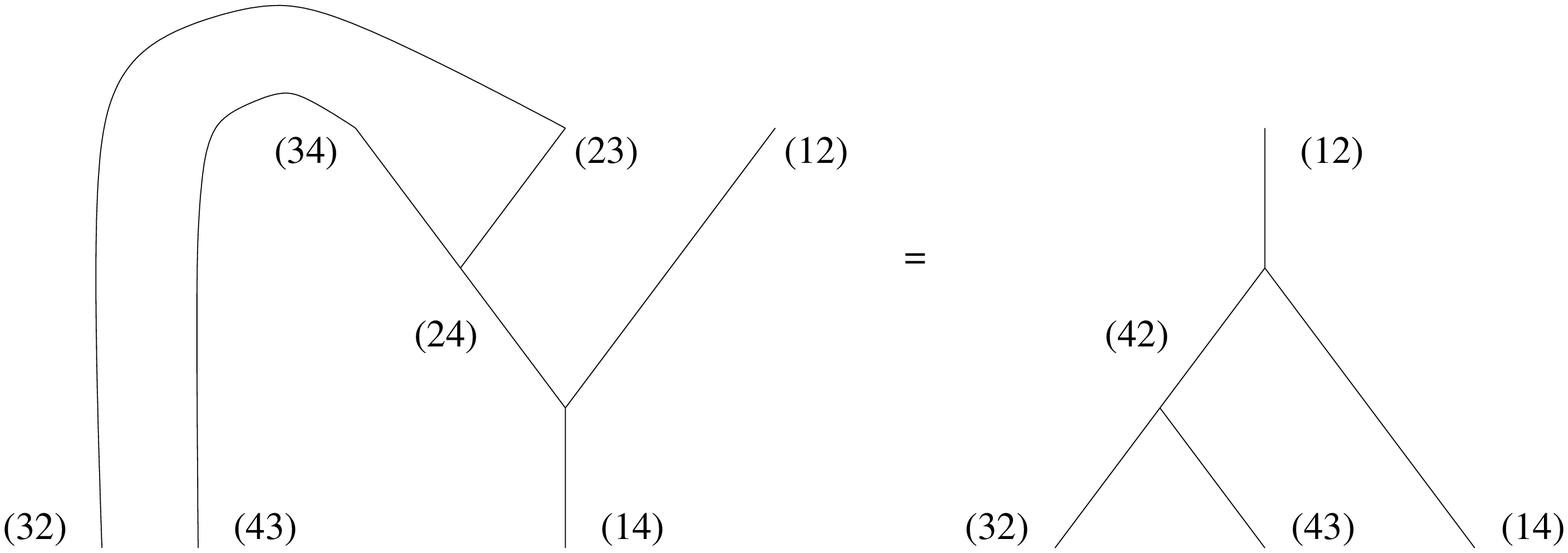}}\hfil}
\vglue 5pt
\centerline{Picture 6. $f_{142}(f_{234}\otimes e_{14})$}
\vglue 13pt
\endinsert 

\noindent Of course the 
cups and caps represent coevaluations $i_{\cdot}$ and evaluations 
$e_{\cdot}$ respectively, defined by the duality on the objects 
in $C$. The way we have drawn our diagrams already shows that we assume 
that the $1$-morphism we obtain in this way is equal to 
$f_{142}(f_{432}\otimes e_{14})$. This is justified by two arguments. 
In the first place we can decompose the $1$-morphism into 
simple $1$-morphisms in $\Hom(e_{12},e_{42}\otimes e_{14})$ and 
$\Hom(e_{42},e_{32}\otimes e_{43})$ by semi-simplicity. Since we 
have shown in lemma 4.1 that $Z(M,T,\ell)$ does not depend on 
the particular choice of simple $1$-morphisms in the various 
isomorphism classes, we may assume that the decomposition gives 
$f_{142}$ and $f_{432}$ actually. 

In the same way we obtain 
a $1$-morphism 
$$e_{12}\mapsto e_{32}\otimes e_{13}\mapsto e_{32}\otimes e_{43}
\otimes e_{14}$$
\noindent from $f_{134}(e_{34}\otimes f_{123})$, which we take 
to be equal to $f_{132}(e_{32}\otimes f_{143})$ (see picture 7). 

\midinsert
\vglue 13pt
\line{\hfil\vbox{\epsfysize=3cm\epsfbox{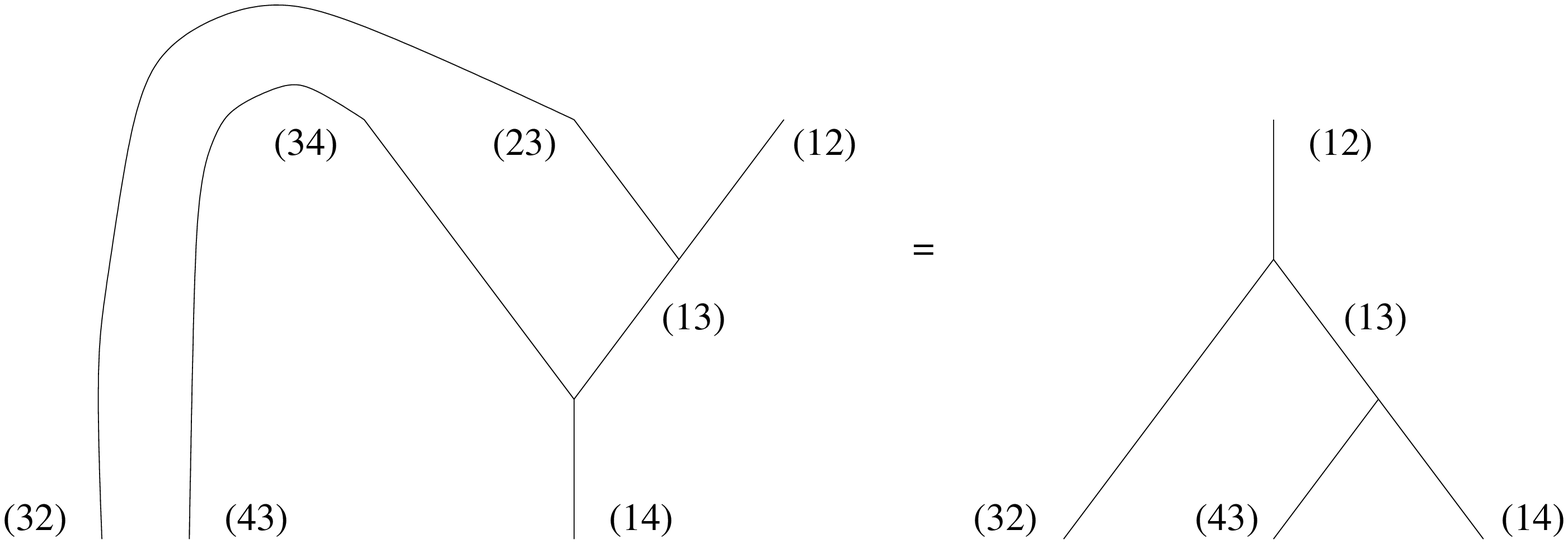}}\hfil}
\vglue 5pt
\centerline{Picture 7. $f_{132}(e_{32}\otimes f_{143})$}
\vglue 13pt
\endinsert

Now by horizontal composition of the identity $2$-morphisms on 
the respective coevaluations and evaluations and $\alpha^{\dag}$ we get the 
image of $\alpha$ under $\sigma'_{1432}$. Explicitly we get 
the $2$-morphism 
$$(i_{e_{32}}\otimes e_{12})(e_{32}\otimes i_{e_{43}}\otimes e_{23}
\otimes e_{12})(e_{32}\otimes e_{43}\otimes f^*_{234}\otimes 
e_{12})$$
$$(e_{32}\otimes e_{43}\otimes e_{24}\otimes i_{e_{42}}
\otimes e_{12})(e_{32}\otimes e_{43}\otimes i^*_{e_{24}}\otimes 
e_{24}\otimes e_{12})(e_{32}\otimes e_{43}\otimes f^*_{124})$$
$$\buildrel T_{e_{42}}\over\Rightarrow 
(i_{e_{32}}\otimes e_{12})(e_{32}\otimes i_{e_{43}}\otimes e_{23}
\otimes e_{12})(e_{32}\otimes e_{43}\otimes f^*_{234}\otimes 
e_{12})(e_{32}\otimes e_{43}\otimes f^*_{124})$$
$$\buildrel \alpha^{\dag}\over\Rightarrow 
(i_{e_{32}}\otimes e_{12})(e_{32}\otimes i_{e_{43}}\otimes e_{23}
\otimes e_{12})(e_{32}\otimes e_{43}\otimes e_{34}\otimes 
f^*_{123})(e_{32}\otimes e_{43}\otimes f^*_{134}).$$
Note that by taking 
the dual of the $1$-morphisms in the diagrams we have to use 
$\alpha^{\dag}$, which corresponds to the fact that $\sigma$ is 
an odd permutation. It is clear that $\sigma'_{1432}$ is a linear 
isomorphism, although it is only an involution up to a 
$2$-isomorphism. But, 
in the end 
when we take the pairing that defines our partition function 
$Z(ijklm)$, these isomorphisms do not harm us because the 
pairing is defined by means of a trace that is invariant under 
conjugation, so we always get the same result anyhow. 
The linear isomorphism $\sigma_{0432}$ is defined in an analogous 
way. As a matter of fact one only has to reverse the morphisms 
in the definition of $\sigma'_{1432}$ and substitute $1$ by $0$. 
Let us now define the linear isomorphism $\sigma'_{0143}
\colon 2H(0134)^*\to 2H(0143)$. Remember that the vector space 
$2H(0134)^*$ was defined by 
$$\THom(f_{014}(f_{134}\otimes e_{01}),f_{034}(e_{34}\otimes f_{013}))$$
\noindent and $2H(0143)$ by 
$$\THom(f_{043}(e_{43}\otimes f_{014}),f_{013}(f_{143}\otimes e_{01})).$$ 
\noindent Picture 8 shows how we map $f_{014}(f_{134}\otimes e_{01})$ 
to $(f_{143}\otimes e_{01})(e_{43}\otimes f^*_{014})$. 

\midinsert
\vglue 13pt
\line{\hfil\vbox{\epsfysize=3cm\epsfbox{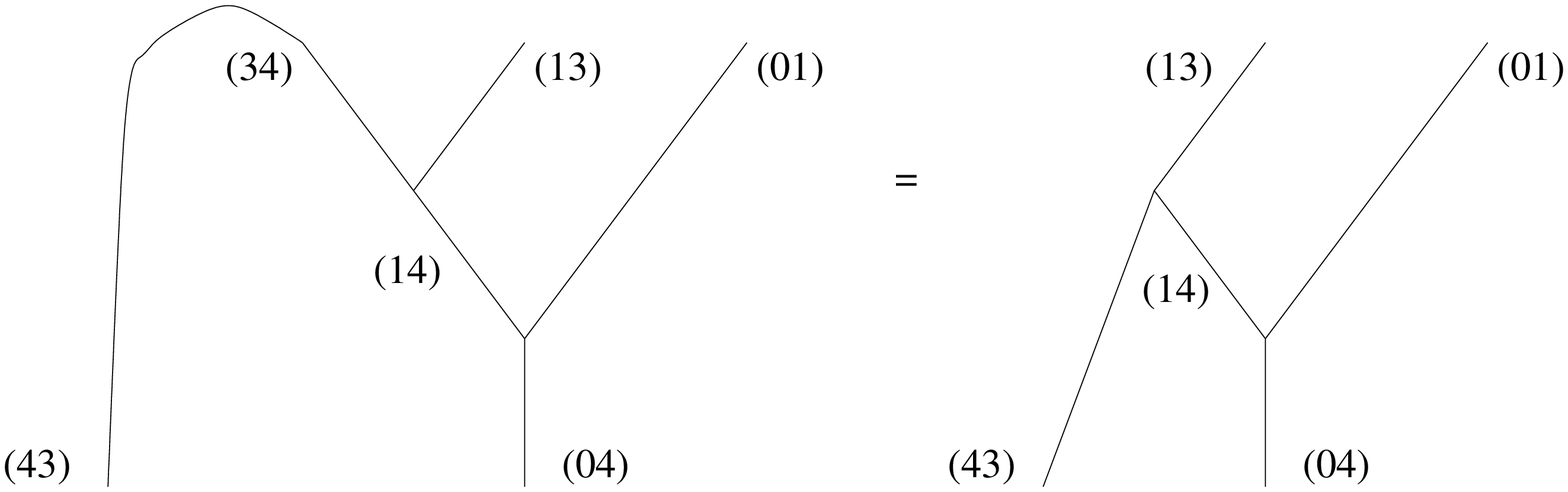}}\hfil}
\vglue 5pt
\centerline{Picture 8. $(f_{143}\otimes e_{01})(e_{43}\otimes f^*_{014})$}
\vglue 13pt
\endinsert

\noindent Picture 9 
shows how we map $f_{034}(e_{34}\otimes f_{013})$ to 
$f^*_{013}f_{043}$. 

\midinsert
\vglue 13pt
\line{\hfil\vbox{\epsfysize=3cm\epsfbox{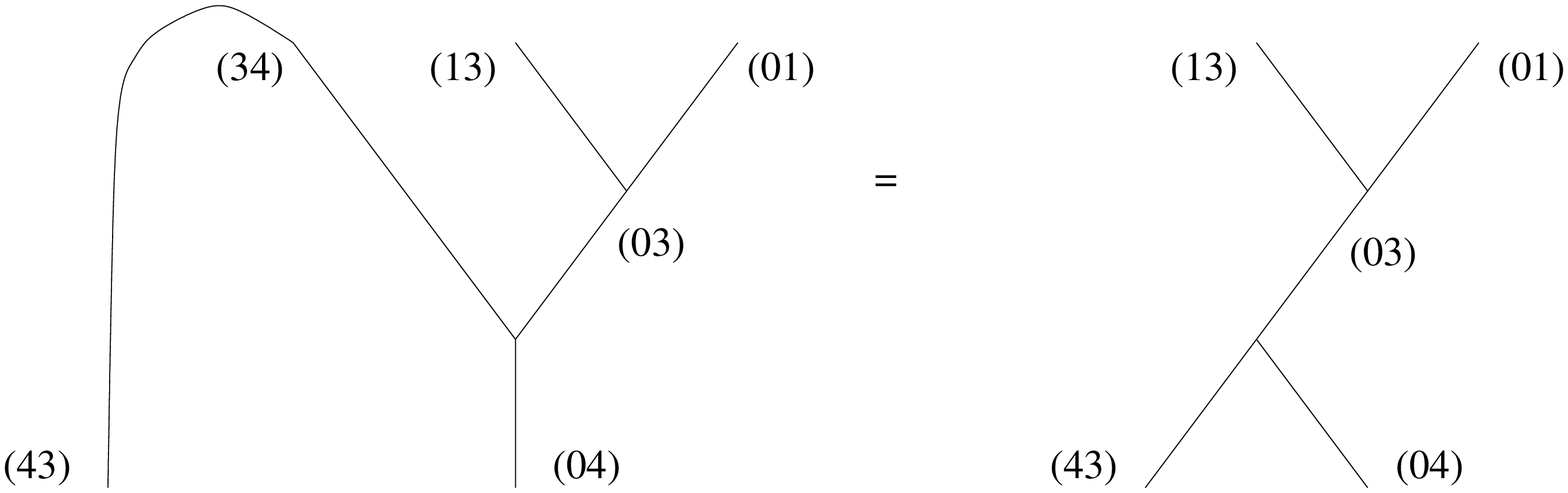}}\hfil}
\vglue 5pt
\centerline{Picture 9. $f^*_{013}f_{043}$}
\vglue 13pt
\endinsert

\noindent So given $\epsilon\in 2H(0134)^*$ we get 
a $2$-morphism 
$$f^*_{013}f_{043}\Rightarrow (f_{143}\otimes e_{01})
(e_{43}\otimes f^*_{014}).$$
\noindent Again we have composed the identities on the coevaluations and 
evaluations with $\epsilon^{\dag}$. Now the image of $\epsilon$ under 
$\sigma'_{0143}$ is defined by the obvious composition 
$$f_{043}(e_{43}\otimes f_{014})\Rightarrow f_{013}f^*_{013}f_{043}(e_{43}\otimes f_{014})$$
$$\Rightarrow f_{013}(f_{143}\otimes 
e_{01})(e_{43}\otimes f^*_{014})(e_{43}\otimes f_{014})  
\Rightarrow f_{013}(f_{143}\otimes e_{01}).$$
\noindent Explicitly we get 
$$f_{043}(e_{43}\otimes f_{014})=(i_{e_{43}}\otimes e_{03})
(e_{43}\otimes f^*_{034})(e_{43}\otimes f_{014})$$
$$\buildrel i_{f_{013}}\over\Rightarrow f_{013}f^*_{013}
(i_{e_{43}}\otimes e_{03})
(e_{43}\otimes f^*_{034})(e_{43}\otimes f_{014})$$
$$\buildrel \otimes^{-1}_{i_{43},f^*_{013}}\over\Rightarrow
 f_{013}(i_{43}\otimes e_{13}\otimes e_{01})(e_{43}\otimes 
e_{34}\otimes f^*_{013})$$
$$(e_{43}\otimes f^*_{034})(e_{43}\otimes 
f_{014})$$
$$\buildrel \epsilon^{\dag}\over\Rightarrow 
f_{013}(i_{43}\otimes e_{13}\otimes e_{01})(e_{43}\otimes 
f^*_{134}\otimes e_{01})$$
$$(e_{43}\otimes f^*_{014})(e_{43}\otimes 
f_{014})$$
$$\buildrel e_{f_{014}}\over\Rightarrow 
f_{013}(i_{43}\otimes e_{13}\otimes e_{01})(e_{43}\otimes f^*_{134}
\otimes e_{01})$$
$$=f_{013}(f_{143}\otimes e_{01}).$$

\noindent Again this is an isomorphism, because $C$ is pivotal. 
Again $\sigma'_{0143}$ is only involutive up to a $2$-isomorphism. 
The other linear isomorphisms $\sigma_{0142}$ and $\sigma'_{0132}$ 
are defined analogously. So now all linear maps in diagram 4.4 are 
defined. Let us prove the commutativity of the diagram. 

\advance\thmno by 1
\proclaim\Lem. With the definitions as above diagram 4.4 is 
commutative.\par
\noindent{\bf Proof}. As I already announced, the proof is 
essentially diagrammatic. I actually worked out all the 
formulas by hand, but these are far too large to fit on ordinary paper. 
Since I have given all the explicit 
isomorphisms and identifications used in the horizontal 
linear isomorphism of our diagram, the reader can work out 
the explicit formulas from them and check that my diagrammatics 
are correct in that way. 

When one works out the image of an element $\alpha\otimes 
\beta\otimes\gamma\otimes\delta\otimes\epsilon\in 
2H(01234)$ under $\sigma_{01432}$, as described above, and 
the right vertical linear map in our diagram explicitly, one 
can read off the diagram in picture 10. 

\pageinsert
\vglue 20pt
\line{\hfil\epsfbox{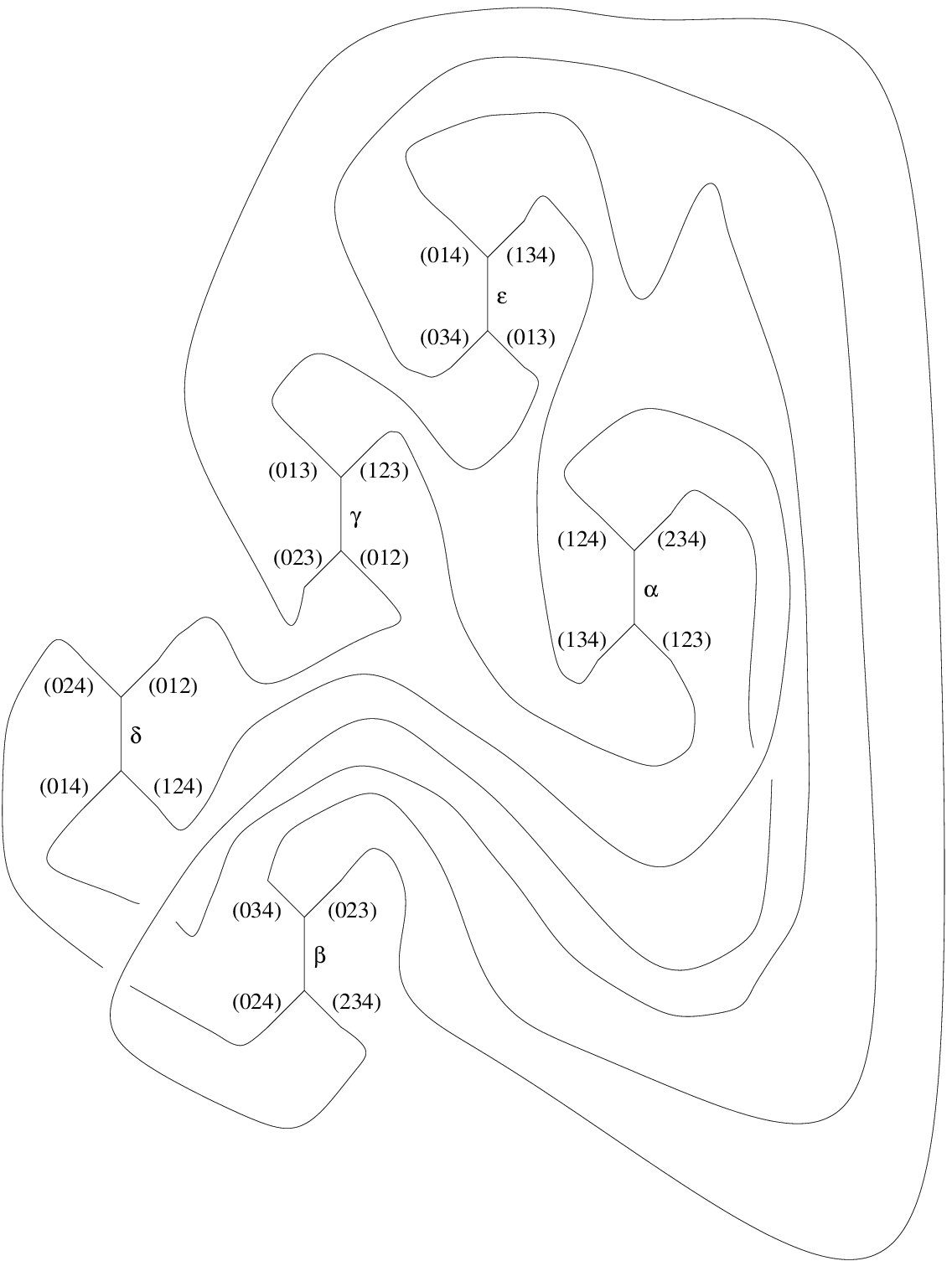}\hfil}
\vglue 5pt
\centerline{Picture 10.}
\vglue 20pt
\endinsert

\pageinsert
\vglue 20pt
\line{\hfil\epsfbox{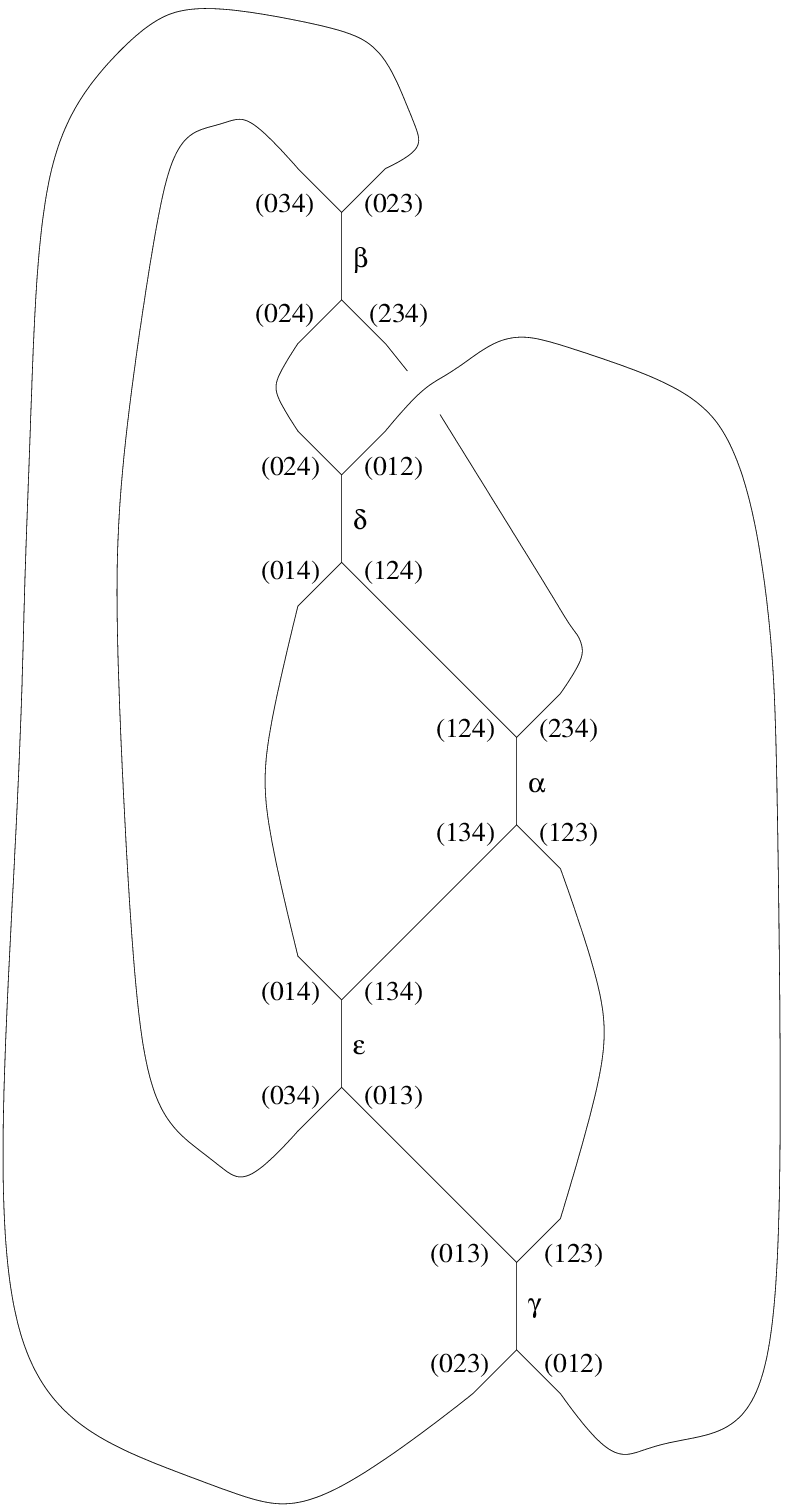}\hfil}
\vglue 5pt
\centerline{Picture 11.}
\vglue 20pt
\endinsert

\noindent We have explained the diagrammatics after our definition 
of $Z(+(ijklm))$ in section 3 and we have explained what kind 
of moves we can apply to them. 
It is now imediately clear that the diagram in picture 
10 can be transformed into the diagram in picture 11. 
We have already explained how these transformations work, but let 
us explain it again in this particular case. For the transformation 
of the diagram in picture 10 to the diagram in picture 11 one only 
has to use the rules for duality on $1$-morphisms (condition 11 in def.~2.3), which 
resemble the same properties for duality on objects in a 
monoidal category, and the 
evident properties of $\otimes_{\cdot,\cdot}$ which precisely 
resemble the properties of a braiding in a monoidal category. So 
far we do not need the condition that $C$ is spherical. But 
in order to transform the diagram in picture 11 into the one 
in picture 1 one definitely needs it. ``Swinging'' around 
the closing strands of our diagram requires both the pivotal 
condition and the spherical 
condition. In order to get the strands labelled by 
$f^*_{023}$ and $f^*_{034}$ to the other side of the diagram 
one needs the pivotal condition. The spherical condition is 
needed because the way in which the caps and cups are 
nested differs in the last two diagrams. This change of nesting 
can only be obtained by changing the left trace functor into 
the right trace functor on some levels of the diagram. All 
these operations of course need to be compatible with the cupping  
and capping, but we have included these compatibilities in 
the axioms of the pivotal and the spherical conditions. 
The proof now finishes with the observation that the last 
diagram is exactly the one describing $Z(01234)$, i.e. 
the left arrow in the diagram of this lemma.\qed  
 
The following theorem is now an immediate result of the previous 
lemmas in this section.

\advance\thmno by 1
\proclaim\Thm. Let $\Phi\colon (M,T,\ell)\to (M',T',\ell')$ be a combinatorial isomorphism 
of labelled triangulated manifolds that is compatible with the 
labellings. Then $Z(M,T,\ell)=Z(M',T',\ell')$.\par

\advance\chno by 1
\thmno=0
\equano=0
\beginsection\the\chno. Invariance under the Pachner moves\par
For the proof of invariance under the Pachner moves we should really 
look at the 4D Pachner moves as equalities between 
series of 3D Pachner moves. This 
idea has been worked out in [11] and can be seen in the figures 
12, 13 and 14. In these figures an arrow should be interpreted as the boundary of the 
$4$-simplex representing a 3D Pachner move. 
The source diagram 
contains the simplicial $3$-complex defining one side of the Pachner 
move and the target diagram contains the simplicial $3$-complex defining 
the other side. For example the arrow labelled by $(01235)$ in the 
first picture represents the $2\rightleftharpoons 3$ 3D Pachner move 
that inserts the edge $(13)$ in the $3$-complex given by the 
two tetrahedra $(0125)$ and $(0235)$ glued over the triangle $(025)$. 

\pageinsert
\vglue 20pt
\line{\hfil\epsfbox{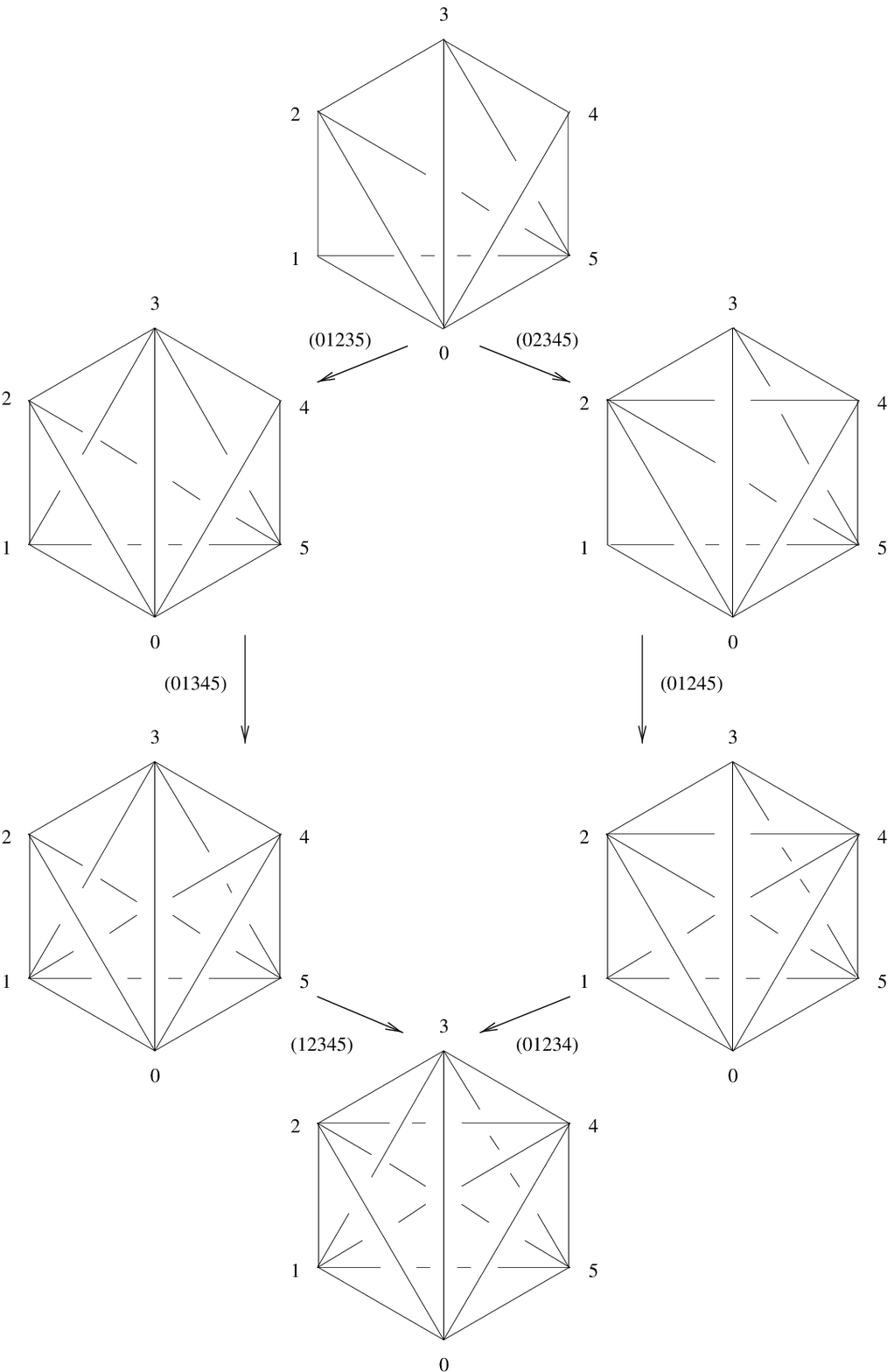}\hfil}
\vglue 5pt
\centerline{Picture 12. Pachner move $3\rightleftharpoons 3$}
\vglue 20pt
\endinsert

\pageinsert
\vglue 20pt
\line{\hfil\epsfbox{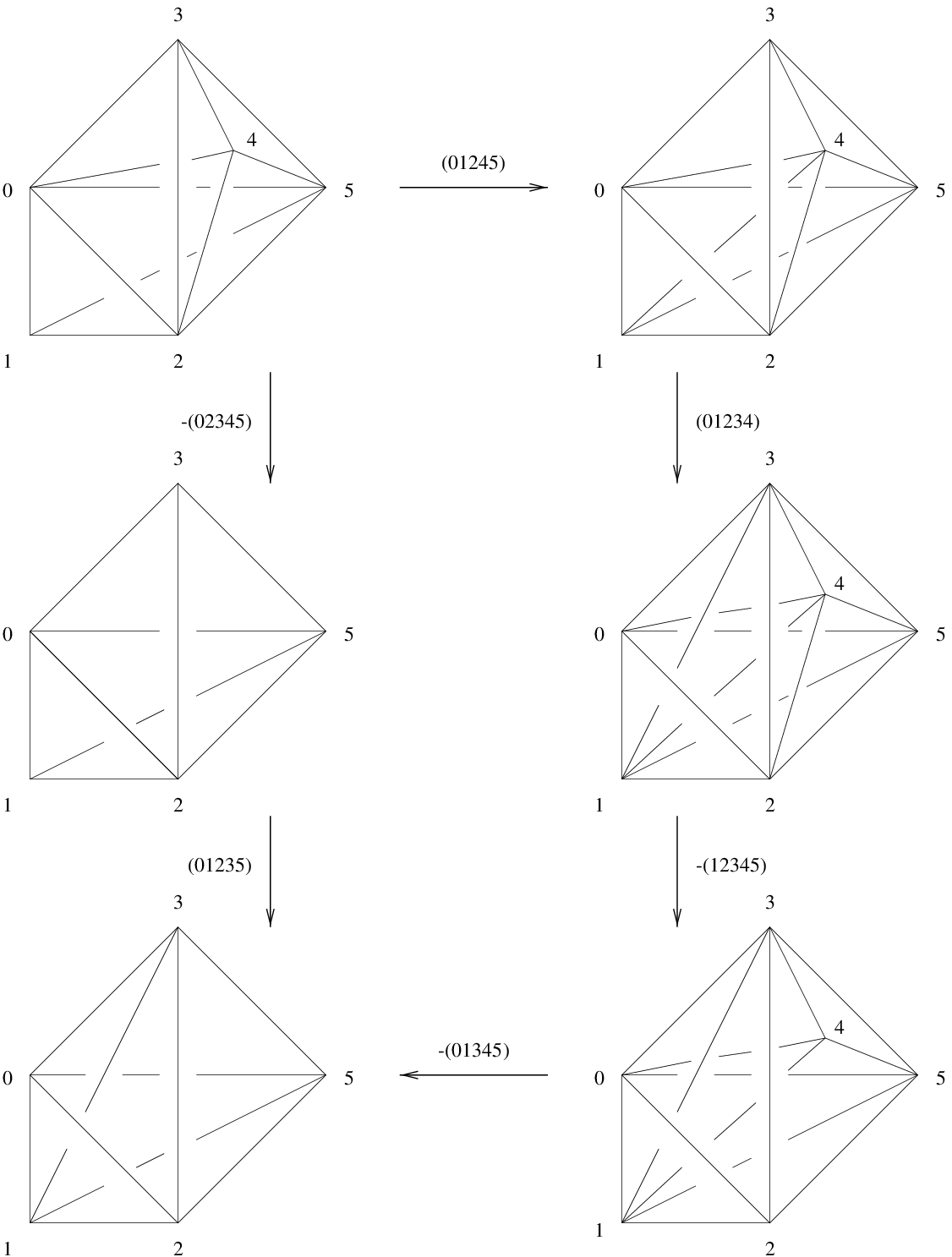}\hfil}
\vglue 5pt
\centerline{Picture 13. Pachner move $2\rightleftharpoons 4$}
\vglue 20pt
\endinsert

\pageinsert
\vglue 20pt
\line{\hfil\epsfbox{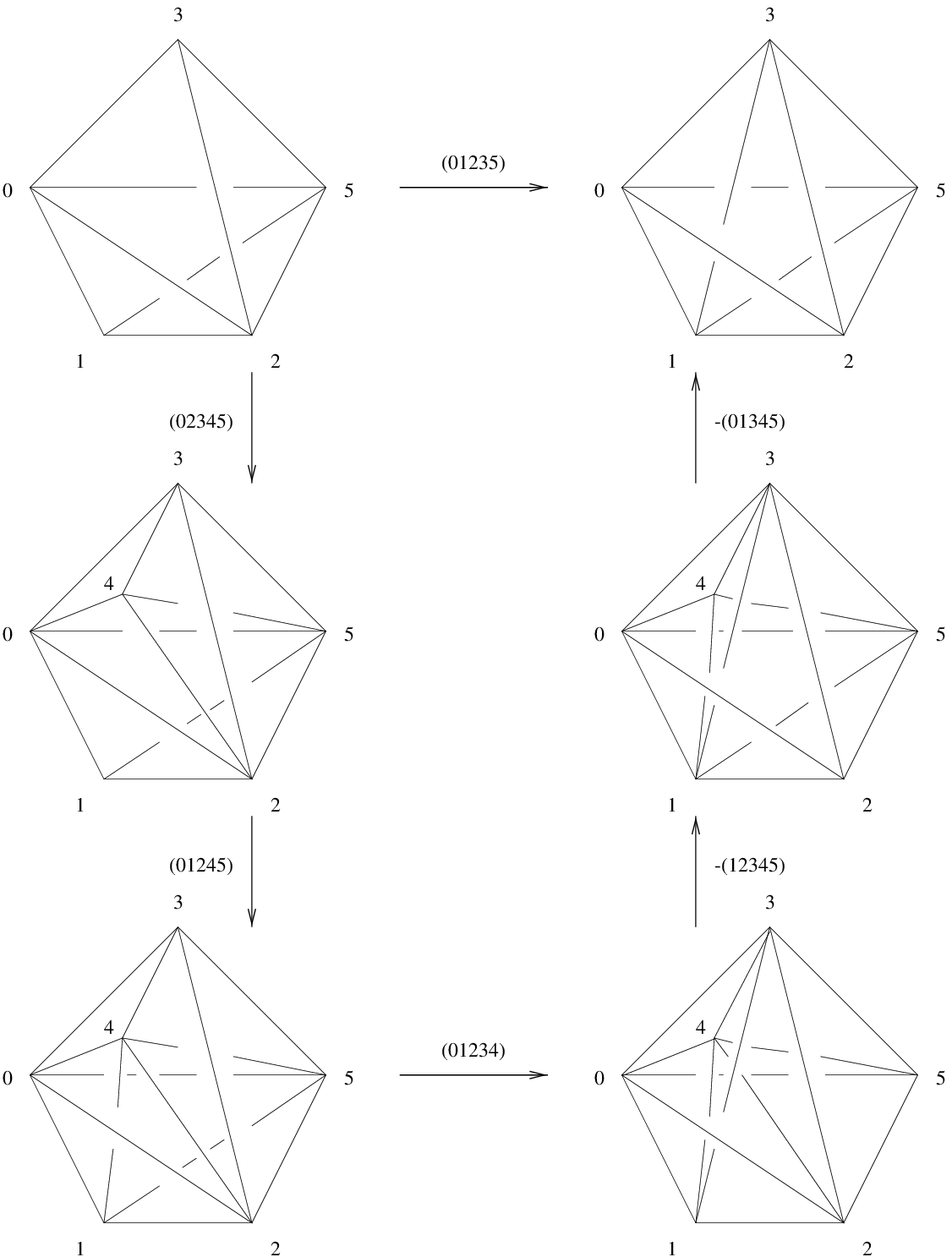}\hfil}
\vglue 5pt
\centerline{Picture 14. Pachner move $1\rightleftharpoons 5$}
\vglue 20pt
\endinsert

These observations about the Pachner moves show that the algebraic 
categorification going from a certain kind of category up to a 
certain kind of $2$-category, as first predicted and sketched in [16],  
goes hand in hand with a geometrical 
kind of categorification. From a very abstract point of view we have 
substituted the identities in the categories which are equivalent  
to the 3D Pachner moves by isomorphisms in the 
$2$-categories which we will prove to satisfy identities equivalent 
to the 4D Pachner 
moves. The $\Phi$ in the crossing lemma 3.1 is the isomorphism that 
substitutes the identity which is equivalent to the 
$2\rightleftharpoons 3$ move of the 3D Pachner moves. Its inverse 
is of course the substitute of the inverse move. The isomorphism 
substituting the $1\rightleftharpoons 4$ move is not so easy to 
describe but will come out of our calculations below. However vague 
these remarks may seem, they describe the deeper reason 
of why everything works as nicely as it does. The notion of 
categorification is really central in this whole setup and causes 
the proofs of invariance under the 4D Pachner moves to become 
almost tautological.

Let $T_1$ and $T_2$ be two triangulations of $M$ that can be obtained 
from one another by one 4D Pachner move. Let $D_1\subseteq T_1$ be the 
simplicial $4$-complex on one side of the Pachner move and $D_2\subseteq 
T_2$ the simplicial $4$-complex on the other. We denote the complement 
of the interior of $D_1$ in $T_1$ by $X$, which by definition is equal to 
the complement of the interior of $D_2$ in $T_2$. Notice that 
$\partial X=\partial D_1=\partial D_2$. Also $D_1\cup D_2$ is the 
boundary of a 5-simplex $(012345)$. Now any labelling of $X\cup\partial 
(012345)$ defines a labelling $\ell_X$ on $X$, a labelling $\ell_1$ 
on $D_1$ and a labelling $\ell_2$ on $D_2$ which are equal on 
intersections. We define $Z(X)$ as the linear map obtained by taking 
the partial trace over all the state spaces of all the tetrahedra in the 
interior of $X$ of the 
tensor product of the partition functions for each labelled $4$-simplex in $X$. 
Define $Z(D_1)$ and $Z(D_2)$ analogously. Then we can decompose our 
state sum $I(M,T_i)=I_C(M,T_i)$ for $i=1,2$ in the following way:
$$I(M,T_i)=K^{-v_X}\tr\bigg(\sum_{l_X}Z(X)\otimes\Big(K^{-v'_i}\sum_{\ell'_i}Z(D'_i)
\prod_{e_{D'_i}}\qd(\ell(e_{D'_i}))^{-1}\prod_{f_{D'_i}}\qd(\ell(f_{D'_i})
)\Big)\bigg)$$
$$\prod_{e_X}\qd(\ell(e_X))^{-1}\prod_{f_X}\qd(\ell(f_X)).$$
\noindent The number $v_X$ is the number of vertices of $X$, $v'_i$ 
is the number of vertices internal to $D_i$ for $i=1,2$. The first 
summation is over all labellings on $X$, the second is over all 
labellings on $D_i$ fixed 
on $\partial D_i$ but ranging over all the simple objects in  
$\cal E$ and all the simple $1$-morphisms in all the ${\cal F}_{ijk}$ 
for all the edges and faces internal to $D_i$. The trace is of course  
the trace over all the state spaces of all the tetrahedra in $\partial X$. 
Therefore proving invariance of our state sum under a 4D Pachner move means showing that the 
following identity holds
$$K^{-v'_1}\sum_{\ell'_1}Z(D'_1)\prod_{e_{D'_1}}\qd(\ell(e_{D'_1}))^{-1}\prod_{f_{D'_1}}\qd(\ell(f_{D'_1}))$$
$$=K^{-v'_2}\sum_{\ell'_2}Z(D'_2)\prod_{e_{D'_2}}\qd(\ell(e_{D'_2}))^{-1}\prod_{f_{D'_2}}\qd(\ell(f_{D'_2})).$$
The lemmas in the rest of this section prove this identity for all 
the 4D Pachner moves.  

The equation proving invariance under the 
$3\rightleftharpoons 3$ move is the analogue of the Biederharn-Elliot equation. In the following lemma $P$ is the transposition $P\colon x\otimes y\mapsto 
y\otimes x.$

\advance\thmno by 1
\proclaim\Lem ($3\rightleftharpoons 3$).
$$\sum_{f_{135}}\qd(f_{135})(1\otimes Z(+(01235)))(P\otimes 1\otimes 1)
(1\otimes Z(+(01345))\otimes 1)$$
$$(P\otimes 1\otimes 1\otimes 1)
(Z(+(12345))\otimes 1\otimes 1\otimes 1)(1\otimes 1\otimes P\otimes 1
\otimes 1)$$
\smallskip
$$=\sum_{f_{024}}\qd(f_{024})(Z(+(02345))(1\otimes 1\otimes P)
(1\otimes Z(+(01245))\otimes 1)$$
$$(1\otimes 1\otimes 1\otimes P)
(1\otimes 1\otimes 1\otimes Z(+(01234))).$$
\par
\noindent{\bf Proof}. Just write down the hexagon of which the 
left-hand side is 
$$\bigoplus_{f_{035},f_{025}}2H(0345)\otimes 2H(0235)
\otimes 2H(0125)$$
$$\Big\downarrow\ 1\otimes\Phi_{01235}$$
$$\bigoplus_{e_{13},f_{035},f_{135},\atop
f_{123},f_{013}}2H(0345)\otimes 2H(1235)\otimes 2H(0135)\otimes 
2H(0123)$$
$$\Big\downarrow\ (P\otimes 1\otimes 1)(1\otimes \Phi_{01345}\otimes 1)$$
$$\bigoplus_{e_{13},e_{14},f_{135},f_{123},\atop
f_{013},f_{145},f_{134},f_{014}}2H(1235)\otimes 2H(1345)\otimes 
2H(0145)\otimes 2H(0134)\otimes 2H(0123)$$
$$\Big\downarrow\ (P\otimes 1\otimes 1\otimes 1)(\Phi_{12345}\otimes 1
\otimes 1\otimes 1)(1\otimes 1\otimes P\otimes 1\otimes 1)$$
$$\bigoplus_{{e_{13},e_{14},e_{24},\atop
f_{245},f_{234},f_{145},f_{124},}\atop
f_{014},f_{134},f_{123},f_{013}}2H(2345)\otimes 2H(1245)\otimes 
2H(0145)\otimes 2H(1234)\otimes 2H(0134)\otimes 2H(0123),$$ 
\noindent and the right-hand side
$$\bigoplus_{f_{035},f_{025}}2H(0345)\otimes 2H(0235)
\otimes 2H(0125)$$
$$\Big\downarrow\ \Phi_{02345}\otimes 1$$
$$\bigoplus_{e_{24},f_{025},f_{245}\atop
f_{024},f_{234}}2H(2345)\otimes 2H(0245)\otimes 2H(0234)\otimes 2H(0125)$$
$$\Big\downarrow\ (1\otimes 1\otimes P)(1\otimes \Phi_{01245}\otimes 1)$$
$$\bigoplus_{e_{24},e_{14},f_{245},f_{024},\atop
f_{234},f_{145},f_{124},f_{014}}2H(2345)\otimes 2H(1245)\otimes 
2H(0145)\otimes 2H(0124)\otimes 2H(0234)$$
$$\Big\downarrow\ (1\otimes 1\otimes 1\otimes P)(1\otimes 1\otimes 1
\otimes \Phi_{01234})$$
$$\bigoplus_{{e_{13},e_{14},e_{24},\atop
f_{245},f_{234},f_{145},f_{124},}\atop
f_{014},f_{134},f_{123},f_{013}}2H(2345)\otimes 2H(1245)\otimes 
2H(0145)\otimes 2H(1234)\otimes 2H(0134)\otimes 2H(0123).$$  
Applying the crossing lemma (lemma 3.1) six times, i.e. for each 
$\Phi_{ijklm}$ separately, shows that the hexagon is commutative. 
The result now follows from restriction to the summands which 
appear in the lemma.\qed

In the next lemma we prove the analogue of the orthogonality 
equation.
\advance\thmno by 1
\proclaim\Lem ({\rm Orthogonality}).  
$$\qd(f_{024})\sum_{\scriptstyle e_{13}, f_{013},\atop 
f_{123}, f_{134}}Z(+(01234))Z(-(01234))\qd(e_{13})^{-1}\qd(f_{013})
\qd(f_{123})\qd(f_{134})$$
$$=\id_{2H(0234)\otimes 2H(0124)}.$$
\smallskip
$$\qd(e_{13})^{-1}\qd(f_{013})\qd(f_{123})\qd(f_{134})\sum_{f_{024}}
Z(-(01234))Z(+(01234))\qd(f_{024})$$
$$=\id_{2H(1234)\otimes 2H(0134)
\otimes 2H(0123)}.$$\par
\noindent{\bf Proof}. This follows from the formulas in the crossing 
lemma (lemma 3.1) for $\Phi$ and $\Phi^{-1}$.\qed

Now the other two Pachner moves follow from the $3\rightleftharpoons 3$ 
lemma and the orthogonality lemma. 

\advance\thmno by 1
\proclaim\Lem ($2\rightleftharpoons 4$). 
$$(Z(-(02345))\otimes 1)(1\otimes 
Z(+(01235))$$
$$=\sum_{e_{14},f_{014},f_{124},f_{145},f_{134}}
\qd(e_{14})^{-1}\qd(f_{014})\qd(f_{124})\qd(f_{145})\qd(f_{134})$$
$$(1\otimes 1\otimes P)(1\otimes Z(+(01245))
\otimes 1)(1\otimes 1\otimes 1\otimes P)$$
$$(1\otimes 1\otimes 1\otimes Z(+(01234)))(1\otimes 1
\otimes P\otimes 1\otimes 1)(Z(-(12345))\otimes 1\otimes 1
\otimes 1)$$
$$(P\otimes 1\otimes 1\otimes 1)(1
\otimes Z(-(01345))\otimes 1)(P\otimes 1\otimes 1).$$\par
\noindent{\bf Proof}. Note that the expressions in this lemma do not 
relate to picture 13 just as directly as the expressions in the lemmas 
about the invariance under the other two 4D Pachner moves do to their 
respective pictures. What I mean is 
that $2H(0345)$ is in the targets of $Z(-(02345))$ and $Z(-(01345))$ 
rather than in their sources. Since the partial trace over the state spaces 
in $\partial X$ in the decomposition of $I(M,T_i)$ is the composition of the 
coevaluation on these state spaces with 
$$\bigg(\sum_{l_X}Z(X)\otimes\Big(K^{-v'_i}\sum_{\ell'_i}Z(D'_i)
\prod_{e_{D'_i}}\qd(\ell(e_{D'_i}))^{-1}\prod_{f_{D'_i}}\qd(\ell(f_{D'_i})
)\Big)\bigg)$$
\noindent and with the evaluation on $2H(0345)$ we can just as well multiply 
both sides of the equation in this lemma by 
$$(\ev_{2H(0345)}\otimes 1\otimes 1\otimes 1).$$
\noindent That is what I had actually done in an earlier version of 
this paper, because than the expressions are exactly the ones one can read 
off from the diagram. However in that case the expressions do not exactly 
correspond to the decomposition of $I(M,T_i)$ into partial traces and so 
I have decided to change it in this version to avoid confusion. The essence 
of the lemma remains the same.

Let us prove the lemma now. 
Multiply each side of the 
$3\rightleftharpoons 3$ equation by 
$$\qd(e_{24})^{-1}\qd(e_{14})^{-1}\qd(f_{234})\qd(f_{245})\qd(f_{014})
\qd(f_{124})$$
$$\qd(f_{145})\qd(f_{134})\qd(f_{035})(Z(-(02345))\otimes 1)$$
\noindent on the left and multiply by  
$$(1\otimes 1\otimes P\otimes 1\otimes 1)(Z(-(12345))
\otimes 1\otimes 1\otimes 1)(P\otimes 1\otimes 1\otimes 1)$$
$$(1\otimes Z(-(01345))\otimes 1)(P\otimes 1\otimes 1)$$
\noindent on the right and sum over all the edges, i.e. simple 
objects, and all the faces, i.e. simple $1$-morphisms, involved.
Using the orthogonality lemma once on the left-hand side and twice 
on the right-hand side we get the $2\rightleftharpoons 4$ equation.
\qed

\advance\thmno by 1
\proclaim\Lem ($1\rightleftharpoons 5$).
$$Z(+(01235))$$
$$=K^{-1}\sum_{{e_{04},e_{14},e_{24},e_{34},e_{45},\atop
f_{014},f_{024},f_{034},f_{045},f_{124},}\atop 
f_{134},f_{145},f_{234},f_{245},f_{345}}
\tr_1\big[(Z(+(02345))\otimes 1)(1\otimes 1\otimes P)(1\otimes Z(+(01245))
\otimes 1)$$
$$(1\otimes 1\otimes 1\otimes P)(1\otimes 1\otimes 1
\otimes Z(+(01234)))(1\otimes 1\otimes P\otimes 
1\otimes 1)(Z(-(12345))\otimes 1\otimes 1\otimes 1)$$
$$(P\otimes 1\otimes 1\otimes 1)(1\otimes Z(-(01345))\otimes 1)
(P\otimes 1\otimes 1)\big]\prod_{i<j\atop i=4\vee j=4}\qd(e_{ij})^{-1}\prod_{i<j<k\atop j=4\vee k=4}\qd(f_{ijk}).$$\par
\noindent{\bf Proof}. Multiply each side of the 
the $3\rightleftharpoons 3$ equation on the right by   
$$K^{-1}\prod_{i<j\atop i=4\vee j=4}\qd(e_{ij})^{-1}\prod_{{i<j<k\atop 
j=4\vee k=4}\atop (ijk)\neq (024)}\qd(f_{ijk})
(1\otimes 1\otimes P\otimes 1\otimes 1)$$
$$(Z(-(12345))\otimes 1\otimes 1\otimes 1)(P\otimes 1\otimes 1\otimes 
1)(1\otimes Z(-(01345))\otimes 1)(P\otimes 1\otimes 1),$$ 
take the trace on the first factor and sum over all the edges and 
faces involved. 

The right-hand side is now equal to the right-hand side of the 
$1\rightleftharpoons 5$ equation.

Using the orthogonality lemma the left-hand side becomes
$${K^{-1}\over \qd(f_{035})}\sum_{e_{04},e_{34},e_{45},\atop
f_{034},f_{045},f_{345}}\qd(e_{04})^{-1}\qd(e_{34})^{-1}\qd(e_{45})^{-1}
\qd(f_{034})$$
$$\qd(f_{045})\qd(f_{345})\tr_1(1\otimes Z(+(01235))).$$
Now use the identity
$$\tr_1(1\otimes Z(+(01235)))=Z(+(01235))\dim(\THom(f_{045}
(e_{45}\otimes f_{034}),f_{035}(f_{345}\otimes e_{03})))$$ 
\noindent and the identity 
$$\sum_{e_{04},f_{034},f_{045}}
\qd(e_{04})^{-1}\qd(f_{034})\qd(f_{045})\dim(\THom(f_{045}
(e_{45}\otimes f_{034}),f_{035}(f_{345}\otimes e_{03})))$$
$$=\qd(e_{35})^{-1}\qd(f_{035})\qd(f_{345}),$$
which will follow from lemma 5.5.
Finally we get
$${K^{-1}\over\qd(f_{035})}\sum_{e_{34},e_{45},\atop
f_{345}}\qd(e_{34})^{-1}\qd(e_{35})^{-1}\qd(e_{45})^{-1}
\qd^2(f_{345})$$
$$\qd(f_{035})Z(+(01235))=Z(+(01235)).$$
\noindent This equality will follow from lemma 5.6.\qed

\advance\thmno by 1
\proclaim\Lem. With the same notation as everywhere in this section 
we have 
$$\sum_{e_{04},f_{034},f_{045}}
\qd(e_{04})^{-1}\qd(f_{034})\qd(f_{045})\dim(\THom(f_{045}
(e_{45}\otimes f_{034}),f_{035}(f_{345}\otimes e_{03})))$$
$$=\qd(e_{35})^{-1}\qd(f_{035})\qd(f_{345}).$$\par
\noindent{\bf Proof}. Consider $f_{035}(f_{345}\otimes e_{03})$. 
Its quantum dimension is equal to 
$$\qd(f_{035})\qd(e_{35})^{-1}\qd(f_{345})$$ 
by lemma 2.25. Now use semi-simplicity to write 
$$1_{f_{035}(f_{345}\otimes e_{03})}=\sum_{e_{04}\atop f_{034},f_{045}}
\alpha_{0345}\cdot{\overline \alpha}_{0345},$$ 
where 
$$\alpha_{0345}\in\THom(f_{035}(f_{345}\otimes e_{03}),
f_{045}(e_{45}\otimes f_{034}))$$ 
\noindent and
$${\overline\alpha}_{0345}\in\THom(f_{045}
(e_{45}\otimes f_{034}),f_{035}(f_{345}\otimes e_{03})).$$
We can always take the $\alpha$'s to be the projections and 
${\overline\alpha}$ the inclusions, so 
$${\overline\alpha}_{0345}\cdot\alpha_{0345}=1_{f_{045}
(e_{45}\otimes f_{034})}.$$
\noindent Therefore we get  
$$\qd(f_{035}(f_{345}\otimes e_{03}))$$
$$=\sum_{e_{04}\atop f_{034},f_{045}}
\langle \alpha_{0345},{\overline \alpha}_{0345}\rangle $$
$$=\sum_{e_{04}\atop f_{034},f_{045}}\langle {\overline \alpha}_{0345}
,\alpha_{0345}\rangle $$
$$=\sum_{e_{04},\atop f_{034},f_{045}}\qd(f_{045}
(e_{45}\otimes f_{034}))\dim(\THom(f_{045}
(e_{45}\otimes f_{034}),f_{035}(f_{345}\otimes e_{03})))$$
$$=\sum_{e_{04},\atop f_{034},f_{045}}
\qd(e_{04})^{-1}\qd(f_{034})\qd(f_{045})\dim(\THom(f_{045}
(e_{45}\otimes f_{034}),f_{035}(f_{345}\otimes e_{03}))).$$\qed

\advance\thmno by 1
\proclaim\Lem. Let $A$ be any simple object in ${\cal E}$ (see def.~2.20). 
Then we obtain the following expression for $K$:
$$K=\sum_{B,C,\atop
f_{A,B\otimes C}}\qd(A)^{-1}\qd(B)^{-1}\qd(C)^{-1}\qd^2(f_{A,B\otimes C}),$$ 
where we sum over all the objects $B,C\in {\cal E}$ and all 
$1$-morphisms 
$f_{A,B\otimes C}\in{\cal F}_{A,B\otimes C}$, the finite basis of 
non-isomorphic $1$-morphisms in $\Hom(A,B\otimes C)$. Note that this implies 
that this expression is independent of the choice of $A\in{\cal E}$.\par
\noindent{\bf Proof}. First let us rewrite the expression in this lemma. 
$$\sum_{B,C,\atop
f_{A,B\otimes C}}\qd(A)^{-1}\qd(B)^{-1}\qd(C)^{-1}\qd^2(f_{A,B\otimes C})$$
$$=\sum_{B}\qd(A)^{-1}\qd(B)^{-1}\sum_{C,f_{A,B\otimes C}}\qd(C)^{-1}
\qd(f_{A,B\otimes C})^2$$
$$=\sum_{B}\qd(A^*)^{-1}\qd(B)^{-1}\sum_{C,f_{A,B\otimes C}}\qd(C^*)^{-1}
\qd(f_{C^*,A^*\otimes B})^2.$$
\noindent The third equality is justified by the results in section 4 where 
we show that there is a bijection between the isomorphism classes of 
simple $1$-morphisms $f_{A,B\otimes C}\colon A\to B\otimes C$ and the 
isomorphism classes of simple $1$-morphisms $f_{C^*,A^*\otimes B}
\colon C^*\to A^*\otimes B$. The definition of the bijection is such that 
the quantum dimensions of corresponding $1$-morphisms are equal. 

The lemma now follows if we can show the following identity:
$$\sum_{X,f_{X,Y\otimes Z}}\qd(X)^{-1}\qd(f_{X,Y\otimes Z})^2=
\qd(Y)\qd(Z).$$
\noindent Note first of all that $\qd(Y)\qd(Z)$ is the quantum dimension 
of $Y\otimes Z$ by the pivotal and spherical conditions, which is by 
definition equal to $\langle 1_{1_{Y\otimes Z}},1_{1_{Y\otimes Z}}\rangle$. 
By horizontal semi-simplicity we can decompose $1_{Y\otimes Z}$ by 
$$1_{Y\otimes Z}=\sum_{X,\atop
f_X,{\overline f}_X}f_X{\overline f}_{X},$$
\noindent where $f_X\colon Y\otimes Z\to X$ and ${\overline f}_{X}
\colon X\to Y\otimes Z$ are the projections onto $X$ and the inclusions 
of $X$ into $Y\otimes Z$ respectively. We have to sum over all $f_X$ and ${\overline f}_X$ 
because the multiplicity of $X$ in $Y\otimes Z$ may be greater than 1 in 
general. By vertical semi-simplicity we can decompose $1_{f_X}$ and  
$1_{{\overline f}_X}$ for each $f_X$ and each ${\overline f}_X$ respectively. 
This can be written as
$$1_{f_X}=\sum_{i,\atop
{\alpha^i_X},{\overline\alpha}^i_X}\alpha^i_X\cdot{\overline\alpha}^i_X,$$ 
\noindent and
$$1_{{\overline f}_X}=\sum_{j,\atop
{\beta^j_X},{\overline\beta}^j_X}\beta^j_X\cdot{\overline\beta}^j_X.$$ 
\noindent In the first decomposition each $\alpha^i_X\colon f_X\to f^i_X$ 
is a projection on a simple $1$-morphism $f^i_X$ and each ${\overline\alpha}
^i_X\colon f^i_X\to f_X$ is an inclusion. The same holds for the second 
decomposition where each $\beta^j_X\colon{\overline f}_X\to {\overline f}^j_X$
is a projection and each ${\overline\beta}^j_X\colon {\overline f}^j_X\to 
{\overline f}_X$ is an inclusion. We can assume that the $\alpha^i_X\cdot 
{\overline\alpha}^i_X$ and the 
$\beta^j_X\cdot{\overline\beta}^j_X$ form ``dual'' bases, in the sense that 
$$(\beta^j_X\circ\alpha^i_X)\cdot({\overline\beta}^j_X\circ
{\overline\alpha}^i_X)$$
$$=(\beta^j_X\cdot{\overline\beta}^j_X)\circ(\alpha^i_X\cdot 
{\overline\alpha}^i_X)=\delta^j_i 1_{1_X},$$
\noindent where $\delta^j_i$ is the Kronecker delta. We can make this 
assumption because there is the non-degenerate pairing 
$$\THom(f_X,f_X)\otimes\THom({\overline f}_X,{\overline f}_X)\to {\bf F}$$
\noindent defined by 
$$\alpha\otimes\beta\mapsto\langle\alpha\circ\beta, 1_{{\overline f}_X f_X}
\rangle.$$
\noindent This last expression is just ``the trace'' of $\alpha\circ\beta$, 
where ``trace'' means the vertical trace-like map followed by the trace 
functor and cupping and capping as everywhere in this paper. 

Now we can write 
$$1_{1_{Y\otimes Z}}=\sum(\alpha^i_X\cdot {\overline\alpha}^i_X)\circ 
(\beta^j_X\cdot {\overline\beta}^j_X)$$
$$=\sum(\alpha^i_X\circ \beta^j_X)\cdot ({\overline\alpha}^i_X\circ 
{\overline\beta}^j_X).$$

As a result we obtain the following equalities.
$$\langle 1_{1_{Y\otimes Z}},1_{1_{Y\otimes Z}}\rangle=\sum\langle\alpha^i_X
\circ\beta^j_X, {\overline\alpha}^i_X\circ{\overline\beta}^j_X\rangle$$
$$=\sum\langle\beta^j_X\circ\alpha^i_X,{\overline\beta}^j_X\circ{\overline
\alpha}^i_X\rangle$$
$$=\sum\langle\beta^i_X\circ\alpha^i_X,{\overline\beta}^i_X\circ{\overline
\alpha}^i_X\rangle$$    
$$=\sum\langle{\overline\beta}^i_X\circ{\overline\alpha}^i_X,\beta^i_X\circ
\alpha^i_X\rangle.$$
\noindent These equalities all follow from the pivotal and spherical 
conditions and the fact that the $2$-morphisms involved are only inclusions 
and projections. Now it is not hard to see that the last expression is 
equal to 
$$\sum \dim(V((f^i_X)^*))^2\langle 1_X,1_X\rangle$$
$$=\sum\qd(X)\dim(V((f^i_X)^*))^2=\sum\qd(X)^{-1}\qd((f^i_X)^*)^2.$$ 
\noindent Here $V((f^i_X)^*)$ is the vector space that appears in the 
equality 
$$(f^i_X)^* f^i_X=V((f^i_X)^*)1_X.$$ 
It is easy to see that the dimension of this vector space is equal to 
the dimension of $V({\overline f}^i_X)$, with  
$${\overline f}^i_X ({\overline f}^i_X)^*=V({\overline f}^i_X)1_X,$$ 
because $f^i_X$ and ${\overline f}^i_X$ are ``dual'' to each other.

The result now follows from the observation that the 
$(f^i_X)^*$ form a basis of $\Hom(X,Y\otimes Z)$.\qed

Thus we can conclude this section by the following theorem, summarizing 
all the results obtained in this paper.

\advance\thmno by 1
\proclaim\Thm. Let $C$ be a non-degenerate finitely semi-simple semi-strict 
spherical $2$-category of non-zero dimension. For any compact closed 
piece-wise linear oriented $4$-manifold $M$ there exists a state sum $I_C(M)$. For any two such manifolds $M$ and $N$ that are PL-homeomorphic we have 
$I_C(M)=I_C(N)$.\par

\thmno=0
\equano=0
\advance\chno by 1
\beginsection\the\chno. Example\par
Let $G$ be a finite group. In this section we explain how to obtain a 
finitely semi-simple non-degenerate semi-strict spherical 
$2$-category $\THilb[G]$ of non-zero dimension. The reason for using this 
notation is that if $\#G=1$ then this $2$-category 
will be just $\THilb_{cc}$ (see example 2.15). The invariant 
$I_{\THilb[G]}(M)$ looks very much like a 
$4$-dimensional version of the Dijkgraaf-Witten invariant [21].  

The objects of $\THilb[G]$ are the elements of the group rig ${\bf N}[G]$, 
which are 
just formal finite linear combinations of the elements of $G$ with 
non-negative integer coefficients. The elements of $G$ are taken to be 
the simple objects in $\THilb[G]$. So $\Hom(g,h)={\bf N}$ if $g=h$ and $\{0\}$ 
otherwise. Here we identify a non-negative integer $n$ with the $1\times 1$ 
matrix with coefficient $n$. 
This means that a $1$-morphism $F$ between $x=n_{i_1}g_{i_1}+\cdots 
+n_{i_k}g_{i_k}$ and $y=m_{j_1}g_{j_1}+\cdots +m_{j_l}g_{j_l}$ 
is a $(m_{j_1}+\cdots+m_{j_l})\times(n_{i_1}+\cdots+n_{i_k})$ 
matrix $(F_i^j)$ where the coefficient $F_i^j$ is 
just a non-negative integer. The composition of two 
$1$-morphisms $F$ and $G$ is defined by the matrix product of $F$ and 
$G$. Here we use the product and the sum of the non-negative integers as the 
operations on the coefficients of the matrices. A $2$-morphism $\alpha$ 
between 
$F\colon x\to y$ and $G\colon x\to y$ is a $(m_{j_1}+\cdots +m_{j_l})\times 
(n_{i_1}+\cdots +n_{i_k})$ matrix $(\alpha_i^j)$ 
where the coefficient $\alpha_i^j$ is a matrix with complex 
coefficients representing a linear map between $F_i^j$ and 
$G_i^j$. The vertical composition of two $2$-morphisms $\alpha$ 
and $\beta$ is defined by the matrix $(\alpha\cdot\beta)$ where the 
coefficient ${\alpha\cdot\beta}_i^j$ is the matrix product of 
${\alpha}_i^j$ and ${\beta}_i^j$. The horizontal composition 
of $\alpha$ and $\beta$ is defined by the matrix product of ${\alpha}$ 
and ${\beta}$, where the matrix product and the matrix sum are the 
operations on the coefficients. 

The tensor product of two objects $x$ and $y$ is just their product in 
${\bf N}[G]$. The tensor products on the $1$-morphisms and the 
$2$-morphisms are 
defined by the tensor products of the respective matrices, just as in $\THilb_{\rm cc}.$ The tensor 
product is not semi-strict in general. Between $g(hk)$ and $(gh)k$ we 
define the {\sl associator} to be the 1-isomorphism represented by 
the $1$-dimensional matrix with coefficient $1$. This means that we 
consider $g(hk)$ and $(gh)k$ only to be isomorphic and not identical! The 
pentagon in the following diagram is not commutative, but holds only up to a 
$2$-isomorphism 
$\pi(g,h,k,l)\in{\bf C}$, which is called the 
{\sl pentagonator}.
 
$$\vbox{\offinterlineskip
\halign{\hfil#\hfil&\quad\hfil#\hfil\quad&\hfil#\hfil\cr
$g(h(kl))$&$\rightarrow\quad (gh)(kl)\quad\rightarrow$
&$((gh)k)l$\cr
\noalign{\vskip7pt}
$\Bigg\downarrow$&&$\Bigg\uparrow$\cr
\noalign{\vskip7pt}
$g((hk)l)$&${\hbox to 50pt{\rightarrowfill}}$&$(g(hk))l$\cr}}$$ 

When ranging over all quadruples of elements in $G$ 
this is equivalent to a 
${\bf N}$-linear function $\pi\colon G\times G\times G\times G\to {\bf C}.$ 
The next order coherence relation that this pentagonator has to satisfy, 
which can be found in [27] for example, is exactly the condition that 
$\pi$ be a $4$-cocycle on the group. We do not write down explicitly this 
coherence condition for a general monoidal $2$-category because below we 
give a direct definition of the state-sum $I_{\THilb[G]}(M)$ and one can 
verify by hand 
that invariance under the $3\rightleftharpoons 3$ Pachner move is equivalent 
to the cocycle condition. In order to apply our construction in this 
particular example one has to strictify this tensor product in order to 
obtain a semi-strict monoidal $2$-category. But as we have mentioned before 
this can always be done, see [26], and we will not work this out in detail 
here. Note that the strictification does not eliminate the pentagonator but 
``hides'' it somehow inside the monoidal structure of the $2$-functor which 
defines the equivalence between the weak monoidal $2$-category and 
the strictified monoidal $2$-category. The coherence condition satisfied by the pentagonator, 
which in this particular example is the cocycle condition, allows us to 
``forget'' the pentagonator while we think abstractly about the construction 
of the state-sum. But once one writes down an explicit example one has to 
unpack the abstract definition of the state-sum and calculate where the 
pentagonator shows up. Of course there are different ways of unpacking, 
according to the different ways of reparenthesizing the tensor product, 
but the coherence condition guarantees that they all give the same answer in 
the end. Below we show how this works out in this example. 
Of course the same happens when 
one considers the Dijkgraaf-Witten invariant as an example of the general 
construction by Barrett and Westbury in [9].     

The dual of an object $x=n_{i_1}g_{i_1}+\cdots+n_{i_k}g_{i_k}$ is defined 
by $x^*=n_{i_1}g_{i_1}^{-1}+\cdots+n_{i_k}g_{i_k}^{-1}$. The coevaluation 
$i_g\colon 1\to gg^{-1}$ on a simple object $g$ is just the 1-dimensional 
matrix with coefficient $1$. The evaluation $e_g\colon g^{-1}g\to 1$ 
is also the 1-dimensional matrix with coefficient $1$. For arbitrary 
objects we obtain the coevaluation and evaluation by extending these 
definitions ${\bf N}$-linearly. The duality on 
$1$-morphisms and $2$-morphisms is defined as in $\THilb_{\rm cc}$ (see example 2.15).

Just as in $\THilb_{\rm cc}$ it is easy to show that this defines the right kind of 
duality and that $\THilb[G]$ is non-degenerate and spherical. $\THilb[G]$ is 
finitely semi-simple by construction and its dimension is equal to 
$|G|$. 

So we know that our construction gives an invariant $I_{\THilb[G]}(M)$, but of 
course it is easy to write down a more explicit formula for it in this 
case. Let $M=(M,T)$ be a triangulated manifold. Label the edges with 
elements of $G$ in such a way that $\ell(ij)\ell(jk)=\ell(ik)$ for any 
triangle $(ijk)$. If this rule is not satisfied for a certain labelling,  
then the Hom-category associated to this triangle is zero, so these 
labellings do not contribute anything to the state-sum. The $4$-simplex 
$(ijklm)$ gets the weight 
$\pi(\ell(ij),\ell(jk),\ell(kl),\ell(lm))$, where $\pi$ is the 
$4$-cocycle defining the pentagonator in $\THilb[G]$. Our state-sum is now 
equal to 
$$I_{\THilb[G]}(M)=|G|^{-v}\sum_{\ell}\prod_{i}\pi(S_i)^{\epsilon_i},$$
\noindent where $v$ is the number of vertices in $T$, the sum is taken over 
all labellings, in each term of the state-sum the product is taken over all 
$4$-simplices $S_i$ in $T$ with a fixed labelling, and 
$\epsilon_i$ is $+1$ if the orientation on $S_i$ induced by the ordering on 
its vertices and the one induced by the orientation of $M$ coincide, and 
$-1$ otherwise. Its easy to show that invariance under the 
$3\rightleftharpoons 3$ move is equivalent to the cocycle condition and that 
invariance under the $2\rightleftharpoons 4$ move also follows directly from 
that. Invariance under the $1\rightleftharpoons 5$ move follows from the 
identity
$$\pi(g,h,k,l)$$
$$=|G|^{-1}\sum_{m\in G}\pi(h,k,l,m)^{-1}\pi(gh,k,l,m)
\pi(g,hk,l,m)^{-1}\pi(g,h,kl,m)\pi(g,h,k,lm)^{-1}.$$ 
\noindent For this calculation one should reorder the vertices so that the 
``new'' vertex in the $1\rightleftharpoons 5$ move becomes the last one in 
the $5$-simplex. This is because in this example we have chosen 
a slightly different convention for the indices of our partition function 
for convenience. 
We know that the state-sum is independent of the ordering 
of the vertices so we are allowed to do this.               

If $M$ is connected and one takes the trivial $4$-cocycle, 
i.e. $\pi(g,h,k,l)=1$ for all 
$g,h,k,l\in G$, then this invariant just counts the number of 
group homomorphisms of $\pi_1(M)$ into $G$, as in the three dimensional 
case (see [21,33,34,38,39]). For a non-trivial cocycle it might become more 
interesting, although for $G={\bf Z}/p{\bf Z}$ the results so far have been 
disappointing. 
In [10] Birmingham and 
Rakowski study state-sums that appear as models for  
{\sl lattice gauge theories with finite gauge group} ${\bf Z}/p{\bf Z}$. 
In dimension 4 they study the state-sum that is similar to $I_{{\THilb}[{\bf Z}/p
{\bf Z}]}(M)$ except that the $4$-cocycle takes   
values in $U(1)$ and ${\bf Z}/p{\bf Z}$ respectively. Unfortunately there are no non-trivial $4$-cocycles of ${\bf Z}/p{\bf Z}$ 
with values in $U(1)$ and for any $4$-cocycle $\pi$ with values in ${\bf Z}/p{\bf Z}$ 
we get $\prod_{i}\pi(S_i)^{\epsilon_i}=1$ for a fixed labelling, so the 
invariant is not very interesting either (see section 5 in [10]). For other 
groups 
I do not know any results about $I_{\THilb[G]}(M)$. In [23] Freed and Quinn provided a detailed analysis of the 
topological quantum field theories (TQFT's) suggested by Dijkgraaf and Witten 
in [21], 
of which the restriction to closed manifolds is the Dijkgraaf-Witten 
invariant. Along the same line one could try to analyse $I_{\THilb[G]}(M)$. 
I must thank John Baez for pointing this out to me.  

Although in some cases the invariant defined by $I_{\THilb[G]}(M)$ might not 
be very sophisticated, it defines a 
non-trivial invariant in general and 
one could interpret $\THilb[G]$ with a 
non-trivial pentagonator as a finite deformation of $\THilb[G]$ with the 
trivial 
pentagonator. This is remarkable because we know that all the 
interesting 3-manifold state-sum invariants, including the so called 
{\sl quantum invariants}, come from 
deformations of monoidal categories. Maybe other interesting (quantum) 
4-manifold 
state-sum invariants will come from deformations of monoidal 
$2$-categories. Does categorification strike again?

As one can notice from the remarks above the example $I_{\THilb[G]}(M)$ is not 
new (see [10,33,34,38,39]). However it is interesting to see that 
this 
invariant is a nice example of the general construction 
presented in this 
paper, just as the Dijkgraaf-Witten invariant is a nice example of the 
general construction presented in [9].

\advance\chno by 1
\thmno=0
\equano=0
\beginsection\the\chno. Concluding remarks\par
In this last section I want to sketch my plans for further research 
based on some remarks concerning the results in this paper. 

First of all let me address the question of examples. Although I have 
given one example of a $2$-category that satisfies my conditions (section 6) 
and suggested that there is probably another one related to the 
categorification of the quantum double of a finite group, we need to 
find more interesting, and more complicated, examples. It is 
extremely hard to find interesting examples of Hopf categories or 
2-categories that are likely to give a new kind of invariant. Any 
``simple'' attempt seems to be doomed to lead to the Euler 
characteristic or, if one is a bit luckier, the signature of the 
$4$-manifolds (see [17,34]), or some other homotopy invariant. 
Crane, Yetter and I are working on the definition 
of monoidal $2$-categories built from representations of quantum 
groups at $q=0$ using their crystal bases and their deformation 
theory. These results will be published in a separate paper. It might 
well be that, in order to get the desired examples of $2$-categories, 
we have to find actual deformations of these $2$-categories and then 
find our way back to a generic $q$. Although this is the hard way of 
trying to find examples it certainly is worth a try.  

Next there is the question of the relation of the construction presented 
in this paper and the ones in [18] and [16]. In order to get back the 
Crane-Yetter invariant out of a construction like the one in this 
paper we would have to relax the definition of a simple object in a 
$2$-category somehow. The way we define things in this paper implies that 
a Vect-linear 
semi-simple monoidal $2$-category with one simple object is just Vect. 
The surprising definition of the dimension of a non-degenerate finitely 
semi-simple semi-strict $2$-category is probably a consequence of this 
restriction, although I do not know how exactly. One could try to take 
$\End(I)$ to be the free Vect-module category generated by an arbitrary 
finitely semi-simple monoidal category $R$, assume 
that each 
$\Hom$-category is a finite dimensional $R$-module category equivalent 
to $R^n$ for some non-negative integer $n$ and 
define a simple object to be an object $A$ such that $\End(A)\simeq R$ 
(here we denote the free Vect-module category generated by $R$ also by 
$R$). 
With such a relaxation 
we would get the Crane-Yetter invariant via our construction for a 
non-degenerate finitely semi-simple semi-strict spherical $2$-category with 
one object (see lemma 2.14). However not all the definitions in 
this paper would still be adequate; for example, one would have to take the 
tensor product over $R$ rather than Vect in definition 2.20. Furthermore 
I do not 
know if all the proofs in this paper 
would still hold after 
such a relaxation. At this point it is not even clear to me that such a relaxation is going 
to be 
relevant for our search for quantum $4$-manifold invariants. In the three 
dimensional case one gets all the known quantum invariants out of 
constructions 
involving monoidal categories where $\End(I)$ is a field 
and not just a commutative ring. A field gives only the trivial 
invariant of $2$-dimensional manifolds when one uses the construction of 
a state-sum out of a finite dimensional semi-simple algebra as defined in [25]. Their state-sum is equal to 1 for every manifold if the semi-simple algebra 
is a field. In a way the Barrett-Westbury invariant [9] is a lift to the 
third dimension of the 
invariant defined in [25]. So it is not clear yet whether, in the process of 
categorification, one needs the 
most general structure for $\End(I)$ on the next level or one can assume that 
it 
is the simplest example possible, a one dimensional $n$-Vector space on level 
$n$, in order to find ``new'' $n+2$-manifold invariants.

Now let me say something about the possible relation between the 
state-sum invariant in this paper and the {\sl Tornado Formula} in [16]. 
First of all let me recall that in the construction of $4$-manifold invariants 
out of Hopf categories there is this difficult question of the right categorified notion of an antipode. Crane and Frenkel do not address this question in 
[16]. Neither 
do they prove invariance of their state-sum under the permutation of vertices. In [31] Neuchl gives a definition of an antipode in a Hopf 
category. It 
is not clear though, at least to me, that the {\sl Tornado Formula} is 
invariant 
under the permutation of vertices, even if one uses Neuchl's definition.    
As I already mentioned in my introduction, the only 
example of a Hopf category of which we are sure that it gives an 
invariant is the categorification of the quantum double of a 
finite group [11,19]. In that case one has to impose some cyclical conditions 
on the cocycles that define the structure isomorphisms. In order to 
resolve the problem of invariance under the permutation of vertices in 
general and 
to establish a concrete relation 
between the construction in this paper and the one in [16], I think 
one should first try to reconstruct the Hopf category out of a spherical 
$2$-category. By Neuchl's result [31] we know that it is going to 
be a bitensor category. The duality used in the construction in this 
paper, which is the right one for the purpose of $4$-manifold invariants, 
will lead to the reconstruction of 
the right notion of the antipode, hopefully the same as Neuchl's, and 
the proof of invariance of the Crane-Frenkel state-sum. 
When we really understand the Crane-Frenkel state-sum, we can try to see 
if the invariant coming from 
the Hopf category is the same as the one coming from the $2$-category 
of its representations. For involutory Hopf algebras and their 
category of representations Barrett and Westbury show in [8] that the 
respective $3$-manifold invariants are the same. 
This result was my original motivation for the ``categorification'' 
of Barrett's and Westbury's construction of $3$-manifold invariants. 

Finally I want to say something about diagrams. If we ever want to 
compute a real quantum invariant of $4$-manifolds we had probably better 
start looking for some diagrammatic way of doing this. Every reader 
will remember the major advantage that the ``skein approach'' 
brought to the computation of $3$-manifold invariants. It even 
enabled Roberts to show that the Crane-Yetter $4$-manifold 
invariant for 
$U_q(sl(2))$ was ``just'' the signature of the $4$-manifolds [36]. 
In this paper I have used some diagrams, but, as I already admitted, 
they are rather poor diagrams and certainly not good enough for 
computational aims. 
For both the construction and the computation of $3$-manifold 
invariants it is very convenient to have the correspondence 
between pivotal categories and labelled oriented planar graphs as shown 
in [24]. Let us formulate this correspondence precisely. A $C$-labelled 
graph is an oriented planar graph with 
its edges labelled 
by objects of $C$ and its vertices labelled by morphisms of $C$, 
such that the source of such a morphism is the tensor product of 
the labels of the ingoing edges and the target the tensor product 
of the labels of the outgoing edges. 
If a graph contains two edges with the same vertices labelled by the objects 
$X$ 
and $Y$ respectively, then we can substitute them by one single edge 
labelled by $X\otimes Y$, and vice versa. If a graph contains an edge 
labelled by $X$, then we can substitute this edge by the edge with 
opposite orientation labelled by $X^*$. We call two $C$-labelled  
graphs {\sl 
essentially equivalent} if one can be obtained from the other by a finite 
number of such substitutions. We now define two $C$-labelled graphs 
$G_1$ and $G_2$ to be equivalent if and only if there are two $C$-labelled 
graphs $G_1'$ and $G_2'$ such that $G_i'$ is isotopic to $G_i$ 
for $i=1,2$ and $G_1'$ is essentially equivalent to $G_2'$. The main 
theorem of Freyd and Yetter in [24] can now be formulated as follows: 
the equivalence classes of $C$-labelled graphs form a pivotal 
category equivalent to $C$. This result enables one to translate algebraic 
manipulations with morphisms in $C$ into diagrammatic moves. As a 
matter of fact Barrett and Westbury [9] do assume, 
without using it explicitly, a similar result for spherical 
categories 
and labelled oriented graphs in $S^2$. 
The question arises what kind of graphs do correspond to spherical 
$2$-categories. If one dualizes 
a tetrahedron and looks at the 2-skeleton of the dual complex, 
one gets a four-valent vertex of a graph with $2$-cells. 
So maybe the answer to my question is:  
``labelled oriented graphs with $2$-cells in $S^3$''. 
Here I am deliberately being 
sloppy in my definition of the sort of graph I am speaking about. Probably 
one has to define it by something like a set of elements, the vertices, 
and a family of two-element subsets, the edges, and a family of 
three- or more-element subsets, the $2$-cells. But I am not trying to 
make a precise conjecture here. I just want to point out a possible topic 
for further research that can lead to a better insight in the relation 
between $4$-dimensional topology and combinatorics and algebra. 
In order to get such a result we could try to study this kind of graph, or hyper-graph, in a  
way similar to that in which Carter, Rieger and Saito [13,14] have studied 
$2$-tangles, and see if they provide the diagrammatical tools for 
the computation of $4$-manifold invariants.     

\beginsection Acknowledgements\par
First of all I want to thank Louis Crane and David Yetter for the 
helpful discussions with them. I also 
want to thank the Mathematics Department of Kansas State University (KSU) 
and Louis Crane in particular for their hospitality during my stay 
in Manhattan, where I wrote this paper. Finally I want to thank 
the Centro de Matematica Aplicada (CMA) and the Funda\c c\~ao 
Luso-Americana para o Desenvolvimento (FLAD) in Portugal for their 
generous financial support which enabled me to visit KSU.

\vglue 12pt
\line{\bf References\hfil}
\vglue 5pt
\frenchspacing

\advance\refno by 1
\item{[\the\refno]} J.~Baez, J.~Dolan, Higher-Dimensional Algebra and 
Topological Quantum Field Theory, {\it Jour. Math. Phys.} {\bf 36} (1995), 6073-6105.
\advance\refno by 1
\item{[\the\refno]} J.~Baez, M. Neuchl, Higher-Dimensional Algebra I: 
Braided Monoidal $2$-Categories, {\it Adv. Math.} {\bf 121} (1996), 196-244. 
\advance\refno by 1
\item{[\the\refno]} J.~Baez, Higher-Dimensional Algebra II: 2-Hilbert 
Spaces, {\it Adv. Math.} {\bf 127} (1997), 125-189.
\advance\refno by 1
\item{[\the\refno]} J.~Baez, J.~Dolan, Higher-Dimensional Algebra III: 
$n$-Categories and the Algebra of Opetopes, {\it Adv. Math.} {\bf 135} (1998), 
145-206.
\advance\refno by 1
\item{[\the\refno]} J.~Baez, L.~Langford, Higher-Dimensional 
Algebra IV: $2$-Tangles, to appear in {\it Adv. Math.}
\advance\refno by 1
\item{[\the\refno]} J.~Baez, L.~Langford, $2$-Tangles, 
{\it Lett. Math. Phys.} {\bf 43} (1998), 187-197.
\advance\refno by 1
\item{[\the\refno]} J.W.~Barrett, B.W.~Westbury, Spherical 
Categories, to appear in {\it Adv. Math.} Preprint available as hep-th/9310164.
\advance\refno by 1
\item{[\the\refno]} J.W.~Barrett, B.W.~Westbury, The Equality 
of $3$-Manifold Invariants, {\it Math. Proc. Camb. Phil. Soc.} {\bf 
118} (1995), 503-510.
\advance\refno by 1
\item{[\the\refno]} J.W.~Barrett, B.W.~Westbury, Invariants of 
Piecewise-Linear $3$-Manifolds, {\it Trans. Amer. Math. Soc.} {\bf 348 (10)} 
(1996), 3997-4022.
\advance\refno by 1
\item{[\the\refno]} D.~Birmingham, M.~Rakowski, On Dijkgraaf-Witten 
Type Invariants, {\it Lett. Math. Phys.} {\bf 37} (1996), 363.
\advance\refno by 1
\item{[\the\refno]} J.S.~Carter, L.H.~Kauffman, M.~Saito, Diagrammatics, 
Singularities, and Their Algebraic Interpretations, 
to appear in {\it Matematica Contemporanea, Sociedade Brasileira de 
Matematica.} Preprint available at http://www.math.usf.edu/$\sim$saito/home.html.
\advance\refno by 1
\item{[\the\refno]} J.S.~Carter, L.H.~Kauffman, M.~Saito, 
Structures and Diagrammatics of Four Dimensional Topological Lattice Field 
Theories, preprint available as GT/9806023.
\advance\refno by 1
\item{[\the\refno]} J.S.~Carter, J.H.~Rieger, M.~Saito, A Combinatorial 
Description of Knotted Surfaces and Their Isotopies, {\it Adv. Math.} {\bf 127
(1)} (1997), 1-51.
\advance\refno by 1
\item{[\the\refno]} J.~Carter, M.~Saito, Knotted Surfaces, Braid Movies, 
and Beyond, in {\it Knots and Quantum Gravity}, ed. J.~Baez, Oxford U. Press 
(1994), 191-229.
\advance\refno by 1
\item{[\the\refno]} S.~Chung, M.~Fukuma, A.~Shapere, The Structure 
of Topological Field Theories in Three Dimensions, {\it Int. J. Mod. Phys.} 
{\bf A9} (1994), 1305-1360.
\advance\refno by 1
\item{[\the\refno]} L.~Crane, I.~Frenkel, Four Dimensional Topological 
Quantum Field Theory, Hopf Categories, and the Canonical Bases, 
{\it Jour. Math. Phys.} {\bf 35} (1994), 5136-5154.
\advance\refno by 1
\item{[\the\refno]} L.~Crane, L.H.~Kauffman, D.N.~Yetter, State-Sum 
Invariants of $4$-Manifolds, {\it JKTR} {\bf 6 (2)} (1997), 177-234.
\advance\refno by 1
\item{[\the\refno]} L.~Crane, D.N.~Yetter, A Categorical Construction 
of 4D  Topological Quantum Field Theories, in {\it Quantum Topology}, 
L.H.~Kauffman and R.A.~Baadhio eds., World Scientific Press (1993), 131-138.
\advance\refno by 1
\item{[\the\refno]} L.~Crane, D.N.~Yetter, Examples of Categorification, 
to appear in {\it Cahiers de Top. Geom. Diff. Cat.}
\advance\refno by 1
\item{[\the\refno]} L.~Crane, D.N.~Yetter, Deformations of (Bi)tensor 
Categories, to appear in {\it Cahiers de Top. Geom. Diff. Cat.}
\advance\refno by 1
\item{[\the\refno]} R.~Dijkgraaf, E.~Witten, Topological Gauge Theories 
and Group Cohomology, {\it Comm. Math. Phys.} {\bf 129 (2)} (1990), 393-429.
\advance\refno by 1
\item{[\the\refno]} J.~Fisher, $2$-Categories and $2$-Knots, {\it Duke Math. 
Jour.} {\bf 75} (1994), 493-526.
\advance\refno by 1
\item{[\the\refno]} D.~Freed, F.~Quinn, Chern-Simons Theory with Finite 
Gauge Group, {\it Comm. Math. Phys.} {\bf 156} (1993), 435-472. 
\advance\refno by 1
\item{[\the\refno]} P.J.~Freyd, D.N.~Yetter, Braided Compact 
Closed Categories with Applications to Low Dimensional Topology, 
{\it Adv. Math.} {\bf 77} (1989), 156-182.
\advance\refno by 1
\item{[\the\refno]} M.~Fukuma, S.~Hosono, H.~Kawai, Lattice Topological 
Field Theory in two dimensions, {\it Comm. in Math. Phys.} {\bf 161} (1994), 
157-176.
\advance\refno by 1
\item{[\the\refno]} R.~Gordon, A.~Powers, R.~Street, Coherence for 
Tricategories, {\it Memoirs of the AMS} {\bf 117 (558)} (1995), 1-81.
\advance\refno by 1
\item{[\the\refno]} M.~Kapranov, V.~Voevodsky, $2$-Categories and 
Zamolodchikov Tetrahedra Equations, in {\it Proc. Symp. Pure Math}. {\bf 56} 
Part 2 (1994), AMS, Providence, 177-260.
\advance\refno by 1
\item{[\the\refno]} V.~Kharlamov, V.~Turaev, On the Definition of the 
$2$-category of $2$-Knots, {\it Amer. Math. Soc. Transl.} {\bf 174} (1996), 
205-221.
\advance\refno by 1
\item{[\the\refno]} G.~Kuperberg, Involutory Hopf Algebras and 
$3$-Manifold Invariants, {\it Int. J. of Math.} {\bf 2 (1)} (1991), 41-66.
\advance\refno by 1
\item{[\the\refno]} L.~Langford, {\it $2$-Tangles as a Free Braided 
Monoidal $2$-category with Duals}, Ph.D. Dissertation, UCR, 1997. Available at http://math.ucr.edu/home/baez/.
\advance\refno by 1
\item{[\the\refno]} M.~Neuchl, {\it Representation Theory  of Hopf 
Categories}, Ph.D. Dissertation, University of Munich, 1997.
\advance\refno by 1
\item{[\the\refno]} U.~Pachner, P.L. Homeomorphic Manifolds are 
Equivalent by Elementary Shellings, {\it European Journal of Combinatorics} 
{\bf 12} (1991), 129-145.
\advance\refno by 1
\item{[\the\refno]} T.~Porter, Interpretations of Yetter's Notion of 
$G$-coloring: Simplicial Fibre Bundles and Non-Abelian Cohomology, 
{\it JKTR} {\bf 5} (1996), 687-720.
\advance\refno by 1
\item{[\the\refno]} T.~Porter, TQFT's from Homotopy $n$-Types, 
to appear in {\it J. London Math. Soc.} Preprint available at: http://www.bangor.
ac.uk/$\sim$mas013/publist.html.
\advance\refno by 1
\item{[\the\refno]} N.Yu.~Reshetikhin, V.G.~Turaev, Ribbon Graphs and 
their Invariants derived from Quantum Groups, {\it Comm. Math. Phys.} {\bf 127 (1)} (1990), 1-26. 
\advance\refno by 1
\item{[\the\refno]} J.~Roberts, Skein Theory and Turaev-Viro 
Invariants, {\it Topology} {\bf 34} (1995), 771-787.
\advance\refno by 1
\item{[\the\refno]} V.G.~Turaev, O.Y.~Viro, State Sum Invariants 
of $3$-Manifolds and Quantum $6j$-Symbols, {\it Topology} {\bf 31 (4)} (1992), 
865-902.
\advance\refno by 1
\item{[\the\refno]} D.N.~Yetter, Topological Quantum Field Theories Associated to Finite Groups and Crossed $G$-Sets, {\it JKTR} {\bf 1 (1)} (1992), 1-20.
\advance\refno by 1
\item{[\the\refno]} D.N.~Yetter, TQFT's From Homotopy $2$-Types, {\it JKTR} 
{\bf 2 (1)} (1993), 113-123.
\advance\refno by 1
\item{[\the\refno]} D.N.~Yetter, State-Sum Invariants of 
$3$-Manifolds Associated to Artinian Semisimple Tortile Categories, 
{\it Topology and Its Applications} {\bf 58 (1)} (1994), 47-80.

\bye